\documentclass[12pt]{article}

\usepackage{microtype}
\usepackage{graphicx} 
\usepackage[a4paper, margin=1in]{geometry}

\newsavebox{\overlongequation}
\newenvironment{overlong}
 {\begin{displaymath}\begin{lrbox}{\overlongequation}$\displaystyle}
 {$\end{lrbox}\makebox[0pt]{\usebox{\overlongequation}}\end{displaymath}}  

\usepackage{relsize}
\usepackage{numprint}
\npthousandsep{\,}
\newcommand{\smallnumprint}[1]{\text{\smaller \numprint{#1}}}

\usepackage{amsmath,amsthm,amsfonts,amssymb}

\usepackage{subcaption}

\usepackage{shuffle}

\usepackage{MnSymbol}
\usepackage{enumitem}
\setlist[enumerate,1]{label={\rm(\arabic*)}}

\usepackage{xcolor}
\usepackage{mathtools}

\definecolor{mycolor}{RGB}{137,57,99}

\usepackage[colorlinks = true,
            linkcolor =  mycolor,
            urlcolor  = mycolor,
            citecolor = mycolor,
            anchorcolor = mycolor,
            ]{hyperref}
\hypersetup{
    pdftitle = {Toric Varieties from Multiple Zeta Values},
    pdfauthor = {Annika Burmester, Steven Charlton, Abhiram Kidambi, Felix Lotter, Rosa Preiß},
    pdfsubject = {The linear relations among fixed weight multiple zeta values (MZVs) are a well-studied subject in number theory. Elements of putative linear bases for this fixed-weight space are, however, not algebraically independent; they satisfy non-linear relations. We study these relations by approximating them through motivic MZVs, formal MZVs and numerical computations. We show that motivic MZVs in fixed weight span a projective toric variety and explore its geometry. We define an outer approximation to the variety of formal MZVs in fixed weight by introducing toric varieties coming from shuffle and stuffle relations and taking their intersection. We conjecture that this intersection in fact agrees with the variety of formal MZVs and hence conjecturally also with the variety of motivic MZVs. This leads to a novel conjecture on a generating set for the extended double shuffle relations, which we checked computationally up to weight 23.}
}
\usepackage{tikz,calc}
\usetikzlibrary{positioning}
\usepackage{tikz-cd}

\usepackage{pgfplots}
\usetikzlibrary{arrows.meta,decorations.markings}

\tikzset{
  directed/.style={
    postaction={decorate},
    decoration={
      markings,
      mark=at position #1 with {
        \arrow{Stealth[length=2.5mm,width=2mm]}
      }
    }
  },
  directed/.default=0.5
}
\pgfplotsset{compat=1.18}
\tikzset{oriented/.style={postaction={decorate}}}

\newcommand{\noaddtocontents}[2]{}
\newcommand\stoptoc{\let\origaddtocontents\addtocontents
    \let\addtocontents\noaddtocontents
 }
\newcommand\resumetoc{\let\addtocontents\origaddtocontents
 }

\newcommand{\gp}{\texttt{Pari/GP}}
\newcommand{\sage}{\texttt{SageMath}}
\newcommand{\macaulay}{\texttt{Macaulay2}}
\newcommand{\github}{\texttt{GitHub}}

\newcommand{\R}{\mathbb{R}}
\newcommand{\Q}{\mathbb{Q}}
\newcommand{\IIS}{\mathsf{IIS}}
\newcommand{\ISS}{\mathsf{ISS}}

\newcommand{\intd}{\,\mathrm d}

\newcommand{\letter}[1]{{\color{cyan}\mathtt{#1}}}
\newcommand{\stuffle}{{\mathbin{*}}}
\newcommand{\tensor}{\otimes}

\newcommand{\mzvb}[1]{\zeta({\textcolor{black}{#1}})}

\newcommand{\ZmC}{\mathcal{Z}^{\mathfrak m}}
\newcommand{\lyndon}{\mathcal{L}}
\newcommand{\zetashreg}{\zeta^{\shuffle,\mathrm{r}}}
\newcommand{\zetastreg}{\zeta^{\stuffle,\mathrm{r}}}
\newcommand{\algshreg}{S^{\shuffle,\mathrm r}}
\newcommand{\algstreg}{S^{\stuffle,\mathrm r}}
\newcommand{\monlyndon}{\mathcal{M}}
\newcommand{\lyndonmone}{\lyndon_{\not=\letter 1}}
\newcommand{\monlyndonmone}{\mathcal{M}_{\not=\letter 1}}
\newcommand{\weight}{\operatorname{wt}}

\newcommand{\true}{$\checkmark$}
\usepackage{utfsym}
\newcommand{\newcrossmark}{\scalebox{0.75}{\usym{2613}}}

\newcommand{\false}{\newcrossmark}
\newcommand{\relvar}{\mathcal R}
\newcommand{\relvarr}{\mathcal R^0}
\newcommand{\relalg}{R}

\newcommand{\Z}{\mathcal Z}

\newcommand{\motZ}{\Z^{\mathfrak m}}
\newcommand{\formZ}{\Z^{\mathfrak f}}

\newcommand{\StC}{\mathcal S^\stuffle}
\newcommand{\StCmone}{\mathcal S^{\stuffle,\mathrm{r}}}
\newcommand{\wt}[1]{\mathrm{wt}(#1)}
\newcommand{\dep}[1]{\mathrm{dep}(#1)}

\newtheorem{theorem}{Theorem}[section]
\newtheorem{conjecture}[theorem]{Conjecture}

\newtheorem{proposition}[theorem]{Proposition}
\newtheorem{corollary}[theorem]{Corollary}
\newtheorem{lemma}[theorem]{Lemma}

\theoremstyle{definition}
\newtheorem{defi}[theorem]{Definition}
\newtheorem{remark}[theorem]{Remark}
\newtheorem{example}[theorem]{Example}
\newtheorem*{notdef}{Notation/Definition}

\newtheorem*{remark*}{Remark}

\DeclareMathOperator{\MT}{MT}
\DeclareMathOperator{\Isom}{Isom}
\DeclareMathOperator{\comp}{comp}
\DeclareMathOperator{\Sym}{Sym}

\newcommand{\mot}{\mathfrak{m}}

\title{Toric Varieties from Multiple Zeta Values}
\author{Annika Burmester, Steven Charlton, Abhiram Kidambi, \\ Felix Lotter, Rosa Prei\ss }
\date{}

\begin{document}

\maketitle

\begin{abstract}
The linear relations among fixed weight multiple zeta values (MZVs) are a well-studied subject in number theory. Elements of putative linear bases for this fixed-weight space are, however, not algebraically independent; they satisfy non-linear relations. We study these relations by approximating them through motivic MZVs, formal MZVs and numerical computations. We show that motivic MZVs in fixed weight span a projective toric variety and explore its geometry. We define an outer approximation to the variety of formal MZVs in fixed weight by introducing toric varieties coming from shuffle and stuffle relations and taking their intersection. We conjecture that this intersection in fact agrees with the variety of formal MZVs and hence conjecturally also with the variety of motivic MZVs. This leads to a novel conjecture on a generating set for the extended double shuffle relations, which we checked computationally up to weight $23$.
\end{abstract}

	{
		\setcounter{tocdepth}{2}
		\tableofcontents
		\setcounter{tocdepth}{3}
	}

\section{Introduction}
An important question in arithmetic geometry and number theory pertains to the irrationality and transcendentality of zeta values \( \zeta(k) := \displaystyle \sum\nolimits_{n=1}^\infty n^{-k} \), for positive integers $k>1$.
Using Euler's result that the even zeta values \( \zeta(2k) \) are rational multiples of \( \pi^{2k} \), for integers \( k \geq 1 \), it follows from the Lindemann--Weierstra\ss{} theorem that the even zeta values are transcendental.
However, for odd zeta values even irrationality is poorly understood -- it is known that there are an infinite number of irrational odd zeta values \cite{ballrioval} including $\zeta(3)$ \cite{apery}. Progress in explicitly proving irrationality of other odd zeta values has been challenging; an important result in this direction is due to Zudilin \cite{Zudilin2001} who proved that at least one of $\zeta(5), \zeta(7), \zeta(9), \zeta(11)$ is irrational.\footnote{There are a growing number of AI assisted manuscripts on irrationality of several odd zeta values which we avoid discussing here.} There is no known statement on transcendentality of odd zeta values. 
 However, the widely believed folklore is that odd zeta values are transcendental and  algebraically independent. 
\begin{conjecture} [Transcendence conjecture] \label{conj:transendence}
The numbers $\pi, \zeta(3),\zeta(5), \zeta(7),\ldots$ are algebraically independent over $\mathbb{Q}$.
\end{conjecture}
To study questions pertaining to algebraic relations between zeta values, it helps to think of zeta values as a part of a larger algebra of \emph{multiple zeta values} (MZVs).  These are generalizations of zeta values given by the series
\begin{equation}\label{eq:sum-def-mzv}
\zeta(k_1,\ldots,k_r)=\sum_{0<n_1<\cdots<n_r} \frac{1}{n_1^{k_1}\cdots n_r^{k_r}} \in \mathbb{R}, \ k_1,\ldots, k_r\in \mathbb{Z}_{>0},\ k_r\in \mathbb{Z}_{\geq 2}.
\end{equation}
We will demand that $k_r\geq2$ for convergence, and a multiset $\underbar{k} = (k_1,\ldots, k_r)\in \mathbb{Z}_{>0},\ k_r\in \mathbb{Z}_{\geq 2}$ is said to be \emph{admissible}.

\begin{remark*}[Notation of MZVs]
    In the literature, there is also the opposite convention of taking the sum \eqref{eq:sum-def-mzv} over $n_1 > \dots > n_r > 0$. In this article we choose the increasing order. The two definitions are of course equivalent under reversal of ${k_1}, \ldots, {k_r}$.
\end{remark*}

    Even and odd zeta values (together known as \emph{single} zeta values) are examples of MZVs with $r=1$. MZVs have a richer algebraic structure and satisfy many linear and algebraic relations. 
\begin{defi}
    Given a MZV $\zeta(\underbar{k})$, $\underbar{k} = (k_1,k_2,\ldots, k_r)$ admissible, we define the weight of the MZV/multiset as $\displaystyle \wt{\underbar{k}} := \sum\nolimits_{i=1}^r k_i$ and its length/depth as $\dep{\underbar{k}} := r$.
\end{defi}
Let us define $\Z$ (resp. $\Z_j$) as the $\mathbb Q$-vector space spanned by all admissible MZVs (resp.\ those in weight $j$):
\begin{align}
    \Z := \left\langle \zeta(\underbar{k})\mid \underbar{k} \text{ admissible}\right\rangle, \qquad     \Z_j := \left\langle \zeta(\underbar{k})\mid \underbar{k} \text{ admissible}, \ \wt{\underbar{k}} = j \right\rangle.
\end{align}

It follows from the definition of MZVs that $\Z$ is in fact a $\mathbb Q$-subalgebra of $\mathbb R$, also see Proposition \ref{prop:zetastufflemorphism}.

\begin{conjecture}[Weight grading conjecture \cite{Hof97}]
\label{conj:weightgrading}
$\displaystyle \Z = \bigoplus\nolimits_j \Z_j$, i.e., the algebra of MZVs is weight graded. 
    \end{conjecture}
    
    A proof of the above conjecture would in particular imply the transcendence of all MZVs over $\mathbb Q$. The question of algebraic independence of zeta values is closely linked to the linear spaces $\Z_j$ and in particular their dimension.
    
\begin{conjecture}[Zagier]
\label{conj:dimQZ}
    Let $\displaystyle \frac{1}{1-t^2 - t^3} = \sum_{m=0}^\infty d_m t^m$. Then, $\dim_{\mathbb Q} \Z_j=d_j$.
\end{conjecture}
It is known from the work of Brown \cite{brown2012mixed}, Deligne--Goncharov \cite{delignegoncharov}, and Terasoma \cite{terasomamzv} that $\dim_{\mathbb Q} \Z_j \leq d_j$. While there is no proof of Conjecture \ref{conj:dimQZ}, there is formal \cite{Blumlein:2009cf} and numerical evidence for this. (See either \cite[Table 1.1]{mzvfromnumberstomotives}, or the GitHub repository linked in Section \ref{sec:code}.)

\subsection*{Overview and main results}
Our main question in this article  concerns the weight-homogeneous spaces $\mathcal Z_j$:

\begin{center}
    \textit{Which polynomial equations with rational coefficients are satisfied by weight $j$ MZVs?} 
\end{center}
\noindent In other words, we ask for the algebraic variety spanned by weight $j$ MZVs over the rational numbers. We approach this question by first setting up the conventions to study these \emph{relation varieties} for MZVs in Section \ref{sec:mzvvariety}. We explain in general how to approximate these varieties via free graded Hopf algebras which surject onto $\Z$. The key point here is that for any such free graded algebra, the associated relation varieties are toric. 
\\[0.2cm] \noindent
Our first approximation is the free graded Hopf algebra of \emph{motivic MZVs}. In Section \ref{sec:mot-var}, we study the associated relation varieties in detail.  
The main results in this section concern their geometric properties: Theorem \ref{thm:mot-dim} determines the dimension of the motivic MZV variety, Theorem \ref{thm:mot-non-normal} the normality properties, and Proposition \ref{prop:not_gen_in_fixed_deg} shows that there is no general bound on the degree of the minimal generators of its ideal. Table \ref{tab:mot} showcases key properties of these varieties up to weight $23$, obtained by symbolic computations in \macaulay.\\[0.2cm]
\noindent
In Section \ref{sec:conj}, we recall the shuffle and stuffle relations satisfied by MZVs, as well as the associated regularizations, and use them to define (regularized) shuffle and stuffle varieties. Our considerations motivate Conjecture \ref{conj:lin-form-approx-eq} which posits a small generating set for the ideal of the extended double shuffle relations, which seems to be novel. We have verified this conjecture computationally up to weight $23$. The conjecture would imply that the intersection of regularized shuffle and stuffle varieties yields the relation variety of formal MZV in the given weight.
\\[0.2cm] \noindent
Motivated by this, we study the regularized stuffle variety in Section \ref{sec:stuf-var}. The main result here is Theorem \ref{thm:unregstuffalg} which studies the geometry of such varieties, and Table \ref{tab:stuffle} showcases again central properties of the variety and its ideal up to weight $14$.
\\[0.2cm] \noindent
In Section \ref{sec:comp}, we compare the intersection of the regularized shuffle and stuffle variety to the double shuffle scheme $\operatorname{DM}$ introduced by Racinet \cite{Racinetthesis}. We explain the connection between MZVs and the theory of \emph{path signatures}. 
\\[0.2cm] \noindent
In Section \ref{sec:computation}, we provide extensive computational and numerical evidence in \gp{} and \sage{} to supplement our main results. We explain our computational workflow in detail. We analyse the existence of what we call \emph{exceptional relations} (which are relations not predicted by our toric formalism) and provide certifiable results (in a few degrees and weights) for the absence of such relations in accordance with the Grothendieck period conjecture.
Details for the \github{} repository hosting code and results is also made available here, including \texttt{Macaulay2} code for evaluating MZVs for the benefit of readers unfamiliar with \gp{} or \sage. The mathematical approach for this evaluation is outlined in Appendix \ref{app:macaulayzeta}.
\\[0.2cm]
\noindent Our results provide a novel viewpoint on the algebraic structures governing MZVs. For future work, it would be desirable to connect our findings to the theory of classical modular forms which play an important role in theory of MZVs (cf.\ \cite{gkz} and \cite[Chapter 5]{ZhaoMZVbook}).

\paragraph*{Acknowledgments:} We thank Dominic Bunnett, David Broadhurst, Herbert Gangl, Kris Shaw, Bernd Sturmfels, and Wadim Zudilin for useful discussions.  AK is funded by the European Union (ERC, MaScAmp, 101167287) and the Max Planck Society. AB is funded by Deutsche Forschungsgemeinschaft (DFG, German Research Foundation) – Project-ID 491392403 – TRR 358. FL and RP acknowledge funding by the Deutsche Forschungsgemeinschaft (DFG, German Research Foundation) – CRC/TRR 388 “Rough Analysis, Stochastic Dynamics and Related Fields” – Project A04, 516748464. AB and AK thank the Erwin Schr\"odinger Institute for its hospitality.  AK thanks ICTP Trieste, Kavli IPMU and the Isaac Newton Institute for Mathematical Sciences for hospitality.  

\paragraph*{LLM and AI usage declaration:} We used ChatGPT 5.6 to help speed up, optimize and parallelize code that was run on a scientific computing cluster, and to write scripts to parse data from and read the MZV data mine. The algorithm to certify absence of exceptional relations within a given coefficient bound (See Section \ref{sec:lll}) was proposed by ChatGPT 6.1 Sol, which superseded an older, restricted computational heuristic we had. We also used ChatGPT 6.1 Sol (in Codex) to help code the fixed-weight relation data mine. We used ChatGPT 5.6 to help generate the \texttt{TikZ} figures (which we then modified) since life is too short to learn \texttt{TikZ}. No other part of this manuscript was generated by AI.
This declaration is \textit{not} an endorsement of OpenAI.

\section{The variety of MZVs and its approximations}
\label{sec:mzvvariety}

From now on, we write $I_j$ for the kernel of the map
\[\Sym \Z_j \to \Z\]
induced by inclusion of $\mathbb Q$-vector spaces. The elements of this ideal are the non-linear relations among weight $j$ MZVs. Geometrically, $I_j$ cuts out an irreducible affine variety in the complex affine space $\Z_j^\vee \otimes_{\mathbb Q} \mathbb C$, with vanishing ideal $I_j \tensor_{\mathbb Q} \mathbb C$. 

\begin{example} \label{ex:I_3}
     There are only two admissible integer tuples in weight $3$, defining MZVs $\zeta(3)$ and $\zeta(1,2)$ which span $\mathcal Z_3$. It is known that $\zeta(3)$ is irrational \cite{apery, beukers1, poorten} and that $\zeta(3) = \zeta(1,2)$ (Euler's identity), so $\mathcal Z_3$ is a one-dimensional vector space. $I_3$ is the kernel of the map $\mathbb Q[x_{\zeta(3)}] \to \Z \subset \mathbb R$ that sends $x_{\zeta(3)}$ to $\zeta(3)$. Determining $I_3$ amounts to deciding transcendentality of $\zeta(3)$, which is included in Conjecture \ref{conj:transendence}.
\\[0.2cm]
    Now consider weight $j=10$. There are $256$ admissible integer tuples, but $\mathcal Z_{10}$ is at most a $7$-dimensional vector space. It is spanned for example by
\[\zeta(10),\zeta(1,9),\zeta(3,7),\zeta(5,5),\zeta(1,1,8),\zeta(1,3,7),\zeta(1,4,5).\]
These $7$ generators satisfy a quadratic relation which is contained in $I_{10}$, see Example \ref{ex:mot-weight-10-expansion}.
\end{example}

Ideally, we would like to determine the ideal $I_j$ and study the associated geometry. However, as the case $j=3$ above already illustrated, this problem is currently out of reach: the transcendence (and even irrationality) of most MZVs remains the subject of deep and long-standing conjectures. To avoid this problem, we study the question instead for closely related graded algebras that surject onto $\mathcal Z$, but are better understood. In this way, we obtain approximations of $I_j$:

\begin{lemma}\label{lem:inner-ideal-approx}
     If $H$ is a graded $\mathbb Q$-algebra and $\psi: H \to \Z$ is an algebra homomorphism that maps $H_j$ into $\Z_j$, then
     \[(\Sym \psi_j)(\ker \phi_{H,j}) \subseteq I_j,\]
     where $\phi_{H,j}$ is the natural map $\Sym H_j \to H$ and $\psi_j$ denotes the restriction of $\psi$ to $H_j$.
\end{lemma}
\begin{proof}
    This follows immediately from commutativity of the diagram
    \[ \begin{tikzcd}[/tikz/baseline=(\tikzcdmatrixname-2-1.base)]
	{\Sym H_j} & H \\
	{\Sym \Z_j} & \Z.
	\arrow[from=1-1, to=1-2]
	\arrow["\Sym \psi_j",from=1-1, to=2-1]
	\arrow["\psi", from=1-2, to=2-2]
	\arrow[from=2-1, to=2-2]
\end{tikzcd} \qedhere \]
\end{proof}

We will study three \emph{free} graded algebras that admit $\Z$ as a quotient. The first is the \emph{algebra $\motZ$ of motivic MZV}, which will be the topic of Section \ref{sec:mot-var}. The iterated integral and iterated series representation of MZVs will then lead us to the study of the \emph{(regularized) shuffle and stuffle algebras} $\algshreg$ and $\algstreg$, and consequently the \emph{algebra of formal MZVs}, $\mathcal Z^{\mathfrak f}$ (as explained in Section \ref{sec:conj}).

Recall that a graded $\mathbb Q$-algebra $H = \bigoplus_{k\geq 0} H_k$ is called \textit{of finite type} if $H_0 = \mathbb Q$ and all $H_k$ are finite-dimensional vector spaces. We introduce the following terminology:

\begin{defi}\label{def:rel_algebra_scheme}
    Let $H = \bigoplus_{k\geq 0} H_k$ be a graded $\mathbb Q$-algebra of finite type. For $j \geq 0$, let again
    $$\phi_{H,j}: \Sym H_j \to H$$
    denote the natural morphism induced by mapping $h \in H_j$ to $h\in H$.\footnote{This is the morphism that is adjoint to the inclusion of vector spaces $H_j \to H$ under the adjunction of $\Sym$ and the forgetful functor from algebras to vector spaces.} We define the $j$-th \textit{relation ideal} as $I_{H,j} := \ker \phi_{H,j}$. This is a homogeneous ideal so that
    \[\relalg_j(H) := \Sym H_j / I_{H,j}\]
    is again a graded algebra. We set
    \[\mathcal{R}^{\mathbb Q}_j(H) := \mathrm{Proj} \ \relalg_j(H) \cong V(I_{H,j}) \subseteq \mathbb P(H_j^\vee).\]
\end{defi}

\begin{proposition}
    The association $H \mapsto \relalg_j(H)$ defines an endofunctor of finite type graded $\mathbb Q$-algebras. If $H$ is a domain, then $I_{H,j}$ is a prime ideal, in particular $\relvar_j^{\mathbb Q}(H)$ is a projective variety.
\end{proposition}
\begin{proof}
    If $f: H \to H'$ is a morphism of graded algebras, it induces a map 
    \[\Sym f_j: \Sym H_j\to \Sym H'_j.\]
    By naturality of $\phi_{H,j}$, $\phi_{H',j} \circ \Sym f_j = f \circ \phi_H$. Thus, $\Sym f_j$ induces a map 
    \[\Sym H_j/\ker \phi_{H,j} \to \Sym H'_j/\ker \phi_{H',j}.\]
    
    The second assertion follows because $\Sym H_j / I_{H,j}$ embeds into $H$, so it is a domain if $H$ is. Then $\relvar_j^{\mathbb Q}(H) = \mathrm{Proj} \ \relalg_j(H)$ is irreducible and reduced.
\end{proof}

For the purpose of this paper, we can restrict ourselves to classical algebraic geometry over $\mathbb C$:
\begin{defi}
    Let $H$ be a graded $\mathbb Q$-domain of finite type. We define the \emph{$j$-th (complex) relation variety} of $H$ as
    \[\relvar_j(H) := V(I_{H,j}) \subseteq \mathbb P(H_j^\vee \tensor_{\mathbb Q} \mathbb C).\]
\end{defi}

The following observation will be central in this paper:

\begin{lemma} \label{lem:free_alg_toric}
    If $H$ is a free algebra, then $R_j(H)$ is a semigroup algebra of an affine semigroup. In particular, $\mathcal R_j(H)$ is a projective toric variety.
\end{lemma}
\begin{proof}
    The morphism $\phi_{H,j}$ factors over the subalgebra $Q_{\leq j} \subseteq H$ generated in degree $\leq j$. Since $H$ is free, so is $Q_{\leq j}$. Thus we may write $Q_{\leq j} \cong \kappa[x_1,\ldots,x_n]$, where $\{x_1,\ldots,x_n\}$ is a homogeneous minimal generating set of $Q_{\leq j}$. Then the monomials of degree $j$ in $x_1,\ldots,x_n$ form a basis of $H_j$, so that 
    \[\Sym(H_j) \cong \mathbb Q [ y_m \ | \ m \text{ monomial in } x_1,\ldots,x_n, \deg m = j].\]
    Thus, the map $\Sym(H_j) \to Q_{\leq j}$ is identified with
    \begin{align*}
         \mathbb Q [ y_m \ | \ m \text{ monomial in } x_1,\ldots,x_n, \deg m = j] &\to \mathbb Q[x_1,\ldots,x_n], \quad 
         y_m \mapsto m
    \end{align*}
    so that $\relalg_j(H) \cong \mathbb Q[m \ | \ m \text{ monomial in } x_1,\ldots,x_n, \  \deg m = j]$.
    The latter ring is isomorphic to the $\mathbb Q$ semigroup algebra of the subsemigroup of $\mathbb Z^n$ generated by the finite set
    \begin{equation*}
        \Big\{\alpha\in\mathbb N_0^{n} \ \Big| \ \sum_{i=1}^n \alpha_i|x_i|=j\Big\},
    \end{equation*}
    where $|x_i|$ denotes the degree of $x_i$ according to the grading of $H$.
    As the (homogeneous) spectrum of the semigroup algebra of an affine semigroup, $\relvar_j(H)$ is toric (see for example \cite[Proposition~1.1.14(b)]{coxlittlescheck}).
\end{proof}

\begin{example}
    For the standard graded free algebra $H:= \mathbb Q[x_1,\ldots,x_n]$, $\relvar_j(H)$ is simply the $j$-th Veronese variety.
\end{example}

In view of the example, one might think of $\relvar_j(H)$ as generalized Veronese varieties, see also \cite{roughveronese}. We will see later that the varieties $\relvar_j(H)$ are not necessarily normal.

The next section is dedicated to the study of the toric varieties $\relvar_j(\motZ)$ associated to the algebra $\motZ$ of so-called \emph{motivic} MZVs.

\section{The variety of motivic MZVs}\label{sec:mot-var}

\begin{defi}[Lyndon word]
    \label{def:lyndon} 
    A Lyndon word is a word that is lexicographically strictly smaller than all its cyclic rotations. 
\end{defi}
\noindent For example, in the alphabet generated by letters $\{\letter{2,3}\}$, the word $\letter{22233}$ is Lyndon, while $\letter{22}$ and $\letter{32}$ are not.
\medskip

\noindent
Following \cite[Theorem~8.1]{brown2012mixed} where Brown establishes the structure of the algebra of motivic MZVs (cf. remark below), we take the following as a definition for convenience:
\begin{defi}
    The (graded) \textit{algebra of motivic MZVs} $\motZ$ is the free commutative $\mathbb Q$-algebra generated by formal symbols $\zeta^{\mathfrak m}(w)$ for Lyndon words $w$ over the alphabet $\{\letter{2}, \letter{3}\}$, where the degree of a generator $\zeta^{\mathfrak m}(w)$ is defined to be $\sum_{i} w_i$.
\end{defi}

To every other integer tuple $(k_1,\ldots,k_r)$ with $k_r \geq 2$ one associates a distinguished polynomial $\zeta^{\mathfrak m}(k_1,\ldots,k_r)$ in this algebra, whose degree is the weight $k_1+\ldots+k_r$.

\begin{remark}
Formally, motivic MZVs are defined using the Tannakian formalism for the category \( \MT(\mathbb{Z}) \) of mixed Tate motives over \( \mathbb{Z} \).  The category \( \MT(\mathbb{Z}) \) comes equipped with two fiber functors, the Betti realization and the de Rham realization \( \omega_{\mathrm{B}}, \omega_{\mathrm{dR}} \colon \ \MT(\mathbb{Z}) \to \mathrm{Vec}_{\mathbb{Q}} \).  There is a comparison isomorphism \( \comp_{\mathrm{B},\mathrm{dR}} \colon \ \omega_{\mathrm{dR}} \otimes_{\mathbb{Q}} \mathbb{C}\to \omega_{\mathrm{B}} \otimes_{\mathbb{Q}} \mathbb{C}  \).  The ring of motivic periods of \( \MT(\mathbb{Z}) \) is then,
\[
    \mathcal{P}_{\MT(\mathbb{Z})}^\mathfrak{m} \coloneqq \mathcal{O}(\Isom_{\MT(\mathbb{Z})}(\omega_{\mathrm{dR}}, \omega_\mathrm{B})) \,, 
\]
the ring of regular functions on the scheme of tensor isomorphisms between \( \omega_{\mathrm{dR}} \) and \( \omega_{\mathrm{B}} \).

Brown \cite[Theorem~3.5]{brown2012decomposition} showed that \( \mathcal{P}_{\MT(\mathbb{Z})}^\mathfrak{m} \cong \mathcal{H}[\mathbb{L}^{-1}] \), where \( \mathbb{L} \) is the Lefschetz motive (the motivic analogue of \( 2 \pi i \), and \( \mathcal{H} \) is the \( \mathbb{Q} \)-algebra spanned by the motivic iterated integrals \( I^\mot(a_0; a_1, \ldots, a_n; a_{n+1}) \) on $ \mathbb{P}^1 \setminus \{ 0, 1, \infty \}$ (motivic lifts of \eqref{eqn:idef}). The motivic MZVs are then (cf. \eqref{eqn:zetaasi})
\[
    \zeta^\mot(k_1,\ldots,k_r) = (-1)^r I^\mot(0; \underbrace{1, 0, \ldots, 0}_{k_1}, \underbrace{1, 0, \ldots, 0}_{k_2},  \ldots, \underbrace{1, 0, \ldots, 0}_{k_r}; 1) \,.
\]
Brown \cite[Theorem~1.1, Theorem~8.1]{brown2012mixed} used this framework to show that the elements \( \{ \zeta^\mot(k_1,\ldots,k_r) \mid k_i \in \{ 2, 3 \} \} \) form a basis for the \( \mathbb{Q} \)-vector space of motivic MZVs and that the motivic MZV associated to Lyndon words in $\letter 2,\letter 3$ are free generators of the algebra, i.e.\ $\mathcal{H} \simeq \motZ$.
\end{remark}

The \textit{period map} $\motZ \to \Z ,\zeta^{\mathfrak m}(w) \mapsto \zeta(w)$ is well-defined, surjective and maps $\motZ_j$ into $\Z_j$. By Lemma \ref{lem:inner-ideal-approx} this implies $I_{\motZ,j} \subseteq I_{j}$. Under the \textit{period conjecture} \cite{Brown2014MotivicPeriods} (injectivity of the period map), we have $\motZ \cong \Z$ and thus $I_{\motZ,j} = I_{j}$.

We will now proceed to study the varieties $\mathcal R_j(\motZ)$. Note that they are toric by Lemma \ref{lem:free_alg_toric}.

\begin{notdef}
    We write $\lyndon(\letter{2},\letter{3})$ for the set of Lyndon words in letters $\{\letter{2},\letter{3}\}$, and $\monlyndon(\letter{2},\letter{3})$ for the set of monomials formed out of the variables $\{\zeta^{\mathfrak m}(w)  \mid|  w\in\lyndon(\letter{2},\letter{3})\}$. We will also refer to the degree of $m \in \motZ$ as the \textit{weight} of $m$.

Furthermore, let us write 
\begin{equation*}
\lyndon(\letter{2},\letter{3})_{\leq k}:=\{w\in \lyndon(\letter 2,\letter 3) \mid \weight(w)\leq k\}
\end{equation*}
and
\begin{equation*}
    \monlyndon(\letter{2},\letter{3})_{k}:=\{m\in \monlyndon(\letter 2,\letter 3) \mid \weight(w)= k\}.
\end{equation*}
\end{notdef}

We have $\Sym(\motZ_j) \cong \mathbb Q[y_m \ | \ m \in \mathcal M(\letter{2},\letter{3})_{j}]$ and $\motZ \cong \mathbb Q[\zeta^{\mathfrak m}(w) \ | \ w \in \mathcal L(\letter{2},\letter{3})]$, and under this identification the map $\phi_{\ZmC,j}$ from Definition \ref{def:rel_algebra_scheme} sends $y_m$ to $m$. According to Definition \ref{def:rel_algebra_scheme}, the variety $\mathcal R_j(\motZ)$ is then cut out by the homogeneous ideal
    \begin{equation*}
        I_{\motZ,j}:=\ker(\phi_{\motZ, j} : \mathbb C[y_m \ | \ m \in \mathcal M(\letter{2},\letter{3})_{j}] \to \mathbb C[\zeta^{\mathfrak m}(w) \ | \ w \in \mathcal L(\letter{2},\letter{3})_{\leq j}])
    \end{equation*}
in $\mathbb P((\mathbb C^{\mathcal M(\letter{2},\letter{3})_{j}})^\vee)$.
Here and in the following we just write $\phi_{\motZ,j}$ even though we mean its base change from $\mathbb Q$ to $\mathbb C$, and restrict its codomain from $\motZ$ to the subalgebra generated by the (finitely many) elements of weight lower or equal $j$. Note that this restriction does not change the kernel since the map is graded and the domain is generated in weight $j$.

\begin{example}
    In the smallest case $j=2$, we have $\mathcal L(\letter 2,\letter 3)_{\leq 2} = \{\letter 2\}$ and $\mathcal M(\letter 2,\letter 3)_2 = \{\zeta^{\mathfrak m}(2) \}$. The map $\phi_{\ZmC,j}: \mathbb C[y_{\zeta^{\mathfrak m}(2)}] \to \mathbb C[\zeta^{\mathfrak m}(2)]$ is an isomorphism, and so $I_{\ZmC, 2} = (0)$ and $\mathcal R_2(\ZmC)$ is just a point. In fact, it turns out that $I_{Z^{\mathfrak m}, j} = (0)$ for all $j<10$, and accordingly, $\mathcal R_j(\ZmC)$ is simply projective space.
\end{example}

By classic toric geometry, the ideal $I_{\ZmC,j}$ is linearly generated by the binomials
$$\{y^u - y^v \ | \ A_{Z^{\mathfrak m},j}u = A_{\ZmC,j}v\}$$
where $A_{\ZmC,j}$ is the linear map $\mathbb Z^{\mathcal M(\letter 2,\letter 3)_j} \to \mathbb Z^{\mathcal L(\letter 2,\letter 3)_{\leq j}}$ induced by mapping a monomial $m$ to its exponent vector.
\begin{example}\label{ex:mot-weight-10}
    Let us determine the variety $\mathcal R_{10}(\ZmC)$. The Lyndon words in weight lexicographic order up to weight 10 are
    \begin{equation*}
        \letter 2, \letter 3, \letter {23}, \letter {223}, \letter {233},\letter {2223},\letter {2233}.
    \end{equation*}
    Writing $x_w := \zeta^{\mathfrak m}(w)$, the monomials in weight 10 are
    \begin{equation*}
        x_{\letter 2}^5,\,x_{\letter 2}^2x_{\letter 3}^2,\,x_{\letter 2}x_{\letter 3}x_{\letter{23}},\,x_{\letter 2}x_{\letter {233}},\,x_{\letter 3}x_{\letter{223}},\,x_{\letter {23}}^2,\,x_{\letter {2233}}.
    \end{equation*}
    In particular, the ambient projective space is $6$-dimensional. Using the given ordering, we can represent the map $A_{\ZmC,j}$ as the following matrix:
    \begin{equation*}
        A_{Z^{\mathfrak m},10}=
        \begin{pmatrix}
            5 & 2  & 1 & 1 & 0 & 0 & 0\\
            0 & 2 & 1 &  0 & 1 & 0 & 0\\
            0 & 0 & 1 & 0 & 0 & 2 & 0\\
            0 & 0 & 0 & 0 & 1 & 0 & 0\\
            0 & 0 & 0 & 1 & 0 & 0 & 0\\
            0 &0 & 0 & 0  & 0 & 0 & 0\\
            0& 0 & 0 & 0& 0 & 0 &1
        \end{pmatrix}
    \end{equation*}
    Its kernel is generated by the vector $(0,1,-2,0,0,1,0)$, and consequently we have 
    \begin{equation}
    \label{eq:exmot10ideal}
        I_{Z^{\mathfrak m},j} = (y_{x_{\letter 2}x_{\letter 3}x_{\letter{23}}}^2 - y_{x_{\letter 2}^2x_{\letter 3}^2} \cdot y_{x_{\letter {23}}^2}).
    \end{equation}
    It follows that $\mathcal R_{10}(\ZmC)$ is a hypersurface, given by a cone over the rational normal curve in $\mathbb P^2$.
\end{example}

\begin{example}\label{ex:mot-weight-10-expansion}
    Replacing $y_m$ by $m$, expanded using the shuffle product, and using the linear relations of MZVs, 
    the quadric from Example \ref{ex:mot-weight-10} becomes a non-trivial relation of seven weight $10$ MZVs:
    \par
    \vspace{-\baselineskip + 1ex}
    \begin{scriptsize}
     \begin{align*}
    & -\frac{3793}{360} \zeta (10) \zeta (1,9)+\frac{3407}{360} \zeta (10) \zeta (3,7)+\frac{27119}{540} \zeta (10) \zeta (5,5)+\frac{6107}{540} \zeta
   (10) \zeta (1,1,8) \\
   &+\frac{3407}{1080} \zeta (10) \zeta (1,4,5) -\zeta (1,9)^2-\zeta (3,7)^2+\frac{224}{9} \zeta (5,5)^2-\frac{4}{9} \zeta
   (1,1,8)^2-\frac{1}{9} \zeta (1,4,5)^2-2 \zeta (1,9) \zeta (3,7) \\
   &-\frac{118}{3} \zeta (1,9) \zeta (5,5)+\frac{2}{3} \zeta (3,7) \zeta
   (5,5)-\frac{4}{3} \zeta (1,9) \zeta (1,1,8)-\frac{4}{3} \zeta (3,7) \zeta (1,1,8)+\frac{94}{9} \zeta (5,5) \zeta (1,1,8)\\
   &-\frac{2}{3} \zeta
   (1,9) \zeta (1,4,5) -\frac{2}{3} \zeta (3,7) \zeta (1,4,5)+\frac{2}{9} \zeta (5,5) \zeta (1,4,5)-\frac{4}{9} \zeta (1,1,8) \zeta
   (1,4,5)-\frac{1012849}{518400}\zeta (10)^2 = 0
 \end{align*}
\end{scriptsize}
\end{example}

Let us define $P_{\ZmC, j}$ as the convex hull of the points $A_{\ZmC, j}(m)$ in $\mathbb Z^{\mathcal L(\letter 2, \letter 3)_{\leq j}} \subseteq \mathbb R^{\mathcal L(\letter 2, \letter 3)_{\leq j}}$, where $m \in \mathcal M(\letter 2, \letter 3)_j$. This is the lattice polytope associated to the (not necessarily normal) projective toric variety $\mathcal R_j(Z^{\mathfrak m})$ in its given embedding. An (embedded) toric variety and its lattice polytope are very closely related. We recall some classical results from toric geometry, known in folklore terms as the 'toric variety -- polytope dictionary'. We refer to \cite{coxlittlescheck} and \cite{telen2022introductiontoricgeometry} for the relevant definitions.

\begin{theorem}\label{thm:toric-recall} 
    Let $X$ be a projective toric variety embedded in $\mathbb P^n$, with vanishing ideal $I:= (y^u - y^v \ | \ A(u)=A(v))$ for some linear map $A:\mathbb Z^m \to \mathbb Z^n$. Let $L$ be the affine sublattice of $\mathbb Z^{n}$ generated by the vectors $A(e_i)$, for $e_i$ the standard basis of $\mathbb Z^m$, and let $P$ be the convex hull of the $A(e_i)$.
    \begin{enumerate}
        \item The dimension of $X$ is the dimension of $P$.
        \item The degree of $X\subseteq \mathbb P^n$ is the normalized volume of $P$ in the affine sublattice $L$.
        \item $X$ is normal if and only if $P$ is very ample in $L$.
        \item $X$ is projectively normal if and only if $P$ is normal in $L$ (or equivalently, the semigroup spanned by the $A(e_i)$ is saturated in the (non-affine) sublattice they generate).
        \item $X$ is smooth if and only if $P$ is smooth in $L$.
    \end{enumerate}
\end{theorem}
\begin{proof}
\begin{enumerate}
\item  follows from \cite[Theorem~1.2.18]{coxlittlescheck}, using that the affine cone over $X$ is the affine toric variety associated to the cone spanned by $P$.
See also \cite[Corollary~3.11]{telen2022introductiontoricgeometry} for an argument via the rank of the lattice.

\item is known as Kushnirenko's theorem and can be found in this form e.g.\ in \cite[Proposition~4.16]{Sturmfels1996}. Also see \cite[Section~3.4]{telen2022introductiontoricgeometry} and \cite[Section~2.2]{IbadanlecturesSotille} for nice expositions.

\item follows from three facts: first, the affine toric varieties corresponding to the semigroups spanned by $A(e_i) - v$ for vertices $v$ of $P$ are coordinate rings of affine patches covering the projective toric variety $X$. This is \cite[Propositions~2.1.8 and 2.1.9]{coxlittlescheck}. Second, by definition $P$ is very ample in $L$ if and only if the semigroups spanned by $A(e_i)-v$ are saturated, which is the case if and only if the associated affine toric varieties are normal (see e.g.\ statement and proof of \cite[Theorem~1.3.5]{coxlittlescheck}). Finally, normality is a local property of varieties.

\item follows again from \cite[Theorem~1.3.5]{coxlittlescheck}: by definition, $X$ is projectively normal if and only if its homogeneous coordinate ring is normal.

\item follows from \cite[Theorem~1.3.12]{coxlittlescheck} by an analogous argument as for (3). We refer to the proof of \cite[Theorem 2.4.3]{coxlittlescheck}.
\qedhere
\end{enumerate}
\end{proof}

\begin{figure}[h]
\centering
    \begin{subfigure}{0.4\linewidth}
        \includegraphics[width=\linewidth]{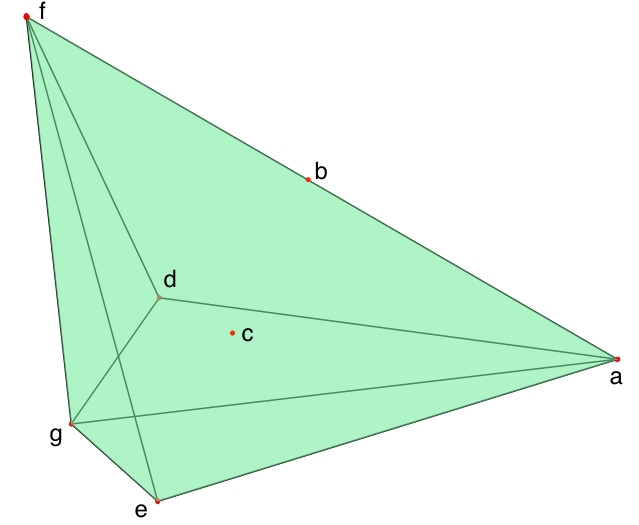}
        \subcaption{}
    \end{subfigure}
    \hspace{3em}
    \begin{subfigure}{0.4\linewidth}    \includegraphics[width=\linewidth]{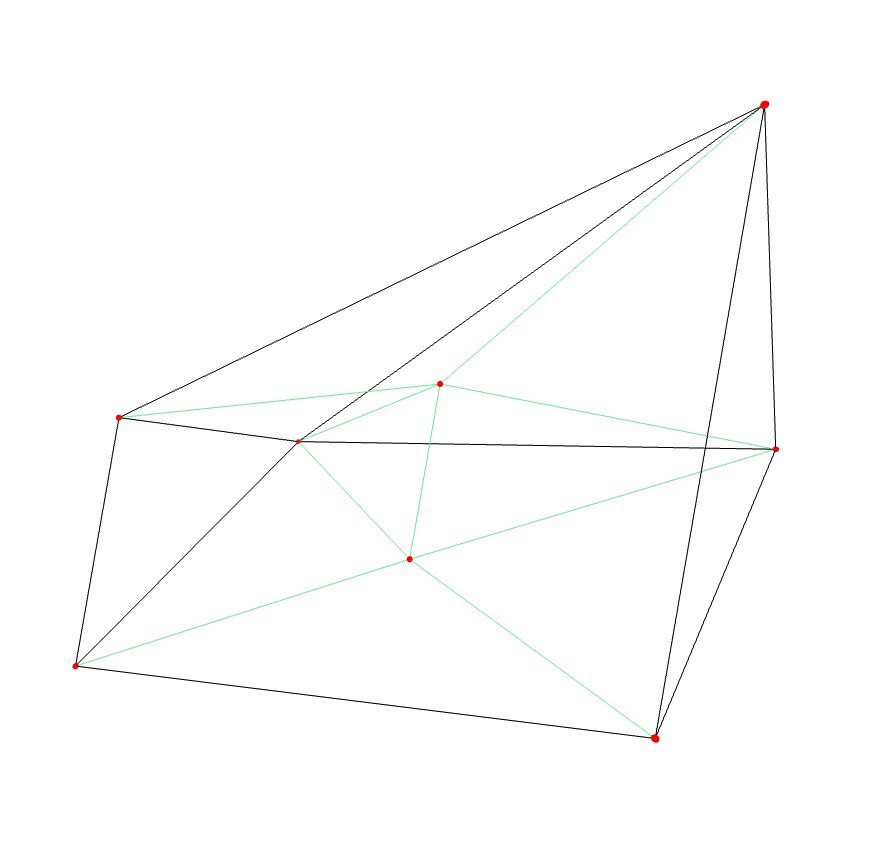}
        \subcaption{}
    \end{subfigure}
    \caption{The three-dimensional polytope associated to the variety $\mathcal R_{12}^0(\ZmC)$ from Example \ref{ex:mot-w12} (a), and a Schlegel diagram of the four-dimensional polytope associated to $\mathcal R_{13}^0(\ZmC)$ (b).}\label{fig:poly-d3}\label{fig:poly-d4}
\end{figure}

\begin{example}\label{ex:mot-w12}
Let us study the variety $\mathcal R_{12}(\ZmC)$. The Lyndon words in weight lexicographic order up to weight 12 are
    \begin{equation*}
        \letter 2, \letter 3, \letter {23}, \letter {223}, \letter {233},\letter {2223},\letter {2233}, \letter{22223},\letter{2333},\letter{22233},\letter{22323}
    \end{equation*}
    Accordingly, writing $x_w := \zeta^{\mathfrak m}(w)$, the monomials in weight $12$ are
    \begin{equation*}
        x_{\letter 2}^6,\,
        x_{\letter 2}^3x_{\letter 3}^2,\,
        x_{\letter 2}^2x_{\letter 3}x_{\letter {23}},\,
        x_{\letter 2}^2x_{\letter {233}},\,
        x_{\letter 2}x_{\letter 3}x_{\letter {223}},\,
        x_{\letter 2}x_{\letter {23}}^2,x_{\letter 2}x_{\letter {2233}},x_{\letter 3}^4,x_{\letter 3}x_{\letter {2223}},x_{\letter {23}}x_{\letter {223}},x_{\letter{22233}},x_{\letter{22323}},
    \end{equation*}
    so the ambient projective space has dimension $11$. Using the given ordering, we can write
    \begin{equation*}
        A_{Z^{\mathfrak m},12}=
        \begin{pmatrix}
            6 & 3 & 2 & 2 &  1 & 1  & 1 & 0 & 0 & 0 & 0 & 0\\
            0 & 2 & 1 & 0 & 1 & 0 & 0 & 4 & 1 & 0 & 0 & 0\\
            0 & 0 & 1 & 0 & 0 & 2 & 0 & 0 & 0 & 1 & 0 & 0\\
            0 & 0 & 0 & 0 & 1 & 0 & 0  & 0 & 0 & 1 & 0 & 0\\
            0 & 0 & 0 & 1 & 0 & 0 & 0 &0 &0 & 0 &0 &0\\
            0 & 0 & 0 & 0 & 0 &0 &0 &0 & 1 & 0 & 0 & 0\\
            0 & 0 &0 & 0 &0 & 0 & 1 & 0 & 0 & 0 &0 &0\\
            0 & 0 &0 & 0 &0 &0 &0 &0  & 0  &0 &0 &0\\
            0 & 0 &0 & 0 &0 &0 &0 &0  & 0  &0 &0 &0\\
            0 & 0 &0 & 0 &0 &0 &0 &0  & 0  &0 &1 &0\\
            0 & 0 & 0 & 0 & 0 & 0 & 0 & 0 & 0 & 0 & 0 & 1
            \end{pmatrix}
    \end{equation*}
    The polytope $P_{Z^{\mathfrak m},12}$ is the convex hull of the columns, and one checks that it is $8$-dimensional. The normalized volume of $P_{Z^{\mathfrak m},12}$ in the lattice affinely generated by the columns of $A_{\ZmC,12}$ is $6$.

    Thus, by Theorem \ref{thm:toric-recall} we know that $\mathcal R_{12}(\ZmC)$ has dimension $8$ and degree $6$ inside the $11$-dimensional ambient space. It is cut out by the four quadrics 
    \[af-b^2,\quad be-c^2, \quad bg-cd, \quad cg - de,\]
    where $a = y_{x_{\letter 2}^6}$, $b = y_{x_{\letter 2}^3x_{\letter 3}^2}$, $c = y_{x_{\letter 2}^2x_{\letter 3}x_{\letter {23}}}$, $d = y_{x_{\letter 2}x_{\letter 3}x_{\letter {223}}}$, $e = y_{x_{\letter 2}x_{\letter {23}}^2}$, $f = y_{x_{\letter 3}^4}$ and $g = y_{x_{\letter {23}}x_{\letter {223}}}$. We express these relations in a certain depth minimal basis in Appendix \ref{sec:wt12deg2ex}.

    Note that in particular only $7$ out of the $12$ variables $y_m$ appear in the relations generating $I_{Z^{\mathfrak m}, j}$; so $\mathcal R_{12}(\ZmC)$ is again a cone, over a three-dimensional toric subvariety $\relvarr_{12}(\ZmC)$ of $\mathbb P^6$. The associated three-dimensional polytope is shown in Figure \ref{fig:poly-d3}.
\end{example}

\begin{notdef}
    As seen in Example \ref{ex:mot-weight-10} and Example \ref{ex:mot-w12}, $\relvar_j(\ZmC)$ is generally a cone over a projective variety in a lower-dimensional projective subspace. We will denote the base of this cone by $\relvarr_j(\ZmC)$ and view it as a subvariety of the smallest linear subspace containing it. We write $P^0_{\ZmC,j}$ for the associated polytope.
\end{notdef}

Let us denote the subalgebra of $\ZmC$ generated in weight $\leq j$ by $\ZmC_{\langle \leq j\rangle}$.

\begin{proposition}\label{prop:rel0-is-rel}
    For $j \geq 11$, we have
    \[\relvarr_j(\ZmC) = \relvar_j(\ZmC_{\langle \leq j-5 \rangle}).\]
\end{proposition}
\begin{proof}
    We need to show that an irreducible relation in $I_{\ZmC,j}$ contains a variable $y_m$ if and only if the variables appearing in the monomial $m$ all have weight $\leq j -5$.

    So let $f$ be an irreducible relation in $I_{\ZmC,j}$. Since this is a toric ideal, we may assume that $f$ is a binomial, $f=g-h$ for monomials $g$ and $h$. Assume $f$ contains a variable $y_m$ such that $m$ contains a generator $z$ of weight $>j-5$. Wlog., $y_m$ appears in $g$. Then $m$ is the unique monomial containing $z$. Since $f\in I_{\ZmC,j}$, we have $\phi_{\ZmC,j}(g) = \phi_{\ZmC,j}(h)$; so since $m$ is unique, $y_m$ must divide $h$. This contradicts irreduciblity.
    
    It remains to show that every variable $y_m$ where $m$ is a monomial of weight $j$ in variables of weight $\leq j-5$ appears in an irreducible relation in $I_{\ZmC,j}$. This follows from the following claim: we can write $m = a \cdot b$ where $a$ and $b$ are monomials of weight $\geq 5$.
    Indeed, if $m$ contains a variable of weight $\geq 5$, take $a$ to be this variable. If not, only the variables $\letter 2$ and $\letter 3$ can appear in $m$, and we may take $a$ to be $\letter 2^3$, $\letter 2\cdot\letter 3$ or $\letter 3^2$ depending on $m$, all of weight $\leq j-5$. In any case, take $b=m/a$ of weight $\geq j-(j-5)=5$.

    Now write $m = ab$ with $a$ and $b$ of weight $\geq 5$. For $j \geq 5$ the vector space $\ZmC_j$ has dimension $\geq 2$, that is, there are at least two different monomials of that weight. So we can choose a monomial $a'\not=a$ of the same weight as $a$ and a monomial $b'\not=b$ of the same weight as $b$. Then we have the non-zero irreducible relation
    \(y_{ab}y_{a'b'} - y_{ab'}y_{a'b} \).
\end{proof}

\begin{table}
    \centering
    \renewcommand{\arraystretch}{1.2}
    \begin{tabular}{c|ccccccccc}
        Property $\backslash$ $j$ & 10 & 11 & 12 & 13 & 14 & 15 & 16 & 17 & 18\\\hline\vspace{-1em} \\ 
        $\mathrm{dim} \ \relvar_j(\ZmC)$ & 5 & 7 & 8 & 11 & 13 & 17 & 21 & 27 & 33\\
        amb.\ $\dim \ \relvar_j(\ZmC)$ & 6 & 8 & 11 & 15 & 20 & 27 & 36 & 48 & 64\\
        $\mathrm{dim} \ \relvarr_j(\ZmC)$ & 1 & 2 & 3 & 4 & 5 & 6 & 8 & 10 & 13\\
        amb.\ $\dim \ \relvarr_j(\ZmC)$ & 2 & 3 & 6 & 8 & 12 & 16 & 23 & 31 & 44\\
        $\mathrm{deg} \ \relvar_j(\ZmC)$ & 2 & 2 & 6 & 7 & 20 & 36 & 98 & 231 & 984\\
        $\relvarr_j(\ZmC)$ is proj.\ normal & \true & \true & \true & \true & \true & \true & \false & \true & \false\\
        $\relvarr_j(\ZmC)$ is normal & \true & \true & \true & \true & \true & \true & \false & \true & \false\\
        $\relvarr_j(\ZmC)$ is smooth & \true & \true & \false & \false & \false & \false & \false & \false & \false\\
        generators of $I_{\ZmC, j}$ & $1_2$ & $1_2$ & $4_2$ & $8_2$ & $18_2$ & $36_2$ & $73_2$ & $140_2$ & $278_2,2_3$\\
    \end{tabular}

    \vspace{1em}
\begin{tabular}{c|ccccc}
        Property $\backslash$ $j$ & 19 & 20 & 21 & 22 & 23\\\hline\vspace{-1em} \\
        $\mathrm{dim} \ \relvar_j(\ZmC)$ & 43 & 53 & 68 & 85 & 109\\
        amb.\ $\dim \ \relvar_j(\ZmC)$ & 85 & 113 & 150 & 199 & 264\\
        $\mathrm{dim} \ \relvarr_j(\ZmC)$ & 16 & 20 & 25 & 32 & 40\\
        amb.\ $\dim \relvarr_j(\ZmC)$ & 58 & 80 & 107 & 146 & 195\\
        $\mathrm{deg} \ \relvar_j(\ZmC)$ & 2585 & 15330 & 71276 & 657464 & 4608050\\
        $\relvarr_j(\ZmC)$ is proj.\ normal & \true & \false & \false & \false & \false\\
        $\relvarr_j(\ZmC)$ is normal & \true & \false & \false & \false & \false\\
        $\relvarr_j(\ZmC)$ is smooth & \false & \false & \false & \false & \false\\
        generators of $I_{\ZmC, j}$ & $514_2,6_3$ & $990_2,9_3$ & $1830_2, 38_3, 28_4$ &  &  \\
    \end{tabular}
    \caption{Properties of the projective varieties $\relvar_j(\ZmC)$ and $\relvarr_j(\ZmC)$. The empty cells in the table represent calculations which did not terminate in reasonable time due to computational limits.}
    \label{tab:mot}
\end{table}

Table \ref{tab:mot} shows interesting properties of the varieties $\relvar_j(\ZmC)$ and $\relvarr_j(\ZmC)$ for small weights $j$. This table was produced by studying the associated polytopes $P_{\ZmC,j}$ and $P^0_{\ZmC,j}$ in \texttt{Oscar.jl} and applying Theorem \ref{thm:toric-recall}.

Let us now compute the dimensions of the varieties $\relvar_j(\ZmC)$ and $\relvarr_j(\ZmC)$ in general.

\begin{lemma}\label{lem:gen-lattice}
    Fix a weight compatible total order of $\lyndon(\letter 2,\letter 3)$.
    For any $j\geq 6$ even or $j\geq 9$ odd and any $w\in\lyndon(\letter 2,\letter 3)$ with $2\leq \weight(w)\leq j-2$,
    \begin{equation*}
        M_j\cap \mathbb Z^{\{v\in\lyndon(\letter 2,\letter 3)|v\leq w\}}, \text{ where } M_j:=\Big\{\sum_i a_i w_i\,|\,a_i\in\mathbb N_0,\,w_i\in\lyndon (\letter 2,\letter 3)_{\leq j-2},\,\sum_i a_i\weight(w_i)=j\Big\}
    \end{equation*}
    generates the lattice
    \begin{equation*}
        L_j\cap \mathbb Z^{\{v\in\lyndon(\letter 2,\letter 3)|v\leq w\}}, \text{ where } L_j:=\Big\{\sum_i a_i w_i\,|\,a_i\in\mathbb Z,\,w_i\in\lyndon (\letter 2,\letter 3)_{\leq j-2},\,\sum_i a_i\weight(w_i)\in j\mathbb Z\Big\}.
    \end{equation*}
    In particular, $M_j$ generates $L_j$.
\end{lemma}
\begin{proof}

    For $j\geq 6$ even, let $\alpha=\tfrac{j}{2}\letter 2$ and $\beta=\tfrac{j-6}{2}\letter 2+2\cdot \letter 3$, both in $M_j$. 
    For $j\geq 9$ odd, let $\alpha=\tfrac{j-3}{2}\cdot\letter 2+\letter 3$,  $\beta=\tfrac{j-9}{2}\letter 2+3\cdot \letter 3$, again both in $M_j$.

    In both cases, $\beta-\alpha=-3 \cdot \letter 2+2\cdot \letter 3$,
    thus $\alpha,\beta$ generate $L_j\cap \mathbb Z^{\{\letter 2,\letter 3\}}$, as whenever for $n,c_1,c_2,d_1,d_2\in\mathbb Z$ we have
    \begin{equation*}
        2c_1+3c_2=nj=2d_1+3d_2,
    \end{equation*}
    then $(d_1-c_1,d_2-c_2)$ is a $\mathbb Z$ multiple of $(3,-2)$

    Let now $w\in\mathcal{L}(\letter 2,\letter 3)$ with $5\leq\weight(w)\leq j-2$ and assume $L_j\cap \mathbb Z^{\{v\in\lyndon(\letter 2,\letter 3)|v<w\}}$ is generated by $M_j\cap \mathbb Z^{\{v\in\lyndon(\letter 2,\letter 3)|v<w\}}$.
    Let $y\in L_j\cap \mathbb Z^{\{v\in\lyndon(\letter 2,\letter 3)|v\leq w\}}$ arbitrary and let $n\in\mathbb Z$ be the coefficient of $w$ in $y$.
    We can find $a,b\in\mathbb N_0$ such that $x=a \cdot \letter{2} + b \cdot \letter{3}+w$ has weight $j$ and thus $x\in M_j$.
    Thus 
    \begin{equation*}
        y=nx+(y-nx)
    \end{equation*}
    is in the lattice generated by $M_j\cap \mathbb Z^{\{v\in\lyndon(\letter 2,\letter 3)|v\leq w\}}$ as $y-nx\in L_j\cap \mathbb Z^{\{v\in\lyndon(\letter 2,\letter 3)|v<w\}}$.
    \end{proof}

    \begin{theorem}\label{thm:mot-dim}
    Let us denote the number of Lyndon words of weight $j$ in letters $\letter{2}$ and $\letter{3}$ by $\lambda_{\letter{23},j}$. Then for $j\geq 10$, we have
    \[\mathrm{dim} \ \relvar_j(\ZmC) = \lambda_{\letter{23},1} + \ldots + \lambda_{\letter{23},j-2} + \lambda_{\letter{23},j} - 1\]
    and for $j \geq 11$,
    \[\mathrm{dim} \ \relvarr_j(\ZmC) = \lambda_{\letter{23},1} + \ldots + \lambda_{\letter{23},j-5} - 1.\]
\end{theorem}
\begin{proof}
    By Theorem \ref{thm:toric-recall} the dimension of $\relvar_j(\ZmC)$ and $\relvarr_j(\ZmC)$ agrees with the dimensions of the polytopes $P_{\ZmC,j}$ and $P^0_{\ZmC,j}$. Thus the statement follows from Lemma \ref{lem:gen-lattice} and Proposition \ref{prop:rel0-is-rel}.
\end{proof}

\begin{lemma}
We have
\begin{align*}
\lambda_{\letter{23},j}= \sum_{\substack{a,b\geq0 \\ 2a+3b=j}} \frac{1}{a+b}\sum_{d\mid \gcd(a,b)} \mu(d)\binom{(a+b)/d}{a/d}= \frac{1}{j}\sum_{d\mid j} \mu\left(\frac{j}{d}\right) \mathsf{Per}(d),
\end{align*}
where $\mu$ is the M\"obius function and $\mathsf{Per}(n)$ the Perrin sequence, i.e., we have $\mathsf{Per}(0)=3$, $\mathsf{Per}(1)=0$, $\mathsf{Per}(2)=2$, and $\mathsf{Per}(n)=\mathsf{Per}(n-2)+\mathsf{Per}(n-3)$. 
\end{lemma}

For the values of $\lambda_{\letter{23},j}$ for small $j$, see the OEIS sequence \href{https://oeis.org/A113788}{A113788}.

\begin{proof}
Let $\mathsf{W}(a,b)$ be the number of words containing exactly $a$ times the letter $\letter 2$ and $b$ times the letter $\letter 3$. Then, we have
\begin{align} \label{eq:number_words_23}
|\mathsf{W}(a,b)|=\binom{a+b}{a}.
\end{align} 
For the first equality, it suffices to show that the number of Lyndon words in $\mathsf{W}(a,b)$ is equal to
\begin{align*} 
\frac{1}{a+b}\sum_{d\mid \gcd(a,b)} \mu(d) \binom{(a+b)/d}{a/d}.
\end{align*}
Any word $w\in \mathsf{W}(a,b)$ can be written as $w=v^d$, where $d\geq1$ and $v\in \mathsf{W}(\frac{a}{d},\frac{b}{d})$ is primitive (i.e., not some power of a proper subword). We deduce
\[
|\mathsf{W}(a,b)|=\sum_{d\mid \gcd(a,b)} \left|\mathsf{W}^{\operatorname{prim}}\left(\frac{a}{d},\frac{b}{d}\right)\right|
\]
where $\mathsf{W}^{\operatorname{prim}}(a,b)$ is the set of primitive words in $\mathsf{W}(a,b)$. Möbius inversion gives
\begin{align} \label{eq:formula_Wprim}
|\mathsf{W}^{\operatorname{prim}}(a,b)|=\sum_{d\mid \gcd(a,b)} \mu(d) \left|\mathsf{W}\left(\frac{a}{d},\frac{b}{d}\right)\right|.
\end{align}
All Lyndon words in $\mathsf{W}(a,b)$ are also contained in $\mathsf{W}^{\operatorname{prim}}(a,b)$ and cyclic permutation of words preserves this space. If we have a primitive word of length $a+b$, then all its $a+b$ cyclic rotations are distinct and by definition exactly one of them is a Lyndon word. Hence, the number of Lyndon words in $\mathsf{W}(a,b)$ is equal to
\[
\frac{1}{a+b} |\mathsf{W}^{\operatorname{prim}}(a,b)|=\frac{1}{a+b}\sum_{d\mid \gcd(a,b)} \mu(d) \binom{(a+b)/d}{a/d},
\]
where the equality follows from \eqref{eq:formula_Wprim} and \eqref{eq:number_words_23}.
This proves the first claimed equality. For the second one see \cite[p. 19]{WaldschmidtMZV}.
\end{proof}
    
    \begin{proposition}
        For $j\geq 5$, the lattice $L_j$ has index $j$ in $\mathbb Z^{\lyndon(\letter 2,\letter 3)_{\leq j-2}}$.
    \end{proposition}
    \begin{proof}
        By definition, $L_j$ is the kernel of the morphism
        \[\mathbb Z^{\lyndon(\letter 2,\letter 3)_{\leq j-2}} \to \mathbb Z / j\mathbb Z, \ w_i \mapsto \mathrm{wt}(w_i) \bmod j\]
        which is surjective for $j\geq 5$.
    \end{proof}
    \begin{corollary}
        For $j\geq 11$ odd, we have
        \begin{equation*}
            \deg \relvar_j(\ZmC)=\deg \relvar_j^0(\ZmC)=(\dim\relvar_j^0(\ZmC))\operatorname{vol}(P^0_{\ZmC,j})/j,
        \end{equation*}
        where $\operatorname{vol}(P^0_{\ZmC,j})$ denotes the
    euclidean volume of $P^0_{\ZmC,j}$ in the weight $j$ hyperplane of $\mathbb R^{\lyndon(\letter 2,\letter 3)_{\leq j-5}}$.
    \end{corollary}

As already observed in Table \ref{tab:mot}, the varieties of motivic MZVs are in general neither projectively normal nor normal. In fact, they are normal only in a finite number of cases:

\begin{theorem}\label{thm:mot-non-normal}
 The varieties $\relvar_j(\ZmC)$ and $\relvarr_j(\ZmC)$ are normal only for $j\leq 15$, $j = 17$ and $j=19$.
\end{theorem}
\begin{proof}
    We show that $\relvarr_j(\ZmC)$ and thus $\relvar_j(\ZmC)$ is not normal for the remaining $j$. Let $z_0$ be a Lyndon word of weight $j-5$ and set $z = z_0 + \letter{23}$. Recall that $A_{\ZmC,j}$ is the linear map $\mathbb Z^{\mathcal M(\letter 2,\letter 3)_j} \to \mathbb Z^{\mathcal L(\letter 2,\letter 3)_{\leq j}}$ induced by mapping a monomial $m$ to its exponent vector. We claim that the semigroup $N$ spanned by the vectors $A_{\ZmC,j}(m)-z$ for monomials $m$ in generators of weight $\leq j-5$ is not saturated in the lattice $L$ it spans. 
    Thus, the polytope $P^0_{\mathcal{Z}^{\mathfrak m},j}$ associated to $\relvarr_j(\mathcal{Z}^{\mathfrak m})$ is not very ample, and we conclude by Theorem \ref{thm:toric-recall} that $\relvarr_j(\mathcal{Z}^{\mathfrak m})$ is not normal.

    For $a \in \mathbb Z^{\mathcal L(\letter 2, \letter 3)_{\leq j}}$, we call $\sum_{w \in \mathcal L(\letter 2, \letter 3)_{\leq j}} a_w \cdot \mathrm{wt}(w)$ the   of $a$. 
    Let us first assume that $j$ is even. Choose Lyndon words $w$, $v$ and $u$ such that $w$ has weight $j-6$, $v$ has weight $\frac{j}{2}$ and $u$ has weight $\frac{j-2}{2}$. This is possible for $j \geq 16$. Set $x := w+v+u+ \letter{3} + 2\cdot \letter{2} - 2z$. Since $x+2z$ has weight $2j$, it lies in the lattice spanned by the $A_{\ZmC,j}(m)$ by Lemma \ref{lem:gen-lattice}. Because of the grading, this implies that $x \in L$. However, we claim that $2x \in N$ while $x \notin N$. Indeed, the elements of $N$ are precisely the sums
    \[y = y_1 + \dots + y_n - n \cdot \letter{23} - n \cdot z_0 + m (\letter{2} + \letter{3} - \letter{23})\]
    where $m,n \in \mathbb N_0$ and $y_1,\ldots,y_n$ are positive vectors of weight $j$ where $z_0$ does not appear. Assuming that $x$ is of this form, we compare coefficients of $z_0$ and $\letter{3}$ and see that we must have $n=2$, $m\leq 1$. But then $y_1, y_2$ must partition the multiset $\{j-6, \frac{j}{2}, \frac{j-2}{2}, 3, 2, 2\}$ (for $m=0$) or $\{j-6, \frac{j}{2}, \frac{j-2}{2}, 5, 2\}$ (for $m=1$) into two sets that sum to $j$. This is impossible since no subset of $\{\frac{j}{2}, \frac{j-2}{2}, 3, 2, 2\}$ sums to $6$ for $j \geq 16$. Thus $x \notin N$. However, we can write
    \[2x = (w + 3 \cdot \letter 2 - z) + (w + 2 \cdot \letter 3 - z) + (2\cdot v - z) + (2\cdot u + \letter 2 - z)\]
    as a sum of generators, and thus $2x \in N$.

    If $j$ is odd, we choose $v$ instead in weight $\frac{j-5}{2}$ and $u$ in weight $\frac{j-7}{2}$. This is possible for $j \geq 21$. Set $x := w + v + u + \letter{23} + \letter{3} + 2 \cdot \letter{2} - 2z$. One argues as in the even case that $x \in L$ and that $x \in N$ implies that we can write $x = (y_1 - z) + (y_2 - z) + m(\letter 2 + \letter 3 - \letter{23})$ with $m \leq 1$ and $y_1,y_2$ positive vectors of weight $j$ where $z_0$ does not appear. Again, this forces $y_1,y_2$ to correspond to partitions of the multiset $\{j-6,\frac{j-5}{2},\frac{j-7}{2},5,3,2,2\}$ (for $m=0$) or $\{j-6,\frac{j-5}{2},\frac{j-7}{2},5,5,2\}$ (for $m=1$). For $j\geq 21$, no subset of $\{\frac{j-5}{2},\frac{j-7}{2},5,3,2,2\}$ sums to $6$ and so this is not possible. Thus $x \notin N$, but  $2x \in N$ since we can write
    \[2x = (w + 3 \cdot \letter 2 - z) + (w + 2 \cdot \letter 3 - z) + (2v + \letter 5 - z) + (2u + \letter {23} + \letter 2 - z)\]
    as a sum of generators.

    The statement thus follows from Table \ref{tab:mot}.
    \end{proof}

    \begin{example}
        For $j=16$, we may choose
        \begin{align*}
            x &= \letter{2233} + \letter{233}+\letter{223}+\letter{3}+2\cdot\letter{2}-2\cdot\letter{2333}-2\cdot\letter{23}
        \end{align*}
        to observe that $x\in L$,
        yet $x\notin N$ since neither
        \begin{align*}
            x' = \letter{2233} + \letter{233}+\letter{223}+\letter{3}+2\cdot\letter{2}
        \end{align*}
        nor
        \begin{equation*}
            x''=\letter{2233} + \letter{233}+\letter{223}+\letter{23}+\letter{2}
        \end{equation*}
        split into two sums (with nonnegative coefficients) each of weight $16$.
    \end{example}

\begin{corollary}
    The variety $\relvarr_j(\ZmC)$ is only smooth for $j\leq 11$.
\end{corollary}
\begin{proof}
    This follows from Table \ref{tab:mot} and Theorem \ref{thm:mot-non-normal}, as smoothness implies normality.
\end{proof}

As already seen in Table \ref{tab:mot}, the ideal $I_{\ZmC, j}$ is in general not generated by quadrics. In fact, something stronger is true:

\begin{proposition}\label{prop:not_gen_in_fixed_deg}
    For every $n$ there is a $j$ such that $I_{\ZmC,j}$ is not generated in degree $\leq n$.
\end{proposition}
\begin{proof}
    The following argument is due to \cite[Proposition~2.9]{colmenarejo2020toric}.

    Fix $n$. Choose $j$ such that there is a rigid $j$-matrix $M$ of size $n+1 \times n+1$, with all entries $\geq 7$. That is, $M$ is a positive integer matrix such that
    \begin{enumerate}
        \item the entries of each row and column of $M$ sum to $j$
        \item any subset of $M$ summing to $j$ is a row or column
        \item all entries of $M$ are distinct.
    \end{enumerate}
    The existence of such a $j$ and $M$ for every $n$ follows from \cite[Lemma 2.8]{colmenarejo2020toric} (multiplying the $M$ constructed there by $7$ and replacing $j$ by $7j$ if needed). Now replace every entry in $M$ by a generator of $\ZmC$ of the corresponding weight.
    This is possible because $\lambda_{\letter{23},j}\geq1$ for $j\geq 7$.
    Then the rows represent monomials $m_1,\ldots, m_{n+1}$ of weight $j$ and the columns represent monomials $m'_1,\ldots,m'_{n+1}$ of weight $j$, and all these monomials are distinct. Thus
    \[f:= y_{m_1} \cdot \ldots \cdot y_{m_{n+1}} -y_{m'_1} \cdot \ldots \cdot y_{m'_{n+1}} \]
    is an irreducible relation of degree $n+1$. Assume that it is generated in degree $\leq n$. Since $I_{\ZmC,j}$ is a toric ideal, its homogeneous components are spanned (as a vector space) by binomials. Thus, we can write $f = \ell_1 (g_1 - h_1) + \ldots + \ell_s (g_s - h_s)$ where $\ell_1,\ldots,\ell_s$ are variables and $g_1,\ldots,g_s,h_1,\ldots,h_s$ are monomials. Assume wlog.\ that $y_{m_1} \cdot \ldots \cdot y_{m_{n+1}} = \ell_1 g_1$. Since $g_1 - h_1 \in I_{\ZmC,j}$ we have $\phi_{\ZmC,j}(\ell_1 g_1) = \phi_{\ZmC,j}(\ell_1 h_1)$. But this means that $\ell_1 h_1$ corresponds to a partition of the entries of $M$ into $n+1$ blocks that sum to $j$, and because $M$ is rigid, this forces $\ell_1 h_1 = y_{m'_1} \cdot \ldots \cdot y_{m'_{n+1}}$. It follows that $\ell_1(g_1 - h_1) = f$ which is impossible since $f$ is irreducible. So $f\in I_{\ZmC,j}$ is not generated in degree $\leq n$.
\end{proof}

\section{A conjecture on formal MZVs}\label{sec:conj}

As MZVs are nested sums, they satisfy \emph{stuffle relations} (for $k_1,k_2\geq2$)
\begin{align*}
\zeta(k_1)\zeta(k_2)=\zeta(k_1,k_2)+\zeta(k_2,k_1)+\zeta(k_1+k_2).
\end{align*}
On the other hand, MZVs also possess an iterated integral representation due to Kontsevich \cite[Preface]{mzvfromnumberstomotives} \cite{zagier1994}
\begin{align}\label{eqn:zetaasi}
\zeta(k_1,\ldots,k_r) =
{(-1)^r I(0; \underbrace{1,0,\ldots,0}_{k_1},\underbrace{1,0,\ldots,0}_{k_2},\ldots,\underbrace{1,0,\ldots,0}_{k_r};1)} 
\end{align}
where
\begin{equation}\label{eqn:idef}
    I(0; x_1,\ldots,x_n;1) = \int_{0 \leq t_1 \leq \cdots \leq t_k \leq 1} \frac{dt_1}{t_1 - x_1} \cdots \frac{dt_k}{t_k - x_k} \,,
\end{equation}
is the standard simplex integral over the family of differential forms \( \frac{dt_i}{t_i - x_i} \) along the straight line path from 0 to 1. Therefore, they also satisfy \emph{shuffle relations} of the form
\begin{align*}
\zeta(k_1)\zeta(k_2)=\sum_{j=1}^{k_1+k_2-2} \left(\binom{j-1}{k_1-1}+\binom{j-1}{k_2-1}\right) \zeta(j,k_1+k_2-j).
\end{align*}
Comparing the shuffle and stuffle relations of MZVs, we get the \emph{double shuffle relations}. To handle non-convergent MZVs $\zeta(k_1,\ldots,k_r)$ where $k_r=1$, formal regularization processes can be applied \cite{Ihara_Kaneko_Zagier_2006}. These processes assign a well-defined value in $\Z$ to each divergent series such that the stuffle and shuffle relations, respectively, remain satisfied. Comparing the product formulas for the shuffle- and stuffle-regularized MZVs yields the \emph{extended double shuffle relations}. By \cite[Conjecture 1]{Ihara_Kaneko_Zagier_2006}, any algebraic relations in the algebra $\Z$ should be a consequence of the extended double shuffle relations.

\begin{example} \label{ex:double_shuffle}
Applying the stuffle and shuffle product formula from above gives
\begin{align*}
\zeta(2)\zeta(3)&=\zeta(3,2)+\zeta(2,3)+\zeta(5)=\zeta(3,2)+3\zeta(2,3)+6\zeta(1,4),
\end{align*}
and we deduce
\[
\zeta(5)=2\zeta(2,3)+6\zeta(1,4).
\]
But for example, the relation $\zeta(3)=\zeta(1,2)$ from Example \ref{ex:I_3} can only be obtained from considering regularized products, as there is no convergent product in weight $3$.
\end{example}

We give an algebraic formulation of the extended double shuffle relations, this will motivate the definition of formal multiple zeta values.
We start by introducing quasi-shuffle Hopf algebras in general.

Let $A$ be any alphabet, and let $\diamond$ be a commutative and associative product on the vector space $\Q A$ spanned by the elements in $A$. The associated \emph{quasi-shuffle product} $\ast_{\diamond}$ to $\diamond$ is recursively defined on the free algebra $\Q\langle A\rangle$ by 
\begin{align*}
v\ast_\diamond 1 &= 1\ast_\diamond v =v, \\
a_1u\ast_\diamond a_2v &= a_1(u\ast_\diamond a_2v)+a_2(a_1u\ast_\diamond v)+ (a_1\diamond a_2)(u\ast_\diamond v),
\end{align*}
for all $a_1,a_2\in A$, and $u,v\in \Q\langle A\rangle$. Furthermore, let $\Delta:\Q\langle A\rangle \to \Q\langle A\rangle \otimes \Q\langle A\rangle$ be the \emph{deconcatenation coproduct}, given by
\[
\Delta(w)=\sum_{uv=w} u\otimes v
\]
for $w\in \Q\langle A\rangle$.

\begin{proposition}[{\cite{Hoffman}}]\label{prop:quasi-sh-hopf}\begin{enumerate}
    \item The tuple $(\Q\langle A\rangle,\ast_\diamond,\Delta)$ is a Hopf algebra.
    \item The pair $(\Q\langle A\rangle,\ast_\diamond)$ is a free algebra generated by the Lyndon words over $A$.
\end{enumerate}
\end{proposition}

\noindent Recall the definition of Lyndon words from Definition \ref{def:lyndon}. We abbreviate $\Q\langle A\rangle_{*_\diamond}=(\Q\langle A\rangle,\ast_\diamond,\Delta)$.
Gradings of $\Q\langle A\rangle_{\ast_\diamond}$ are one-to-one with gradings of the non-unital algebra $(\Q A,\diamond)$. Proposition \ref{prop:quasi-sh-hopf} together with Lemma \ref{lem:free_alg_toric} has the following important consequence.

\begin{proposition}
Assume $(\Q A,\diamond)$ is a graded non-unital algebra and consider $\Q\langle A\rangle_{\ast_\diamond}$ with the induced grading. Then for every $j$, $\relvar_j(\Q\langle A\rangle_{\ast_\diamond})$ is a projective toric variety.
\end{proposition}

If the diamond $\diamond$ is zero everywhere,
\begin{align} \label{eq:shuffle}
a_1\diamond a_2=0,\qquad a_1,a_2\in A,
\end{align}
we call the corresponding quasi-shuffle product the \emph{shuffle product} and denote it by $\shuffle$.

\begin{remark}
    If $A = \{1,\ldots,d\}$ and we endow $\Q A$ with the trivial grading concentrated in degree $1$, the variety $\relvar_j(\Q\langle A\rangle_\shuffle)$ is the \textit{universal variety} $\mathcal{U}_{d,j}$ from \cite{bib:AFS2019, colmenarejo2020toric}.
\end{remark}

Given an alphabet $A$, it follows from Proposition \ref{prop:quasi-sh-hopf} that the algebras $\mathbb Q\langle A\rangle_{\shuffle}$ and $\mathbb Q\langle A\rangle_{\ast_\diamond}$ are always isomorphic. In fact, there is a natural choice for this isomorphism, known as the \textit{Hoffman exponential} \cite{Hoffman}. If $\mathbb Q\langle A\rangle_{\ast_\diamond}$ is graded, this is a graded isomorphism.

\begin{corollary} \label{cor:quasi_sh_var_iso}
Fix an alphabet $A$ and $j\in \mathbb N$. Then all the varieties $\relvar_j(\Q\langle A\rangle_{\ast_\diamond})$ associated to graded non-unital algebras $(\mathbb Q A, \diamond)$ are isomorphic.
\end{corollary}

Now we introduce the two specific quasi-shuffle products, which are needed to formulate the extended double shuffle relations from the above algebraically. The first one is the \emph{shuffle product} $\shuffle$ defined on the specific alphabet $X=\{\letter 0, \letter 1\}$. The Hopf algebra $\Q\langle X\rangle_{\shuffle}$ is graded for any grading of $\Q X$, in the following we choose the one where every letter has degree $1$. We write 
\[
S^\shuffle:=\Q\langle X\rangle_\shuffle /\langle \letter 0 \rangle.
\]
Note that $S^\shuffle$ can be identified with the space spanned by words which start with $\letter 1$.

The second one in the \emph{stuffle product} $\stuffle$ defined on the specific alphabet $A=\mathbb{N}$, which is associated to the diamond
\begin{align} \label{eq:stuffle}
\letter i \diamond \letter j = \letter i \letter + \letter j.
\end{align}
If we grade $\Q A$ by $(\Q A)_i = \Q \letter i$, then $\Q\langle \mathbb{N}\rangle_\stuffle$ is a graded Hopf algebra. Let
\[
S^\stuffle=\Q\langle \mathbb{N}\rangle_\stuffle.
\]

The algebras $S^\stuffle$ and $S^\shuffle$ describe the two product expressions of MZVs, coming from the iterated sum representation in the definition and the iterated integral representation in \eqref{eqn:zetaasi}. To make this precise, we restrict to admissible indices, i.e., indices where the corresponding MZVs converge. So let $S^{\stuffle,\operatorname{adm}}$ be the subalgebra of $S^\stuffle$ spanned by words which do not end in $\letter 1$. Similarly, we let $S^{\shuffle,\operatorname{adm}}$ be the subalgebra of $S^\shuffle$ spanned by words which do not end in $\letter 1$.
\begin{proposition}
\label{prop:zetastufflemorphism}
\begin{enumerate}
    \item There is a surjective algebra morphism
    \begin{align} \label{eq:stuffle_mzv}
    \zeta^\stuffle: S^{\stuffle,\operatorname{adm}}\to \mathcal{Z},  \quad \letter{k}_1 \cdots \letter{k}_r \mapsto \zeta(k_1,\ldots,k_r).
    \end{align}
    \item There is a surjective algebra morphism
    \begin{align} \label{eq:shuffle_mzv}
    \zeta^\shuffle: S^{\shuffle,\operatorname{adm}}\to \mathcal{Z},  \quad \letter 1  \letter 0^{k_1 - 1} \cdots \letter 1 \letter 0^{k_r - 1} \mapsto \zeta(k_1,\ldots,k_r).
    \end{align}
\end{enumerate}
\end{proposition}

\begin{proof} See \cite[Proposition~1.148, 1.156]{mzvfromnumberstomotives}.
\end{proof}

\begin{remark*}
    Note that here $\zeta^\stuffle$ and $\zeta^\shuffle$ do not denote the \emph{regularized MZVs}, which is sometimes common in the MZV literature. Instead, we will refer to these regularized maps as $\zetastreg$ and $\zetashreg$, see Proposition \ref{prop:zeta-reg} below.
\end{remark*}

Consider the graded vector space isomorphism
\begin{align*}
    \gamma: S^\stuffle &\to S^\shuffle,\\
    \letter{i}_1\cdots \letter{i}_k &\mapsto \letter 1  \letter 0^{i_1 - 1} \cdots \letter 1 \letter 0^{i_k - 1}.
\end{align*}
Then we have $\zeta^\stuffle = \zeta^\shuffle \circ \gamma$. Note however that $\gamma$ is not a morphism of algebras. By comparing the stuffle products coming from \eqref{eq:stuffle_mzv} and the shuffle products coming from \eqref{eq:shuffle_mzv}, we get the \emph{finite double shuffle relations} for MZVs. Those are not enough to cover all relation among MZVs as indicated in Example \ref{ex:double_shuffle}, hence we introduce \emph{regularization}. This means, we extend the maps in \eqref{eq:stuffle_mzv} resp. \eqref{eq:shuffle_mzv} to the whole algebras $S^\stuffle$ resp. $S^\shuffle$ such that the stuffle resp. shuffle product is still satisfied. One verifies as in \cite[Theorem~1.168, 1.170]{mzvfromnumberstomotives} that this can be done for any choice of the image of $\letter 1$, hence we choose this image to be $0$.

\begin{proposition}\label{prop:zeta-reg}
    \begin{enumerate}
    \item Define an algebra homomorphism 
    \[\zetastreg: S^\stuffle \to \Z\] 
    on Lyndon words by $\letter 1 \mapsto 0$, $\letter 1 \not= w \mapsto \zeta^\stuffle(w)$. This is well-defined if the alphabet $\mathbb N$ is ordered by $\letter 1 < \letter 2 < \ldots$. Then for any $w\in S^{\stuffle,\operatorname{adm}}$ we have $\zetastreg(w) = \zeta^\stuffle(w)$.
    \item  Define an algebra homomorphism 
    \[\zetashreg: S^\shuffle \to \Z\] 
    on Lyndon words by $\letter 1 \mapsto 0$, $\letter 1 \not= w \mapsto \zeta^{\shuffle}(w)$. This is well-defined if the alphabet $\{\letter 0,\letter 1\}$ is ordered by $\letter 1 < \letter 0$. For any $w\in S^{\shuffle,\operatorname{adm}}$ we have $\zetashreg(w) = \zeta^\shuffle(w)$. 
\end{enumerate}
\end{proposition}

The maps $\zetastreg$ and $\zetashreg$ are called the \textit{stuffle} and the \textit{shuffle} regularization of MZVs, respectively. 
  
\begin{proof}
By Proposition \ref{prop:quasi-sh-hopf} this is well-defined in both cases. By Proposition \ref{prop:zetastufflemorphism}, it suffices to show that the subrings $S^{\stuffle,\operatorname{adm}}$ resp.\ $S^{\shuffle,\operatorname{adm}}$ are generated by all Lyndon words except $\letter 1$ (in $S^{\stuffle}$ resp.\ $S^{\shuffle}$) (cf. \cite[Theorem~3.1, Theorem~4.1]{Hof97}). We demonstrate the argument for $S^{\stuffle,\operatorname{adm}}$. Clearly, the subring generated by the Lyndon words is contained in $S^{\stuffle,\operatorname{adm}}$. Conversely, write $p \in S^{\stuffle,\operatorname{adm}}$ as a polynomial in Lyndon words and assume that it involves the variable $\letter 1$. Then consider it as a polynomial in $\letter 1$ with coefficients in polynomials of the remaining Lyndon words. Since by the choice of ordering on the alphabet no other Lyndon word ends on $\letter 1$, this leading term can not simplify to a word not ending on $\letter 1$. For the same reason, it can not simplify with any polynomials of lower degree in $\letter 1$, contradicting the choice of $p$.
\end{proof}

Now by definition, the map $\zetastreg$ factors over the algebra $\algstreg := S^\stuffle / \langle \letter 1 \rangle$ and the map $\zetashreg$ factors over $\algshreg := S^\shuffle / \langle \letter 1 \rangle$. We call them the \textit{regularized stuffle} and \textit{regularized shuffle algebra}, respectively.
Any element in $S^\stuffle$ can be written as a polynomial in $\letter 1$ with coefficients in $S^{\stuffle,\operatorname{adm}}$ and any element in $S^\shuffle$ can be written as a polynomial in $\letter 1$ with coefficients in $S^{\shuffle, \operatorname{adm}}$, see \cite[Theorem~1.168, 1.170]{mzvfromnumberstomotives}. Hence, we get isomorphisms 
\[
\algstreg\simeq S^{\stuffle,\operatorname{adm}},\qquad \algshreg\simeq S^{\shuffle,\operatorname{adm}}.
\]

Importantly, shuffle and stuffle regularization are not naively compatible: $\zetastreg \neq \zetashreg \circ \gamma$. One needs to replace $\gamma$ by another comparison map, which was introduced in \cite{Ihara_Kaneko_Zagier_2006} and which we recall now in our algebraic setup.

Let $\lambda_k \in \mathbb Q\langle \letter 0,\letter 1\rangle$ be the coefficient of $T^k$ in the power series
\begin{align} \label{eq:def_lambda}
\exp_\shuffle\left(\sum_{k\geq 2} \frac{(-1)^{k-1}}{k} \letter 1 \letter 0^{k-1} \cdot T^k \right) \in \mathbb Q_\shuffle\langle \letter 0, \letter 1\rangle[\![ T ]\!]  
\end{align}
where $\exp_\shuffle$ denotes the shuffle exponential.
\begin{defi} \label{def:rho}
    The \textit{comparison isomorphism} is the morphism
    \begin{align*}
         \rho: S^{\stuffle} &\to S^{\shuffle} \\
         w &\mapsto \sum_{w =v \letter 1^k} \lambda_k \shuffle \gamma(v)
    \end{align*}
    of graded vector spaces.
\end{defi}

We emphasize that $\rho$ is not an algebra homomorphism. In fact, the failure of $\rho$ to be an algebra homomorphism is exactly what gives the extended double shuffle relations.

\begin{lemma}
\label{lem:rhoisom}
    The map $\rho$ is indeed an isomorphism of $\mathbb Q$-vector spaces. Its inverse is given by 
    \begin{align*}
    \rho^{-1}: S^{\shuffle} &\to S^{\stuffle}, \\
    w&\mapsto \sum_{w=v \letter 1^k} \gamma^{-1}(\mu_k\shuffle v).
    \end{align*}
     where $\mu_k\in \Q\langle \letter 0, \letter 1\rangle$ is the coefficient of $T^k$ in
\begin{align*}
\exp_\shuffle\left(\sum_{k\geq2} \frac{(-1)^k}{k}\letter 1 \letter 0^{k-1} \cdot T^k\right) \in \Q_\shuffle\langle \letter 0, \letter 1\rangle[\![ T ]\!].
\end{align*}
\end{lemma}
\begin{proof}
See Appendix \ref{app:proofrho}.
\end{proof}

\begin{proposition}\label{prop:reg-map}
    We have
    \[\zetastreg = \zetashreg \circ \rho.\]
\end{proposition}
\begin{proof}
See \cite[Theorem~1]{Ihara_Kaneko_Zagier_2006} and \cite{BoutetDeMonvel2004DivergentMultizeta}.
\end{proof}

Comparing stuffle and shuffle products of MZVs with the help of Proposition \ref{prop:reg-map} gives the \emph{extended double shuffle relations} among MZVs.

\begin{example}
    By Proposition \ref{prop:reg-map}, we have
    \begin{align*}
    0 &= \zetastreg(\letter 1) \cdot \zetastreg(\letter 2) = \zetastreg(\letter 1 \stuffle \letter 2) = \zetastreg(\letter{12} + \letter{21} + \letter{3}) = \zetashreg(\letter{110} + \letter{101} + \letter{100})\\
    &= \zetashreg(\letter{110} + (\letter{10} \shuffle \letter{1} - 2 \cdot \letter{110}) + \letter{100}) = \zetashreg(-\letter{110} + \letter{100}) = -\zeta(1,2) + \zeta(3).
    \end{align*}
    which is Euler's identity. The following example is a bit more involved: we have
    \[\zetastreg(\letter{211}) = \zetashreg(\letter{1011} - \frac{1}{2} \letter{10} \shuffle \letter{10}) = \zetashreg(3 \cdot \letter{1110} - \frac{1}{2}\letter{10} \shuffle \letter{10}) = 3 \zeta(1,1,2) - \frac{1}{2} \zeta(2)^2\]
    where we used the comparison map in the first equality and shuffle regularization in the second, while on the other hand
    \[\zetastreg(\letter{211}) = \zetastreg(\letter{112} + \letter{13} + \frac{1}{2} \letter{4}) = \zeta(1,1,2) + \zeta(1,3) + \frac{1}{2}\zeta(4).\]    
\end{example}

It is conjectured in \cite{Ihara_Kaneko_Zagier_2006} that all algebraic relations in the algebra $\mathcal{Z}$ of MZVs come from the extended double shuffle relations. This directly motivates the definition of formal MZVs.

\begin{defi}
The algebra of \emph{formal MZVs} is given by
\[
\mathcal{Z}^{\mathfrak f} := S^\stuffle / \langle \letter 1, I_{\mathrm{EDS}}\rangle \cong \algstreg / \overline{I_{\mathrm{EDS}}}.
\]
where $I_{\mathrm{EDS}} := \langle v \stuffle w - \rho^{-1}(\rho(v) \shuffle \rho(w))\mid v,w\in \Q\langle \mathbb{N}\rangle\rangle$ is the ideal of extended double shuffle relations and $\overline{I_{\mathrm{EDS}}}$ its image in $\algstreg$.
\end{defi}

We do not need all pairs $v,w \in \Q\langle \mathbb{N}\rangle$ to generate $I_{\mathrm{EDS}}$. Several smaller sets generating this ideal are given in \cite[Theorem~2]{Ihara_Kaneko_Zagier_2006}. We propose a new conjecture on a small generating set for this ideal in Conjecture \ref{conj:lin-form-approx-eq}.

Note that by definition of the ideal $I_{\mathrm{EDS}}$, the datum of an algebra homomorphism $\formZ \to R$ for any $\mathbb Q$-algebra $R$ is precisely the datum of a map $\alpha: S^{\stuffle} \to R$ that maps $\letter 1$ to $0$ such that $\alpha\circ \rho^{-1}: S^\shuffle \to R$ is an algebra homomorphism as well. In other words, morphisms $S^\stuffle \to R$ that factor through $ \formZ$ behave algebraically like MZVs, see Proposition \ref{prop:reg-map}. 

Also note that by definition of $\formZ$, the map $S^\shuffle \to \formZ, \ w \mapsto \overline{\rho^{-1}(w)}$ is a surjective algebra homomorphism.

As motivic multiple zeta values satisfy the extended double shuffle relations, there are surjective algebra morphisms
\begin{align*}
\formZ\to \motZ \to \Z,\quad \zeta^{\mathfrak f}(k_1,\ldots,k_r)\mapsto \zeta^{\mathfrak m}(k_1,\ldots,k_r) \mapsto \zeta(k_1,\ldots,k_r)
\end{align*}
where $\zeta^{\mathfrak f}(k_1,\ldots,k_r)$ denotes the class of $\letter k_1 \ldots \letter k_r$ in $\formZ$.
So Lemma \ref{lem:inner-ideal-approx} implies $I_{\formZ,j}\subseteq I_{\motZ,j} \subseteq I_j$. By \cite[Conjecture 1]{Ihara_Kaneko_Zagier_2006} we should have $\formZ\simeq\Z$ and thus $I_{\formZ,j}=I_j$.
Therefore, we would like to study the relation variety for $\mathcal{Z}^{\mathfrak f}$. However, for this we would need to know at least the dimensions of the graded components $\mathcal{Z}_j^{\mathfrak f}$, which are the subject of an open conjecture \cite[Section~5.6]{BCK2024FormalMZV}. Thus, we take another approach: we study the relation varieties of $\algstreg$ and $\algshreg$ separately and then compute their intersection in the appropriate ambient space. 

To make this precise, consider the following commutative diagram:
\begin{equation}\label{eq:sh-st-diagram}
    \begin{tikzcd}
	{\Sym S_j^{\stuffle}} & {\Sym \algstreg_j} & {\algstreg} & \\
	&&& \formZ \\
	{\Sym S_j^{\shuffle}} & {\Sym \algshreg_j} & {\algshreg}
	\arrow[two heads, from=1-1, to=1-2]
	\arrow["\hspace{2.5em} \cong \ \Sym \rho_j"{description}, from=1-1, to=3-1]
	\arrow[from=1-2, to=1-3]
	\arrow[two heads, from=1-3, to=2-4]
	\arrow[two heads, from=3-1, to=3-2]
	\arrow[from=3-2, to=3-3]
	\arrow[two heads, from=3-3, to=2-4]
\end{tikzcd}
\end{equation}
Let us write $\psi_\ast$ for the map $\Sym S_j^\stuffle \to \algstreg$ and $\psi_\shuffle$ for the map $\Sym S_j^\shuffle \to \algshreg$ in the diagram. Then the kernel of the map $\Sym S_j^\stuffle \to \formZ$ contains both $\ker(\psi_\ast)$ and $(\Sym \rho_j^{-1})(\ker \psi_\shuffle)$. Moreover, we have $\ker \psi_\ast = J_{\letter 1, \stuffle, j} + I_{\algstreg,j}$ and $\ker \psi_\shuffle = J_{\letter 1, \shuffle, j} + I_{\algshreg,j}$ where 
\[J_{\letter 1,\stuffle, j} := \langle(\letter 1\stuffle S^\stuffle)_j\rangle \quad \text{and} \quad J_{\letter 1,\shuffle, j} := \langle(\letter 1\shuffle S^\shuffle)_j\rangle \]
are the ideals generated by the linear forms given by the weight $j$ stuffle and shuffle multiples of $\letter 1$, respectively.

Thus, if $p$ denotes the projection $\Sym S_j^\stuffle \to \Sym \formZ_j$, we obtain
\begin{equation}\label{eq:inner-form-approx}
    J_{\letter 1,\stuffle,j} + I_{\algstreg,j} + (\Sym \rho_j^{-1})(J_{\letter 1,\shuffle,j} + I_{\algshreg,j})  \subseteq p^{-1}(I_{\formZ,j})
\end{equation}

since $\Sym S_j^\stuffle \to \formZ$ agrees with $\phi_{\formZ,j} \circ p$. This inspires the following conjecture:

\begin{conjecture}\label{conj:inner-form-approx-eq}
    The inclusion in \eqref{eq:inner-form-approx} 
    is an equality.
\end{conjecture}

Because on the left hand side of \eqref{eq:inner-form-approx} $J_{\letter 1,\stuffle,j}$ and $J_{\letter 1,\shuffle,j}$ are generated by linear forms, while $I_{\algstreg,j}$ and $(\Sym \rho_j^{-1})(I_{\algshreg,j})$ are generated in degree $\geq 2$, the equality in \eqref{eq:inner-form-approx} would imply the following:
\begin{conjecture}\label{conj:lin-form-approx-eq}
We have
    \[S^\stuffle_j / \big((\letter 1 \stuffle S^\stuffle)_j +  \rho_j^{-1}((\letter 1 \shuffle S^\shuffle)_j))\big) \simeq \formZ_j \]
for all $j$. In other words,
    \[S^\stuffle/(\langle \letter 1 \rangle_\stuffle + \rho^{-1}(\langle \letter 1 \rangle_\shuffle)) \simeq \formZ.\]
    Note that this implies in particular that
    \[I_{\mathrm{EDS}} = \langle \letter 1 \rangle_\stuffle + \rho^{-1}(\langle \letter 1 \rangle_\shuffle).\]
\end{conjecture}
We have verified Conjecture \ref{conj:lin-form-approx-eq} computationally up to weight $j=23$. To our knowledge this generating set for the EDS-relation has not been suggested in the literature so far.

\begin{proposition}
    Conjecture \ref{conj:inner-form-approx-eq} and Conjecture \ref{conj:lin-form-approx-eq} are equivalent.
\end{proposition}
\begin{proof}
    As discussed above, Conjecture \ref{conj:inner-form-approx-eq} implies Conjecture \ref{conj:lin-form-approx-eq}. For the converse, note that assuming Conjecture \ref{conj:lin-form-approx-eq} we have
    \[\formZ = \algstreg\Big/ \ \overline{\rho^{-1}(\langle \letter 1\rangle_\shuffle)}.\]
    Set $I' := J_{\letter 1,\stuffle,j} + I_{\algstreg,j} + (\Sym \rho_j^{-1})(J_{\letter 1,\shuffle} + I_{\algshreg,j})$. Then by definition
    \[R := (\Sym S_j^\stuffle) / I' \cong \algstreg_{\langle j \rangle} \ \Big/ \ \overline{(\Sym \rho_j^{-1})(J_{\letter 1,\shuffle,j} + I_{\algshreg,j})}\]
    where $\algstreg_{\langle j\rangle}$ denotes the subring generated in degree $j$. The surjection $\Sym S^\stuffle \to \formZ$ factors through $R$; we need to show that the map $R \to \formZ$ is an isomorphism. For this it suffices to show that the kernel of $\algstreg_{\langle j\rangle} \to R$ contains the restriction of $\overline{\rho^{-1}(\langle \letter 1\rangle_\shuffle)}$ to $\algstreg_{\langle j\rangle}$. Pick $x \in \algstreg_{\langle j\rangle}$ from this restricted ideal, and let $x'$ denote a lift to $\Sym S^\stuffle_j$. Then by assumption $(\Sym \rho_j)(x')$ maps to an element of $\langle \letter 1\rangle_\shuffle$ under $\phi_{S^\shuffle,j}$. But this means $(\Sym \rho_j)(x')$ is in the kernel of the map $\Sym S_j^\shuffle \to \algshreg$, which is $J_{\letter 1,\shuffle,j} + I_{\algshreg,j}$.
\end{proof}

Motivated by Conjecture \ref{conj:inner-form-approx-eq} we now want to study the relation varieties $\relvar_j(\algstreg)$ and $\relvar_j(\algshreg)$. The following fact is convenient:

\begin{proposition}\label{prop:ssh-sst-iso}
    There is an isomorphism $\psi: S^\shuffle \to S^\stuffle$ of graded algebras, inducing an isomorphism $\algshreg \to \algstreg$.
\end{proposition}
\begin{proof}
By Proposition \ref{prop:quasi-sh-hopf}, $\Q\langle \letter 0, \letter 1 \rangle_\shuffle$ is freely generated by the Lyndon words in $\letter 0, \letter 1$ with respect to the ordering $\letter 1 < \letter 0$. Hence, $S^\shuffle$ is freely generated by the same Lyndon words expect the word $\letter 0$. Similarly, $S^\stuffle$ is freely generated by the Lyndon words in $\mathbb N$ with respect to the reversed ordering $\letter 1 < \letter 2 < \letter 3 < \cdots$. So it suffices to give a bijection between the Lyndon words in $\letter 0, \letter 1$ (except $\letter 0$) in weight $j$ and the Lyndon words in $\mathbb N$ of weight $j$.

As there are the same number of Lyndon words of weight $j\geq2$ in $\mathbb{N}$ and $\letter 0, \letter 1$, it suffices to check that $\gamma$ maps any Lyndon word in $\mathbb{N}$ to a Lyndon word in $\letter 0, \letter 1$ under the map $\rho$ from Definition \ref{def:rho}. So, let $\letter{i}_1\cdots \letter{i}_r \in \mathbb Q\langle \mathbb{N}\rangle$ be a Lyndon word.

Assume there is a right subword $v$ of the image $w := \letter 1 \letter 0^{i_1-1} \letter 1 \cdots \letter 1\letter 0^{i_r-1}$ with $v\leq w$ lexicographically. Then this subword has to start with $\letter 1$ by our choice of order. Thus there is some $i_j$ with $v = \letter 1 \letter 0^{i_j -1} \letter 1 \cdots \letter 1 \letter 0^{i_r - 1}$. But then again by the choice of order $v\leq w$ implies $\letter {i}_j \cdots \letter {i}_r \leq \letter {i}_1 \cdots \letter {i}_r$ lexicographically which is impossible since $\letter {i}_1 \cdots \letter {i}_r$ is a Lyndon word.\qedhere

\end{proof}

\noindent In particular, $\relvar_j(\algstreg) \cong \relvar_j(\algshreg)$ under the linear isomorphism $\mathbb P((\algstreg_j)^\vee) \cong \mathbb P((\algshreg_j)^\vee)$ induced by $\psi$. Thus it suffices to study the relation varieties associated to $\algstreg$, which will be the content of the next section.

\section{The regularized stuffle variety}\label{sec:stuf-var}

Motivated by Section \ref{sec:conj} we now study the relation varieties $\relvar_j(\algstreg)$.

\begin{notdef}
    Similar as before, we write $\mathcal L(\mathbb N)$ for the set of Lyndon words over the alphabet $\mathbb N$, ordered by $\letter 1 < \letter 2 < \ldots$. Write $\lyndonmone(\mathbb N) := \lyndon(\mathbb N) - \{ \letter 1\}$. We identify $\StC \cong \mathbb C[y_w \ | \ w \in \mathcal L_{\not= \letter 1}(\mathbb N)\}$ and write $\monlyndonmone(\mathbb N)$ for the set of monomials formed out of the variables $\{x_w \ | \ w\in \mathcal L_{\not= \letter 1}(\mathbb N)\}$. Note that the degree (or weight) of a variable $x_{\letter i_1 \ldots \letter i_k}$ is $i_1 + \dots + i_k$.\par
    We write
    \begin{equation*}
    \lyndonmone(\mathbb N)_{\leq k}:=\{w\in \lyndonmone(\mathbb N) \ | \ \weight(w)\leq k\}
    \end{equation*}
    and
    \begin{equation*}
        \monlyndonmone(\mathbb N)_{k}:=\{m\in \monlyndonmone(\mathbb N) \ | \ \deg(m)= k\}.
    \end{equation*}
\end{notdef}

We have $\Sym(\StCmone_j) \cong \mathbb C[y_m \ | \ m \in \monlyndonmone(\mathbb N)_j]$ and defining $\StCmone_{\langle \leq j \rangle}$ as the subalgebra of $\StCmone$ generated in degree $\leq j$, we see that $\relvar_j(\StCmone)$ is cut out in $\mathbb P(\StCmone_j)$ by the homogeneous ideal
$$I_{\StCmone,j} = \ker(\phi_{\StCmone,j} : \mathbb C[y_m \ | \ m \in \monlyndon(\mathbb N)_j] \to \StC_{\langle \leq j \rangle}).$$

\begin{lemma}\label{lem:stuffle-amb-dim}
    The dimension of $\mathbb P(\StCmone_j)$ is $2^{j-2} - 1$.
\end{lemma}
\begin{proof}
    Since $\StCmone$ is spanned as a $\mathbb C$-vector space by words not ending on $\letter 1$, the dimension of $\StCmone_j$ is the number of compositions of $j$ (without zeroes) minus the number of compositions of $j-1$. Since compositions of $j$ are in bijection with subsets of $\{1,\ldots,j-1\}$, there are $2^{j-1} - 2^{j-2}$ such compositions.
\end{proof}

Next, let us again determine dimension and (projective) normality of the varieties $\relvar_j(\StCmone)$. They are again cones over a smaller base variety, which we denote $\relvarr_j(\StCmone)$.

Let $\lambda_{\mathbb N,k}$ denote the cardinality of $\lyndon(\mathbb N)_k$. 

\begin{lemma}
    We have
\begin{equation*}
    \lambda_{\mathbb N,k} = \frac{1}{k}\sum_{i|k}\mu(k/i)(2^i-1).
\end{equation*}
\end{lemma}
\begin{proof}
    This follows by a M\"obius inversion argument as in \cite[Corollary~4.14]{Reutenauer}. More concretely, the proof of Lemma \ref{lem:stuffle-amb-dim} implies
    \[\frac{1}{2} \left(1 + \frac{1}{1-2t}\right) = 1 + \sum_{i=1}^\infty 2^{k-1} t^k = \prod_{i=1}^\infty \frac{1}{(1-t^i)^{\lambda_{\mathbb N, i}}}\]
    since $\mathbb Q\langle \mathbb N\rangle$ is isomorphic to $\mathbb Q[x_w\mid w\in\mathcal{L}(\mathbb{N})]$ as a graded vector space by Proposition \ref{prop:quasi-sh-hopf}, and the two sides of the equation are the respective Hilbert series.
    Taking logarithmic derivatives on both sides and multiplying by $t$ one obtains
    \[\frac{t}{(1-t)(1-2t)} = \sum_{i=0}^\infty i \lambda_{\mathbb N,i} \frac{t^i}{1-t^i}.\]
    Comparing coefficients, we see that \(\displaystyle 2^i - 1 = \sum_{\ell=0}^{i-1} 2^\ell = \sum_{\ell | i} \ell \lambda_{\mathbb N, \ell} \)
    and thus by M\"obius inversion
    \[\lambda_{\mathbb N,i} = \frac{1}{i} \sum_{\ell | i} \mu(i/\ell) (2^\ell - 1). \qedhere\]
\end{proof}

Similar as before, let $A_{\StCmone,j}$ denote the linear map $\mathbb Z^{\monlyndonmone(\mathbb N)_j} \to \mathbb Z^{\lyndonmone(\mathbb N)_{\leq j}}$ induced by mapping a monomial $m$ to its exponent vector. Let us write 
\[M_j := \{A_{\StC,j}(e_m) \ | \ m\in \monlyndonmone(\mathbb N)_j\} =\Big\{\sum_i a_i w_i\,|\,a_i\in\mathbb N_0,\,w_i\in\lyndonmone (\mathbb N)_{\leq j},\,\sum_i a_i\weight(w_i)=j\Big\}.\]
We denote the polytope given by the convex hull of $M_j$ by $P_{\StCmone,j}$. Again, we study the variety $\relvar_j(\StCmone)$ via Theorem \ref{thm:toric-recall}.

\begin{theorem}\label{thm:unregstuffalg}
    \begin{enumerate}
        \item For $j\geq 6$ even or $j\geq 9$ odd, $M_j\cap \mathbb Z^{\lyndonmone(\mathbb N)_{\leq j-2}}$ generates the sublattice 
        \[ L_j=\Big\{\sum_i a_i w_i\,|\,a_i\in\mathbb Z,\,w_i\in\lyndonmone (\mathbb N)_{\leq j-2},\,\sum_i a_i\weight(w_i)\in j\mathbb Z\Big\}\] of index $j$ in $\mathbb Z^{\mathcal L(\mathbb N)_{\leq j-2}}$.
        \item For $j\geq 11$, we have $
            \relvarr_j(\StCmone)=\relvar_j(\StCmone_{\langle \leq j-3\rangle})$
        \item For $j\geq 8$, we have
    \(\displaystyle\dim \relvar_j(\StC) = \lambda_{\mathbb N,j}+\sum_{i=2}^{j-2} \lambda_{\mathbb N, i} - 1\) and \(\displaystyle\dim \relvarr_j(\StC) = \sum_{i=2}^{j-3} \lambda_{\mathbb N, i} - 1.\)
        \item $\relvarr_j(\StCmone)$ is not normal for $j=16, 18$ and $j \geq 20$.
    \item For every $n$ there is a $j$ such that $I_{\StCmone,j}$ is not generated in degree $\leq n$.
    \end{enumerate}
    
\end{theorem}

\begin{proof}
    \begin{enumerate}
        \item The proof that $M_j$ generates $L_j$ is completely analogous to Lemma~\ref{lem:gen-lattice}.

        $L_j$ has index $j$ as by definition, $L_j$ is the kernel of the morphism
        \[\mathbb Z^{\lyndonmone(\mathbb N)_{\leq j-2}} \to \mathbb Z / j\mathbb Z, \ w_i \mapsto \mathrm{wt}(w_i)\!\!\!\!\mod j\]
    which is surjective for $j\geq 1$.
    \item This can be argued as in Proposition \ref{prop:rel0-is-rel}, using that now there are $2$ distinct monomials (both of them generators) of weight $3$.
    \item The first assertion follows from 1.\ and Theorem~\ref{thm:toric-recall}.
    Note that the induction argument used to prove point 1.\ shows also that $M_j\cap \mathbb Z^{\lyndonmone(\mathbb N)_{\leq j-3}}$ generates $L_j\cap \mathbb Z^{\lyndonmone(\mathbb N)_{\leq j-3}}$,
    which then implies the second assertion using point 2. for $j\geq 11$. For $j=8,9,10$ the second assertion can be computationally verified.
    \item Analogous to the proof of Theorem \ref{thm:mot-non-normal}. We just replace $\lyndon(\letter 2,\letter 3)_{\leq j}$ by $\lyndonmone(\mathbb N)_{\leq j}$, $\monlyndon(\letter 2,\letter 3)_{\leq j}$ by $\monlyndonmone(\mathbb N)_{\leq j}$, $A_{\ZmC,j}$ by $A_{\StCmone,j}$ and $\letter{23}$ by $\letter 5$.
    \item Completely analogous to the proof of Proposition \ref{prop:not_gen_in_fixed_deg}. However, since $\StCmone$ has a generator in each weight $\geq 2$,  the entries of the rigid matrices may now be chosen from all positive integers $\geq 2$. \qedhere
    \end{enumerate}
\end{proof}
For a similar theorem concerning the case of the whole stuffle algebra, see Appendix \ref{sec:unregstufthm}.

\begin{table}[!htbp]
    \centering
    \renewcommand{\arraystretch}{1.2}
    \begin{tabular}{c|ccccccccc}
        Property $\backslash$ $j$ & 6 & 7 & 8 & 9 & 10 & 11 & 12 & 13 & 14\\\hline\vspace{-1em} \\ 
        $\mathrm{dim} \ \relvar_j(\StCmone)$ & 14 & 28 & 50 & 94 & 167 & 310 & 558\\
        amb.\ $\dim \ \relvar_j(\StCmone)$ & 15 & 31 & 63 & 127 & 255 & 511 & 1023 & 2047 & 4095\\
        $\mathrm{dim} \ \relvarr_j(\StCmone)$ & 1 & 4 & 11 & 20 & 38 & 68 & 124 \\
        amb.\ $\dim \relvarr_j(\StCmone)$ & 2 & 7 & 24 & 53 & 126 & 269 & 589 & \\
        $\mathrm{deg} \ \relvar_j(\StCmone)$ & 2 & 4 & 64 & 1128 &495616 & \\
        $\relvarr_j(\StCmone)$ is proj.\ normal &\true & \true & \true &\true &\true &\true & \true & \true & \true \\
        $\relvarr_j(\StCmone)$ is normal &\true & \true &\true &\true & \true & \true & \true & \true &\true \\
        $\relvarr_j(\StCmone)$ is smooth &\true & \true & \false & \false & \false & \false \\
        generators of $I_{\StCmone, j}$ & $1_2$ & $6_2$ & $48_2$ & $282_2$ & $1442_2$ & $7533_2$
    \end{tabular}
    \caption{Properties of the projective varieties $\relvar_j(\StCmone)$ and $\relvarr_j(\StCmone)$.
    Note that  $\relvarr_j(\StCmone)$ is also projectively normal (and thus normal) for $j=15,17,19$. The empty cells in the table represent calculations which did not terminate in reasonable time due to computational limits.}
    
    \end{table}

\section{Comparison of formalisms}\label{sec:comp}

\subsection{Connection to the double shuffle group}

Racinet \cite{Racinetthesis} studied the algebra $\formZ$ of formal MZVs by considering the associated affine group scheme $\operatorname{DM}_0$.

\begin{defi} \label{def:DM0}
For any $\mathbb{Q}$-algebra $R$, the set $\operatorname{DM}(R)$ consists of all non-commutative power series $\phi\in R\langle\langle \letter 0,\letter 1\rangle\rangle$ such that
\begin{enumerate}
    \item $(\phi\mid \letter 0)=(\phi\mid \letter 1)=0 $,
    \item $\phi$ is grouplike for the coproduct $\Delta_{\shuffle}$ dual to the shuffle product,
    \item $\phi_\ast$ is grouplike for the coproduct $\Delta_\stuffle$ dual to the stuffle product.
\end{enumerate}
Here, we set 
\begin{align} \label{eq:powseries_reg}
\phi_\ast=\pi(\phi)\exp\left(\sum_{n\geq2} \frac{(-1)^{n-1}}{n} (\phi\mid \letter 1\letter 0^{n-1}) \letter 1^n\right)  \in R\langle\langle \mathbb{N}\rangle\rangle
\end{align}
and $\pi:\R\langle \langle \letter 0,\letter1 \rangle\rangle \to R\langle\langle \mathbb{N}\rangle\rangle$ sends any word starting with $\letter 0$ to $0$ any $\letter 1 \letter 0^{k_1-1}\cdots \letter 1 \letter 0^{k_r-1}$ to $\letter{k}_1\cdots \letter{k}_r$. The set $\operatorname{DM}_0(R)$ is the subset of all $\phi\in\operatorname{DM}(R)$ additionally satisfying
\begin{itemize}
    \item[4.] $(\phi\mid \letter 1\letter 0)=0 $.
\end{itemize}
\end{defi}

The definition of $\phi_\ast$ in \eqref{eq:powseries_reg} is exactly the reformulation of the map $\rho$ in Definition \ref{def:rho} on power series.
\begin{proposition}
Let $\phi\in R\langle\langle \letter 0, \letter 1\rangle\rangle$ be grouplike. Then, we have for any $w\in \Q\langle \mathbb{N}\rangle$
\begin{align*}
 (\phi_\ast\mid w)=  (\phi\mid \rho(w)).
\end{align*}
\end{proposition}

\begin{proof}
Fix $\phi\in R\langle\langle \letter 0, \letter 1\rangle\rangle$ grouplike and denote by $e_k\in R$ the coefficient of $T^k$ in 
\[
\exp\left(\sum_{n\geq2} \frac{(-1)^{n-1}}{n}(\phi\mid \letter 1 \letter 0^{n-1})T^n\right).
\]
As $\phi$ is grouplike, the map $\Q\langle \letter 0, \letter 1\rangle\to R,\ u\mapsto (\phi\mid u)$ is an algebra morphism and we obtain from the definition of $\lambda_k$ from \eqref{eq:def_lambda}
\begin{align*}
\sum_{k\geq 0} (\phi \mid \lambda_k) T^k=\left(\phi\ \middle|\ \exp_{\shuffle}\left(\sum_{k\geq2} \frac{(-1)^{k-1}}{k}\letter 1 \letter 0^{k-1}T^k\right)\right)=\exp\left(\sum_{k\geq2} \frac{(-1)^{k-1}}{k} (\phi\mid \letter 1\letter 0^{k-1}) T^k\right).
\end{align*}
This implies $(\phi\mid \lambda_k)=e_k$ for all $k\geq0$. Hence, we calculate with Definition \ref{def:rho}
\begin{align*}
(\phi\mid \rho(w))=\sum_{w=v\letter 1^k} (\phi\mid \lambda_k\shuffle \gamma(v))=\sum_{w=v\letter 1^k} (\phi\mid \gamma(v))(\phi\mid \lambda_k)=\sum_{w=v\letter 1^k} (\pi(\phi)\mid v)e_k=(\phi_\ast\mid w).
\end{align*}
In the third step we used that $\pi$ is dual to $\gamma$ under the pairing $(-\mid-)$, and in the last step we applied the definition of $\phi_\ast$.
\end{proof}

The assignment $\operatorname{DM}:R\mapsto \operatorname{DM}(R)$ is a functor from the category of commutative, unital $\Q$-algebras into the category of sets. By construction, it is represented by the algebra $\formZ$, i.e., we have bijections
\[
\operatorname{DM}(R)\simeq \operatorname{Hom}_{\Q\operatorname{
-alg}}(\formZ,R).
\]
In this bijection, a power series $\phi=\sum_{w\in \{\letter 0, \letter 1\}^*} (\phi\mid w) w$ maps to the algebra morphism
\begin{align} \label{eq:DM_Homalg}
\formZ\to R,\quad \zeta^{\mathfrak f}(w) \mapsto (\phi\mid w).
\end{align}
In the same way, we get a functor $\operatorname{DM}_0:R\mapsto \operatorname{DM}_0(R)$ represented by $\formZ/(\zeta^{\mathfrak f}(2))$.

\begin{theorem}[{\cite[Theorem~III]{Racinetthesis}}] $\operatorname{DM}_0$ is an affine group scheme.
\end{theorem}

We explain how the shuffle and stuffle varieties from the previous subsection relate to $\operatorname{DM}$.

First, we consider the regularized shuffle algebra $\algshreg$ and the corresponding graded $\mathbb Q$-algebra $R_j(\algshreg)$ from Definition \ref{def:rel_algebra_scheme}. To relate $R_j(\algshreg)$ to the affine group scheme $\operatorname{DM}_0$, we want to consider instead of the projective scheme the affine scheme defined by $R_j(\algshreg)$. We denote this by $\mathcal{R}_j^{\mathbb Q, \operatorname{aff}}(\algshreg) := \mathrm{Spec} \ R_j(\algshreg)$. Recall that this variety naturally comes with an embedding into the affine space $(\algshreg_j)^\vee$, where it is cut out by $I_{\algshreg,j}$.

We equip $\mathbb Q\langle\langle \letter 0, \letter 1\rangle\rangle$ with a scheme structure by writing
\[\mathbb Q\langle\langle \letter 0, \letter 1\rangle\rangle := \varprojlim_i \mathbb Q\langle \letter 0, \letter 1\rangle^\vee_{\leq i}.\]
In other words, $\mathbb Q\langle\langle \letter 0, \letter 1\rangle\rangle := \mathrm{Spec}(\Sym \mathbb Q\langle \letter 0, \letter 1\rangle)$. Now consider the composition

\[\begin{tikzcd}
    \Sym \mathbb Q\langle \letter 0,\letter 1\rangle \rar & \Sym \algshreg \rar & \algshreg
\end{tikzcd}\]
corresponding to the morphisms of schemes

\[\begin{tikzcd}
    \mathrm{Spec} \ \algshreg \rar & (\algshreg)^\vee \rar & \mathbb Q\langle\langle \letter 0,\letter 1\rangle\rangle
\end{tikzcd}\]

where now $(\algshreg)^\vee := \varprojlim_i (\bigoplus_{j\leq i} \algshreg_j)^\vee$. One can check that this is a closed immersion which is given on $R$-valued points as the inclusion
\[(\mathrm{Spec} \ \algshreg)(R) = \{\phi \in R\langle\langle \letter 0, \letter 1 \rangle\rangle \mid (\phi \mid \letter 0) = (\phi \mid \letter 1) = 0, \ \Delta_{\shuffle}(\phi)=\phi\otimes\phi \} \subseteq R\langle\langle \letter 0,\letter 1\rangle\rangle.\]
These are the power series satisfying 1 and 2 in Definition \ref{def:DM0}. On the other hand, consider the composition
\[\begin{tikzcd}
    \Sym S^\shuffle \rar{\Sym \rho^{-1}} & \Sym S^\stuffle \rar & \Sym \algstreg \rar & \algstreg
\end{tikzcd}\]
and the associated morphisms of schemes
\[\begin{tikzcd}
     \mathrm{Spec} \ \algstreg \rar & (\algstreg)^\vee \rar &(S^\stuffle)^\vee \rar & (S^\shuffle)^\vee \subseteq \mathbb Q\langle\langle \letter 0,\letter 1\rangle\rangle
\end{tikzcd}\]
This is again a closed immersion.
By construction, an element $\phi \in (\mathrm{Spec} \ \algshreg)(R)$ lies in $(\mathrm{Spec} \ \algstreg)(R)$ if and only if it satisfies 3 in Definition \ref{def:DM0}. In other words, we have
\[\mathrm{DM} =\mathrm{Spec} \ \algshreg \cap  \mathrm{Spec} \ \algstreg\]
with the intersection taken in $\mathbb Q\langle\langle\letter0, \letter 1\rangle\rangle$.
Note that the intersection scheme is precisely the vanishing locus in $(S^\shuffle)^\vee$ of the ideal $(\operatorname{Sym} \rho_j)(\ker \psi_\ast)+\ker \psi_\shuffle$ introduced below \eqref{eq:sh-st-diagram}.
Write $p_j$ for the projection $\mathbb Q\langle\langle\letter0, \letter 1\rangle\rangle \to \mathbb Q\langle\langle\letter0, \letter 1\rangle\rangle_j$. Then we obtain
\begin{equation}\label{eq:rel_sh/st_variety_DM0}
    p_j(\mathrm{DM}) \subseteq p_j(\mathrm{Spec} \ \algshreg) \cap  p_j(\mathrm{Spec} \ \algstreg)
\end{equation}
where we mean the scheme-theoretic image.
Now observe that by definition 
\[p_j(\mathrm{Spec} \ \algshreg) = \relvar^{\mathbb Q,\mathrm{aff}}_j(\algshreg) \quad \text{and} \quad p_j(\mathrm{Spec} \ \algstreg) = \relvar^{\mathbb Q,\mathrm{aff}}_j(\algstreg)\]
under the appropriate linear embeddings $(\algshreg_j)^\vee \to \mathbb Q\langle\langle \letter 0,\letter 1\rangle\rangle_j$ and, via $\rho_j$, $(\algstreg_j)^\vee \to \mathbb Q\langle\langle \letter 0,\letter 1\rangle\rangle_j$. If the operations given by intersection and projecting to the $j$-th component commute, we get equality in \eqref{eq:rel_sh/st_variety_DM0}. This is equivalent to Conjecture \ref{conj:inner-form-approx-eq}.

Note that all in all we have the inclusions
\[\relvar_j(\motZ) \subseteq \relvar_j(\formZ) \subseteq  \relvar_j(\algshreg) \cap \relvar_j(\algstreg).\]
Conjecturally, all three varieties are isomorphic.

\begin{remark}
(i) From an explicit calculation of the right hand side of \eqref{eq:rel_sh/st_variety_DM0}, one obtains restrictions on what explicit elements in $\operatorname{DM}(R)$ for specific $R$ look like in low weights. For example, there is no explicit element known in the case $R=\Q$.
\vspace{0.1cm}\\
(ii) By \cite[Corollary~0.3]{FuMZVAss}, the pro-unipotent Grothendieck--Teichmüller group $\operatorname{GRT}_1$ embedds into $\operatorname{DM}_0$. Actually, it is conjectured that the double shuffle group $\operatorname{DM}_0$ and $\operatorname{GRT}_1$ coincide. Thus, the right hand side of \eqref{eq:rel_sh/st_variety_DM0} should also describe $\operatorname{GRT}_1$ and might help to understand elements in there.\footnote{A preprint posted by OpenAI on October 6, 2026 announced a proof that the Lie algebra of $\operatorname{GRT}_1$ is freely generated by one element in each odd weight $\geq3$, i.e., that the Lie algebra of $\operatorname{GRT}_1$ coincides with the one of the motivic MZVs from Section \ref{sec:mot-var}. The preprint \cite{willwacher} provides a simplified version of the OpenAI proof. In \cite{hirose} an independent reproof of this statement is given, using the Hirose--Sato confluence relations.} 
\end{remark}

\subsection{Connection to path signatures and time series signatures}\label{rem:iss-and-iis}

    Our work in this article is heavily inspired by how varieties of iterated-integral signature tensors \cite{bib:AFS2019, colmenarejo2020toric,  roughveronese} and iterated-sums \cite{bellingeri2025discrete} were studied in previous work.
    We want to illustrate for the interested reader how the present setting of MZVs fits into this context.

    Given a bounded variation (i.e. finite length) path $X:[a,b]\to \mathbb  R^2$,
    the \emph{iterated-integrals signature} $\IIS(X)$ is the infinite collection of values
    \begin{equation}\label{eq:IIS}
    \langle\IIS(X),\letter{i}_1\cdots\letter{i}_s\rangle:=\int_{0\leq t_1\leq \cdots\leq t_s\leq 1} \intd X_{i_1}(t_1)\,\cdots \,\intd X_{i_s}(t_s)
    \end{equation}
    for arbitrary $k\in\mathbb N_0$ and letters $\letter{i}_j$ from the alphabet $\{\letter 0, \letter{1},\dots,\letter{d-1}\}$.
    We refer the reader to \cite{chevyrevkormilitzin2016}, \cite{bib:AFS2019} and \cite{friz2010multidimensional} for more details.
    Note that when $X$ is piecewise continously differentiable, we may rewrite the Riemann-Stieltjes integrals as Riemann integrals over the derivatives,
    \begin{equation*}
        \int_{0\leq t_1\leq \cdots\leq t_s\leq 1} \intd X_{i_1}(t_1)\,\cdots \,\intd X_{i_s}(t_s)=\int_{0\leq t_1\leq \cdots\leq t_s\leq 1}X_{i_1}'(t_1)\cdots X_{i_s}'(t_s)\intd t_1\cdots \intd t_s.
    \end{equation*}

    We generalize this setup to accommodate the integrals representation of MZVs. 
    Instead of a finite length path defined on a compact time domain, we take the unbounded, infinite length (but still locally bounded variation) path
    \begin{equation}\label{eq:defn_Y_log}
        Y:(0,1)\to \mathbb R^2,\quad t\mapsto (Y
    _0(t),Y_1(t)):=(\log(t),-\log(1-t))
    \end{equation}
    with domain an open time interval.
    The main conceptual difference thus lies in the fact that, while in the usual setting of the IIS we deal with proper Riemann--Stieltjes that always exist when we restrict to sufficiently locally regular paths, we now have improper integrals whose convergence may fail.

    \begin{figure}[h]
\centering
\begin{tikzpicture}
\pgfmathsetmacro{\R}{5}
\pgfmathsetmacro{\S}{ln(exp(\R)-1)}
\begin{axis}[
  width=7cm, height=7cm, axis equal image,
  axis lines=middle, axis line style={-},
  xmin=-5, xmax=1.5, ymin=-1.5, ymax=5,
  xlabel={}, ylabel={},
  xtick={-4,-3,-2,-1,0,1}, ytick={-1,0,1,2,3,4}, xticklabel=\empty, yticklabel=\empty,
  tick label style={font=\small}, clip=false
]

\addplot[teal!75!black, ultra thick, directed=0.2] coordinates {(0,0) (-\R,0)};
  \addplot[blue!70!black, ultra thick,
    domain=-\S:\S, samples=300, variable=\s,directed=0.53]
    ({-ln(1+exp(-\s))},{ln(1+exp(\s))});
\addplot[orange!85!black, ultra thick, directed=0.8] coordinates {(0,\R) (0,0)};

\end{axis}
\end{tikzpicture}
\caption{In blue: the path $Y:(0,1)\to\mathbb R^2$, $t\mapsto(\log(t),-\log(1-t))$. Using segments of infinite length (green, yellow), it can be closed 'at infinity' to a loop.}\label{figure:loop_at_infinity}
\end{figure}
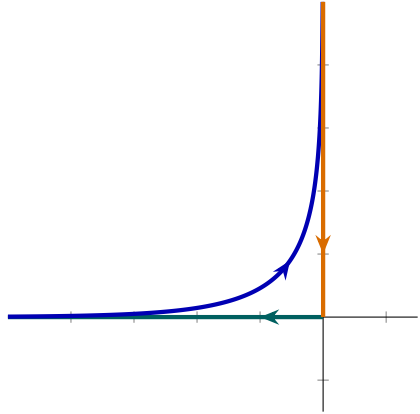
    So let us now extend the previous definition \eqref{eq:IIS}, towards a notion of iterated-improper-integral-signature.    
        Let $X:\,I\to \mathbb R^d$ be a continuous map of locally bounded variation\footnote{We could go to even more generality here. The Riemann--Stieltjes integral locally exists for functions of locally bounded $p$-variation for any $1\leq p<2$. This is known as Young integration due to \cite{youngintegration}, see \cite{friz2010multidimensional} for details. Another possibility is to generalize to discontinuous paths.} from a non-empty, non-singleton interval $I\subseteq \mathbb R$ to $d$-dimensional real space. 
        Let $s\in \mathbb N_0$, $\letter{i}_j\in\{\letter{0},\letter{1},\dots,\letter{d-1}\}$. 
        We write
        \begin{equation*}
            \langle\IIS(X),\letter{i}_1\cdots\letter{i}_s\rangle:=\int_{t_s\in I} 
            \int_{\substack{t_{s-1}\in I,\\ t_{s-1}\leq t_s}} 
            \int_{\substack{t_{s-2}\in I,\\ t_{s-2}\leq t_{s-1}}}
            \cdots
            \int_{\substack{t_1\in I,\\t_1\leq t_2}}\intd X_{i_1}(t_1)\cdots \intd X_{i_s}(t_s)
        \end{equation*}
        whenever all of the improper Riemann--Stieltjes integrals involved converge absolutely. 
    That said, we have the following.
    \begin{proposition}
        Let $Y$ be defined as in \eqref{eq:defn_Y_log}. 
        For any $k_1,\dots,k_r\geq 1$, $k_r\geq 2$, we have
        \begin{equation*}
            \zeta(k_1,\dots,k_r)=\langle\IIS(Y),\letter{1}\letter{0}^{k_1-1}\letter{1}\letter{0}^{k_2-1}\cdots\letter{1}\letter{0}^{k_r-1}\rangle
        \end{equation*}
    \end{proposition}
    \begin{proof}
        Clearly, $Y_0'(t)=\tfrac{1}{t}=\omega_0(t)$ and $Y_1'(t)=\tfrac{1}{1-t}=\omega_1(t)$.
        Note that a Riemann--Stieltjes integral $\int_a^b f\intd g$ equals $\int_a^b f(t)g'(t)\intd t$ for continuously differentiable $g$, and our $Y$ is indeed continuously differentiable.
        By a comparison with \eqref{eqn:zetaasi},
        the statement follows.
    \end{proof}

A similar approach works for time series, i.e.\ we switch from a continuous to a discrete time domain. 
The classical setting there (see \cite{diehleftapia2020acta}) is a time series $Z$, that is any map $Z:\, J\to \mathbb R^d$ where $J$ is some finite integer range, usually by convention starting from $0$ or $1$. 
Automatically, due to finiteness of $J$, $Z$ is bounded and more specifically of bounded variation.

For simplicity, let us restrict to the case of $d=1$, sufficient for our purposes here.
Then given such a $Z$, the iterated-sums signature $\ISS(X)$ is the infinite collection of values
\begin{equation*}
    \langle\ISS(Z),\letter k_{1}\cdots\letter k_{r}\rangle:=\sum_{a<n_1<\cdots<n_r<b-1}(\Delta_{n_1} Z)^{k_1}\,\cdots\,(\Delta_{n_r} Z)^{k_r}
\end{equation*}
where $a=\min J$, $b=\max J$ and $\Delta_n Z:=Z(n+1)-Z(n)$.

However, the time series we are interested in here is unbounded and defined over an infinite time interval, namely
\begin{equation}\label{eq:W_series}
     W:\mathbb N\to \mathbb R,\quad n\mapsto\sum_{k=1}^{n}\tfrac{1}{k}.
\end{equation}
Therefore, let us now consider any map $Z:\, J\to \mathbb R$ from a non-empty integer interval $J\subseteq \mathbb Z$ to $d$-dimensional space. 
    Let $r\in\mathbb N_0$, $k_1,\ldots,k_r\in\mathbb N$. Then we put
    \begin{equation*}
    \langle\ISS(Z),\letter k_1\cdots\letter k_r\rangle=\sum_{\{n_1<\cdots <n_r\}\subseteq J}(\Delta_{n_1} X)^{k_1}\cdots(\Delta_{n_r} X)^{k_r}
    \end{equation*}
    for all words $\letter k_1\cdots\letter k_r$ for which the series converges absolutely.

Once again, we note that the definition can be readily generalized to all $d\in\mathbb N$.

\begin{proposition}
        Let $W$ be defined as in \eqref{eq:W_series}. 
        For any $k_1,\dots,k_r\geq 1$, $k_r\geq 2$, we have
        \begin{equation*}
            \zeta(k_1,\dots,k_r)=\langle\ISS(W),\letter k_1\cdots\letter k_r\rangle
        \end{equation*}
\end{proposition}
\begin{proof}
    Immediate from $\Delta_n 
    W=\tfrac{1}{n}$ by comparison with \eqref{eq:sum-def-mzv}.
\end{proof}

\begin{example}
    For the purpose of this example, let us also define  $\langle\IIS(X),\letter{i}_1\cdots\letter{i}_k\rangle$ and $\langle\IIS(Z),\letter{k}_1\cdots\letter{k}_r\rangle$ in the case that we have unambiguous divergence towards either $-\infty$ or $+\infty$.
   Both 
    \begin{equation*}
      \langle \IIS(Y),\letter 0\rangle=\int_0^1\intd Y_{0}(t_1)=\lim_{t\nearrow 1}Y_0(t)-\lim_{t\searrow  0}Y_0(t)=+\infty
    \end{equation*}
    and 
    \begin{equation*}
    \langle\IIS(Y),\letter 1\rangle=\int_0^1\intd Y_{1}(t_1)=\lim_{t\nearrow 1}Y_1(t)-\lim_{t\searrow 0} Y_1(t)=+\infty
    \end{equation*}
    diverge, which is the fact that the path $Y$ possesses infinite increment.
$\langle \IIS(Y),\letter 0\rangle$ does not have an MZV interpretation, while $\langle\IIS(Y),\letter 1\rangle$ is $\zeta(1)$ which diverges logarithmically.
Similarly, 
\begin{equation*}
\langle \ISS(W),\letter 1\rangle=\sum_{n_1=1}^\infty\Delta_{n_1}W=\sum_{n_1=1}^\infty\tfrac{1}{n_1}=+\infty=\zeta(1),
\end{equation*}
meaning the time series $W$ also possesses infinite increment. 
\end{example}

Finally, we conjecture that the limiting iterated-integrals signature of $Y$ pre- and post-concatenated with infinite length linear paths towards a loop starting and ending in the origin (see Figure \ref{figure:loop_at_infinity}) exists for all indices $\letter j_1\cdots\letter j_r$, and yields the shuffle regularization of MZVs with the convention $\zetashreg(\letter 1)=0$.
To formulate this, let for any $l\in\mathbb N$ the continuous finite length path $Y_l:\,[-1,2]\to\mathbb R^2$ be given as
\begin{equation*}
    Y_l(t):=
    \begin{cases}
        \big(\log(-t),0\big), & \text{if }\, t\in[-1,-\tfrac{1}{l}),\\
        \big(\log(\tfrac{1}{l}),-\tfrac{1}{2}(lt   +1)\log(1-\tfrac{1}{l})\big),& \text{if }\, t\in [-\tfrac{1}{l},\tfrac{1}{l}],\\
        \big(\log(t),-\log(1-t)\big),& \text{if }\, t\in(\tfrac{1}{l},1-\tfrac{1}{l}),\\
        \big(\tfrac{1}{2}(-l(t-1)+1)\log(1-\tfrac{1}{l}),-\log(\tfrac{1}{l})\big),& \text{if }\, t\in [1-\tfrac{1}{l},1+\tfrac{1}{l}],\\
        \big(0,-\log(t-1)\big),&\text{if }\,t\in (1+\tfrac{1}{l},2].
        
    \end{cases}
\end{equation*}

Clearly, $\lim_{l\to \infty}Y_l(t)=Y(t)$
for $t\in (0,1)$.
On $[-1,0)$, $Y_l$ converges pointwise to an infinite line segment from the origin $(0,0)$ to the point at infinity $(-\infty,0)$. 
On $(1,2]$, $Y_l$ converges pointwise to an infinite line segment from the point at infinity $(0,+\infty)$ to the origin $(0,0)$. 
\begin{conjecture}
    We have
    \begin{equation*}
        \lim_{l\to\infty}\,\langle\IIS(Y_l),w\rangle=\zetashreg(w)
    \end{equation*}
    for any word $w$ in letters $\{\letter 0,\letter 1\}$.
\end{conjecture}
The idea is that obviously $\langle\IIS(Y_l),\letter 0\rangle=\langle\IIS(Y_l),\letter 1\rangle=0$ (since $Y_l$ is a loop) and that we \emph{expect} that $\langle \IIS(Y_l),\letter 1 w\letter 0\rangle$ converges to $\langle \IIS(Y),\letter 1w\letter 0\rangle=\zetashreg(\letter 1w\letter 0)$ for any admissible word $\letter 1w\letter 0$ as $l\to \infty$.
The latter needs a proof of course\footnote{We strongly believe this is true and think we have an argument but will leave it open for now.}, but morally speaking, the fact that we integrate against the first component first (which is constant zero on the infinite linear segment in the beginning) and the zeroth component last (which is constant zero on the infinite linear segment in the end) should lead to the two infinite linear segments not contributing to the iterated integral in the limit.
Provided that this holds, we obtain the statement for all words via shuffle relations.
Indeed, as used before in Section~\ref{sec:conj}, any word in letters $\{\letter 0,\letter 1\}$ is a shuffle polynomial in $\letter 0$, $\letter 1$ and words of the form $\letter 1w\letter 0$.

It is not clear to us what could be a similar conjecture for the $\ISS$ and $\zetastreg$.
If we were to naively close $W$ 'at infinity' to a loop based in the origin in the most direct way, this would correspond to the limiting sequence
\begin{equation*}
    W_l(n):=\begin{cases}
        \sum_{k=1}^{n}\tfrac{1}{k},& n\leq l,\\
        0,& n>1.
    \end{cases}
\end{equation*}
But then $\lim_{l\to\infty} \langle\ISS(W_l),\letter{2}\rangle=\lim_{l\to\infty}\sum_{k=1}^l\tfrac{1}{k^2}+\lim_{l\to\infty}(\sum_{k=1}^l\tfrac{1}{k})^2=\zeta(2)+(+\infty)=+\infty$.

\section{Computational approach}
\label{sec:computation}
In this section, we address our main question computationally. We compute and tabulate a list of polynomial equations satisfied by MZVs in weight $j$ and also the transformations that allow us to express them in different bases. Before we do this, we give a short summary of the \emph{MZV data mine} \cite{Blumlein:2009cf}. The reason for this is that we compute an exhaustive list of polynomial equations in a given weight (similar in spirit to the MZV data mine), and we use the datamine \cite{Blumlein:2009cf} to transform the relations we find to other bases. 
\subsection{Summary of the multiple zeta value data mine}
The MZV data mine is a formal computational approach to relations between MZVs where the one of the main computational problems studied is to determine linear relations at a fixed weight and depth and to reduce the corresponding MZVs to explicit and ``\emph{nice}" bases. Using extended double shuffle relations, the datamine formally determines the set of all relations between MZVs at a fixed weight and depth. The choice of bases in the data mine are the \emph{Fibonacci basis} and the \emph{Lyndon basis}. The Fibonacci basis is a basis at fixed weight consisting of nested sums of that weight. The {Lyndon basis} is a basis for the graded ring of MZVs in which as many elements as possible are represented as products of lower weight basis elements, and are generated by Lyndon words \cite[Section~2]{Blumlein:2009cf}. 

The computational side of the datamine is implemented in \texttt{FORM}.
The main result is an explicitly proven reduction of all MZVs up to and including weight \(22\) to bases whose dimensions agree with the Zagier conjecture (\ref{conj:dimQZ}) (as well as the \emph{Broadhurst--Kreimer} conjecture \cite[Equation 7]{BroadhurstKreimer}, which refines Conjecture \ref{conj:dimQZ} to include the depth filtration).
The data mine further uses modular and depth limited computations to provide further evidence at larger weights. The data mine also tests Hoffman's \(\{2,3\}\) basis in the aforementioned range of weight.

In the computational approach presented here, we  do not attempt to reproduce or utilize the complete \texttt{FORM} reduction system of \cite{Blumlein:2009cf}.  Instead, we focus on a fixed weight and generate a monomial basis of Lyndon words in a $\{2,3\}$ alphabet, compute generating toric relations (which can be distinguished from additional numerical MZV relations, if any), and then transform these relations to different bases by an exact change of basis using the datamine and \sage.\footnote{The motivic implementation in \href{https://doc.sagemath.org/html/en/reference/modfrm/sage/modular/multiple_zeta.html}{\sage{}} uses a semi-numerical evaluation routine in conjunction with a motivic change of basis.} We further use the datamine to re-write these relations in a depth minimal basis, but we refrain from doing so for all examples here, and instead refer the reader to the \github{} repository linked in Section \ref{sec:code}.

\subsection{Numerical approaches for toric geometry of MZVs}

\subsubsection{Main computational question and result}
Recall our main question: Determine the set of all $\mathbb Q$-relations of degree $d$ between weight $j$ MZVs. In the previous section we have identified relations that arise from toric geometry. We now want to numerically verify the absence of additional relations to a certain extent, which is expected from the period conjecture. Moreover, we want to transform the toric relations to relations in other (putative) linear bases.

\subsubsection{Computational strategy}
The computation has three key components:
\begin{enumerate}
    \item Construct the free Lie algebra from Lyndon words in $\{2,3\}$. Use this as a basis to construct a monomial basis to find the toric MZV relations at weight \(j\), degree $d$.
    \item Numerically search for relations as well as check that there are no further relations beyond toric relations. 
    \item Transform the Lyndon basis to other bases.
\end{enumerate}

\begin{figure}[htbp]
\centering
\resizebox{\textwidth}{!}{\begin{tikzpicture}[
  >=Latex,
  node distance=9mm and 13mm,
  box/.style={
    draw,
    rounded corners=2pt,
    align=center,
    minimum height=18mm,
    text width=38mm,
    inner sep=4pt
  },
  gp/.style={box,fill=cyan!7},
  pari/.style={box,fill=purple!8},
  sage/.style={box,fill=cyan!8},
  exact/.style={box,fill=gray!8},
  io/.style={box,fill=orange!8},
  output/.style={box,fill=blue!7},
  arrow/.style={->,thick},
  dashedarrow/.style={->,thick,dashed}
]

\node[io] (input) {Input: \\Weight \(j\), Degree \(d\)};

\node[gp,right=of input] (model) {\gp: \\Find toric relations in degree $d$};
\node[io,below=of input] (m2input) {Input: Weight $j$};
\node[gp,below=of model] (m2) {\texttt{Macaulay2}: Find generating toric relations};
\node[sage,right=of model] (hoff) {\gp: \\Numerical checks \&  certification };
\node[gp,right=of hoff] (rewrite) {\sage:\\ Basis transformation};
\node[io,below=of rewrite] (normalize) {Output: \\Integer relations};
\draw[arrow] (input) -- (model);
\draw[arrow] (m2input) -- (m2);
\draw [dashedarrow, thick] (m2) -| (hoff);

\draw[dashedarrow] (model) -- (hoff);
\draw[arrow] (hoff) -- node[above,sloped] {} (rewrite);
\draw[arrow] (rewrite) -- (normalize);
\end{tikzpicture}
}
\caption{Computational flowchart: Solid arrows denote exact combinatorial and
motivic computations. The dashed arrows denote numerical computations.}
\label{fig:globalworkflow}
\end{figure}
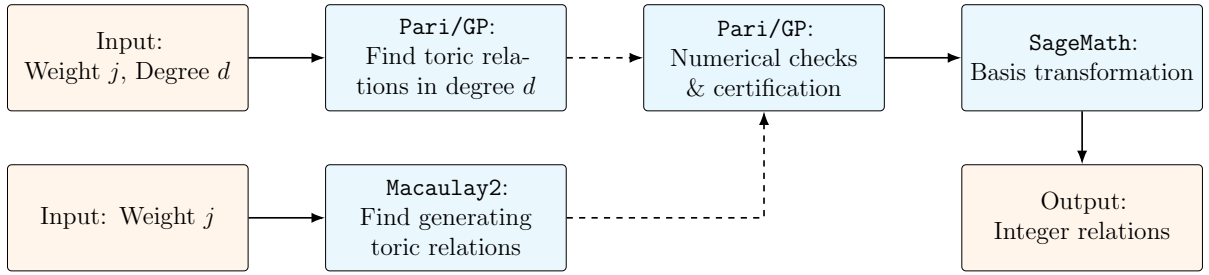

Figure~\ref{fig:globalworkflow} shows the computational flowchart. We use \gp{} to construct the coordinate the monomial basis and its exact formal relations, as well as for fast numerics which we need to rule out \emph{exceptional relations}, which are relations not predicted in the toric geometry setting (See Section \ref{sec:exceptional}). This certification is done using the \texttt{LLL} algorithm \cite{lll}.  We use \sage{} to compute exact rational change-of-basis matrices in the motivic MZV algebra. These matrices can then be exported back to \gp, where all degree-\(d\) relations are rewritten by
ordinary exact polynomial substitution.

The computational steps from (1) to (2) above is combinatorial in nature and can be programmed in any computer algebra system, while going from (2) to (3) is an exact map which we compute in \sage. 

\subsubsection{Details of computational approach}
\label{sec:computationalapproach}
\begin{enumerate}[label=(\textbf{Part \arabic*}), ]
\item This part is a short re-cap of the toric formalism explained in Section \ref{sec:mot-var}. We however explain the key computational aspects again since (a) we implement this in afresh \gp{} (which is more frequently used in the MZV community), and (b) it helps readability of the coding workflow as well as help understand the code. 

We start with a {formal Lyndon word generator}, in which symbols attached to strict Lyndon words are treated as independent indeterminates. Lyndon words are primitive words that are lexicographically smaller than all their proper cyclic shifts. They are used in this case to construct a basis of free Lie algebras. Let $L_1, \ldots, L_r$ represent the fundamental Lyndon generators constructed from the alphabet $\{2,3\}$ (which eventually makes a connection with the Hoffman basis). We first compute strict Lyndon words in the
alphabet \(\{2,3\}\). 
\\[0.2cm] We then construct polynomials  comprised of monomials \(X_i\) of total weighted degree \(j\) in those formal Lyndon generators. For each strict Lyndon word \(\lambda\) over the ordered alphabet \(2<3\), we introduce a formal variable \(L_\lambda\).  At a fixed coordinate weight \(j\), we take all strict Lyndon words of weights at most \(j\) and enumerate each exponent vector
\( \displaystyle
E_i=(e_{i1},\ldots,e_{ir})\in\mathbb N^r,
\ 
\sum\nolimits_{a=1}^r \operatorname{wt}(\lambda_a)e_{ia}=j,
\)
exactly once.  The corresponding monomial is
\(\displaystyle
X_i=\prod\nolimits_{a=1}^r L_{\lambda_a}^{e_{ia}},
\)
where  the symbols \(L_\lambda\) are treated as formal variables. These monomials are ordered by length and then lexicographic order.
To compute the exact relations, we consider a map \[
\Psi_j:\Q[X_1,\ldots,X_N]\longrightarrow
\Q[L_{\lambda_1},\ldots,L_{\lambda_r}],
\qquad
X_i\longmapsto L^{E_i}.
\]
For an exponent vector \(\alpha=(\alpha_1,\ldots,\alpha_N)\) with \(|\alpha|=d\),
\(\displaystyle
\Psi_j(X^\alpha)=L^{\sum_i\alpha_iE_i}.
\)
Two degree-\(d\) monomials in $X_1,\ldots, X_N$ have the same image when their image of the exponent vectors agree.
We can then enumerate the degree-\(d\) monomial basis, and construct the set of monomials grouped by the image above. By choosing the first monomial in each set as a representative,\footnote{Note that the monomials are ordered, as mentioned before.} for every other monomial in the same group we store the difference with that representative. For a set of size \(m\), the \(m-1\) differences form a basis of the kernel on that set, and  therefore we get a basis of $(\ker\Psi_j)_d$, as well as formal toric MZV relations at weight $j$ and degree $d$.  For example, for a few cases of degree $d$ and weight $j$, we compute that 
\[
\begin{array}{c|c|c|c|c}
(j,d)&N&\#\{X^\alpha:|\alpha|=d\}&\#\{\Psi_j(X^\alpha)\}&\dim(\ker\Psi_j)_d\\ \hline
(12,2)&12&78&74&4\\
(14,2)&21&231&213&18\\
(14,3)&21&1771&1437&334
\end{array}
\]

A table with only $\dim(\ker\Psi_j)_d$ (i.e., the number of relations) for various values of $j,d$ is shown in Table \ref{tab:numericsnumrel}.

If \(\lambda_a\) is the Lyndon word associated with \(L_{\lambda_a}\), the formal Lyndon basis can be mapped onto a MZV i.e, 
\(L_{\lambda_a}\longmapsto \zeta(\lambda_a).
\)
Consequently
\(\displaystyle
X_i\longmapsto \prod\nolimits_a\zeta(\lambda_a)^{e_{ia}}.
\)
Every formal relation that we get here is a tautological relation. These formal relations can however be mapped to a different basis. In this paper, we will eventually choose the Hoffman basis (which is also why we chose the Lyndon basis in $\{2,3\}$).  

\item More importantly, we need to use numerics to detect relations that are not toric in nature. When we constructed the toric relations in the \textbf{Part 1} of the computational program, we looked at those relations that come from $\ker(\Psi_j)$ at a degree $d$. These were tautological in the sense that they arise from trivial cancellations of monomials. However, there could be relations that do not arise in this way. These are what we call \emph{exceptional relations}. To exclude these (to a certain extent) we use the method described in Section \ref{sec:lll}.

\item We can now express the relations in other bases of weight $j$ MZVs. 
For a monomial
\(
X_i=\prod_a\zeta(\lambda_a)^{e_{ia}},
\)
we can construct
\( \displaystyle
\varphi(X_i)=\prod\nolimits_a\varphi(\zeta(\lambda_a))^{e_{ia}}
\)
directly in the \(F\)-algebra \cite{brown2012decomposition}. The auxiliary \(F\)-algebra is the graded shuffle algebra
\[
F
=
\mathbb Q[f_2]
\otimes
\mathbb Q\langle f_3,f_5,f_7,\ldots\rangle_{\shuffle},
\]
where the symbol $f_2$ is the ordinary polynomial generator and the symbols $f_{1+2\mathbb{Z}_{>0}}$ are the odd generators of the formal shuffle algebra. The auxiliary $F$ algebra provides an explicit basis of alphabets for the algebra \(\mathcal H\) of motivic MZVs.  The computational advantage of using the $F$-algebra is that each fixed weight component is a finite-dimensional \(\mathbb Q\)-vector space, so motivic identities and changes of basis reduce to exact rational linear transformations. $\varphi$ as mentioned above is a non-canonical graded algebra isomorphism from the algebra of motivic MZVs to the $F$-algebra
\(
\varphi:\mathcal H \xrightarrow{\sim}  F.
\) In a nutshell, \(\varphi\) sends a motivic MZV to its exact \(f\)-alphabet expansion, and this is implemented in \sage. We refer the reader to \cite{brown2012decomposition} for more details concerning $\varphi$.
For an ordered single MZV basis\footnote{Not to be confused with single zeta.} \(B=(B_1,\ldots,B_m)\), we construct two matrices $P_B$ (square, non-singular) and $V_X$. \(P_B\) is constructed in such a way that the vectors \(\varphi(B_j)\) are columns of \(P_B\), and \(V_X\) similarly has vectors \(\varphi(X_i)\) as columns.  We then compute the matrix \(A_B\) as
\(
P_BA_B=V_X.
\) Consequentially, \( X = A_B B\). $X$ is the basis formed from the monomials $X_i$.  Once we have computed the transformation of bases, we can perform a substitution to compute all formal relations in the target basis. 
The target basis can be a depth minimal basis (See Appendix \ref{sec:depthminimal} for an example).\\[1em]
As an example, we will consider the case when \(B\) is the Hoffman basis. In this basis, every weight-\(j\) coordinate \(X_i\) is expressed as a rational linear combination of single weight-\(j\) (motivic) MZVs whose multiset entries lie in \(\{2,3\}\). These are the Hoffman basis/words. The Hoffman basis is a basis for (motivic) MZVs \cite[Corollary~7.5]{brown2012mixed}.
We start with constructing the Hoffman basis (with some well defined choice of lexicographic ordering - 
often the default in a computer algebra system. It is however important to pay attention to this ordering.)
\[
 H_j=\big\{ (s_1,\ldots,s_r):s_i\in\{2,3\},\ \sum\nolimits_i s_i=j\big\}.
\]
In our computations, we use an explicit parser between \gp{} and \sage{} at this stage.\footnote{Some of this functionality is also available in \sage. We refer the reader to the \sage{} \href{https://doc.sagemath.org/html/en/reference/modfrm/sage/modular/multiple_zeta.html}{documentation on MZVs}. For example, in \sage, we can use the functions \texttt{phi()} and \texttt{phi\textunderscore as\textunderscore vector()}, which implements \(\varphi\).}

\end{enumerate}

\subsubsection{Certifying absence of bounded integer relations with \texttt{LLL}}
\label{sec:lll}

Let $x_1,\ldots,x_n$ be exact real numbers, such as monomials in multiple
zeta values. Given an integer $C\ge1$, we wish to prove that
\[
 a\in\mathbb Z^n,\quad |a_i|\le C,\quad \sum_{i=1}^n a_i x_i=0
 \quad\Longrightarrow\quad a=0,
\]
    i.e., the absence of non-trivial integer relations between $x_1,\ldots,x_n$ whose coefficient vector $a$  satisfies the given coefficient bound. For coefficients $a_i, \ i = 1,\cdots,n$ with 23 or fewer decimal digits, we have $C=10^{23}-1$. A general computational strategy to detect integer relations between real numbers is to approximate them to high precision and use the \texttt{LLL} algorithm \cite{lll}. However, \texttt{LLL} does not necessarily find the shortest coefficient vector. Thus, failing to find a relation within the given bound using \texttt{LLL} does not prove the absence of such relations. However, \texttt{LLL} can be used to obtain a lower bound on the coefficients of \emph{every} integer relation, as we explain in the following. The method used here seems to be folklore, however we were not able to find an exact reference. We make no claim of originality.

Start by choosing $\ell_i,u_i, D$ such that
\(\displaystyle
 {\ell_i} \big/ {D} \le x_i\le {u_i} \big/ {D},
 \ \ell_i,u_i\in\Z, \ D\in\mathbb Z_{>0}.
\)
Choose two positive integers $S,K$ and $n$ integers $r_i$ such that $r_i$ is an approximation for $Sx_i$ for all $i = 1,\cdots, n$.
Let
\(\displaystyle
 \eta_i=\max\{|Dr_i-S\ell_i|,\ |Dr_i-Su_i|\},
\)
and $ \displaystyle H:=\sum\nolimits_{i=1}^n\eta_i$.
Then $|r_i-Sx_i|\le\eta_i/D$. Consider the $(n+1)\times n$ integer matrix
\(\displaystyle
 B=\begin{pmatrix}K I_n\\ r_1\ \cdots\ r_n\end{pmatrix}\) and the lattice \( L\subset \R^{n+1}\) generated by its column vectors. If $a$ is a nonzero relation, then $B\cdot a\ne0$. Let us suppose that $|a_i|\le C$. We then have the inequality
\[
 \left|\sum_i a_i r_i\right|
 =\left|\sum_i a_i(r_i-Sx_i)\right|\le\frac{CH}{D},
\]
from which it follows that
\begin{equation}\label{eq:upper}
 \|B \cdot a\|_2^2
 =K^2\sum_i a_i^2+\left(\sum_i a_i r_i\right)^2
 \le C^2\left(K^2n+\frac{H^2}{D^2}\right)=:R^2.
\end{equation}
Thus any relation within the coefficient bound would force the existence a nonzero lattice
vector of length at most $R$.

\begin{lemma}[Gram--Schmidt orthogonalization]\label{lem:lower}
Let $b_1,\ldots,b_n$ be a lattice basis and let $b_1^*,\ldots,b_n^*$ be its
Gram--Schmidt orthogonalization. Every nonzero vector $v$ in the lattice satisfies
\[
 \|v\|_2^2\ge\min_{1\le i\le n}\|b_i^*\|_2^2.
\]
\end{lemma}

\begin{proof}
Write $v=\sum_i z_i b_i$ with $z_i\in\Z$, and let $j$ be the largest index
with $z_j\ne0$. The vector $b_j^*$ is orthogonal to all original basis vectors.
The orthogonal projection of $v$ onto its span is therefore $z_j b_j^*$.
Since $|z_j|\ge1$,
\[
 \|v\|_2^2\ge z_j^2\|b_j^*\|_2^2
 \ge\min_i\|b_i^*\|_2^2.\qedhere
\]
\end{proof}

\noindent Note however that this lemma applies to \emph{any} lattice basis (not just \texttt{LLL} reduced bases).

The \texttt{LLL} algorithm constructs a reduced lattice basis which is related to the original basis by some unimodular transformation $U$, i.e., $\det U=\pm1$. This means that the columns of $BU$ form a basis of the same lattice $L$. We can construct the integer Gram
matrix $G=(BU)^T(BU)$, and if $\Delta_i$ is the determinant of its leading
$i\times i$ principal submatrix, with $\Delta_0=1$, then these determinants are
    positive, and are related to the squared Gram--Schmidt lengths by
\[
 \|b_i^*\|_2^2=\frac{\Delta_i}{\Delta_{i-1}}.
\]

\begin{theorem}
If the following inequality
\begin{equation}\label{eq:certificate}
 D^2\Delta_i>
 \Delta_{i-1}C^2\bigl(K^2nD^2+H^2\bigr)
\end{equation}
holds for every $i=1,\ldots,n$,  then there is no nonzero integer relation among the $x_i$ with
$|a_i|\le C$.
\end{theorem}

\begin{proof}
Inequality~\eqref{eq:certificate} asserts that every squared Gram--Schmidt
length exceeds $R^2$. Lemma~\ref{lem:lower} then makes every nonzero vector
in $L$ longer than $R$. This contradicts~\eqref{eq:upper} for any 
relation within the coefficient bound.
\end{proof}

In our case, the numbers $x_1,\ldots,x_n$ are chosen as a putative basis of $(\Z_{\langle j\rangle})_{jd}$ for a given $j$ and $d$, where $\Z_{\langle j\rangle}$ denotes the $\mathbb Q$-subalgebra of $\Z$ generated in degree $j$. The value of $n$ for given $j$ and $d$ is the value $\mathcal H(R_j(\motZ))(d)$ of the Hilbert function of the algebra $R_j(\motZ)$ defined in Definition \ref{def:rel_algebra_scheme}. Also see Table \ref{tab:numericsnumrel} below. We use the algorithm described in \cite{akhilesh} (as well as Theorems 3 and 5) to obtain approximations $x_i \sim \ell_i/D, u_i/D $ with certified containments $x_i \in [\ell_i/D,u_i/D]$ for every $i$.
\subsection{Results}
\subsubsection{Number of relations}
In Table \ref{tab:numericsnumrel}, we provide data for the number of relations computed from \textbf{Part 1} of the computational approach described in Section \ref{sec:computationalapproach}.  Note that the number of linearly independent relations for given $j$ and degree $d$ is by definition the value of the Hilbert function \[\mathcal H(I_{\motZ,j})(d) := \dim \ (I_{\motZ,j})_d.\]
Note that we have
\begin{equation}
\label{eq:relhilbfcb}
    \mathcal H(I_{\motZ,j})(d) = \binom{\dim(\motZ_j) + d-1}{d} - \mathcal H(R_j(\motZ))(d),
\end{equation} with $R_j(\motZ)$ as in Definition \ref{def:rel_algebra_scheme}. In the cases where $\mathcal R_j(\motZ)$ is normal (see Theorem \ref{thm:mot-non-normal}), a classical theorem in toric geometry says that $\mathcal H(R_j(\motZ))$ agrees with the Ehrhart polynomial of the polytope $P_{\motZ,j}$ for sufficiently large $d$, see \cite[Proposition~9.4.3]{coxlittlescheck}.
\begin{table}[h!]
\centering
\begin{tabular}{||l|l|l|l|l||}
\hline
$j$\textbackslash{}$\deg$ & $\deg = 2$ & $\deg = 3$                 & $\deg = 4$ 
& $\deg = 5$
\\ \hline
$10$ & $1$ & $7$ &$28$  &  $84$  \\
$11$ & $1$ & $9$ &$45$  & $165$\\
$12$ & $4$ & $46$ &  $285$&   $1266$     \\
$13$ & $8$ & $116$ & $899$ & $4948$      \\
$14$ &$18$ & $334$ &   $3248$   &  $22092$     \\
$15$ &$36$ & $862$ & $10723$  &    $92506$         \\
$16$ &$73$ &   $2254$ & $35782$ &  $390578$ \\
$17$ &$140$&    $5594$  & $114317$      &             \\
$18$ & $278$   &   $14466$ &  $382003$    &             \\
$19$ &  $514$   &   $34790$                &      &             \\
$20$ &  $990$   &   $87431$                &      &             \\
$21$ & $1830$  &                   &      &             \\
$22$ & $3413$ &                   &      &             \\
$23$ & $6294$ &                   &      &             \\
\hline
\end{tabular}
    \caption{Table of number of MZV relations at weight $j$ and degrees $2,3,4,5$ that were obtained computationally. Due to the exponential complexity of the algorithms, the computing power to determine the number of relations as degree increases increases substantially. These numbers can be checked against \eqref{eq:relhilbfcb}.}
    \label{tab:numericsnumrel}
\end{table}

\subsubsection{Exceptional relations}
\label{sec:exceptional}
 An \textit{exceptional relation} is a relation that is not predicted by the toric formalism explained in Section \ref{sec:mot-var}, based on Brown's motivic result \cite{brown2012decomposition} and the Grothendieck period conjecture. 
We tested for the existence of exceptional relations numerically (as explained in Part 3 of Section \ref{sec:computationalapproach} and Section \ref{sec:lll}) and can certify that there are no such integer relation vectors whose entries have fewer than 23 decimal digits in the cases listed in Table \ref{tab:noexceptional1}, in accordance with the Grothendieck period conjecture.

\begin{table}[!h]
    \centering
    \begin{tabular}{c|c}
        degree & weights checked \\ \hline
        2 & 10, 11, 12, 13, 14, 15\\
        3 & 10, 11, 12\\
        4 & 10
    \end{tabular}
    \caption{Degrees and weights in which we certified the absence of exceptional relations within the coefficient bound specified above.}
    \label{tab:noexceptional1}
\end{table}

\subsubsection{Code and examples}
\label{sec:code}
All code is available at \url{https://github.com/abhirammk/ToricMZVCode}. We also give a thoroughly worked out example for the case of weight $j = 14, \deg = 2$, which can be found in the \texttt{Results} subfolder in the above \github{} repository. The repository also contains a datamine for the generating relations expressed in a depth minimal basis, as well as in the Hoffman basis.

\bibliography{refs}

@article{bib:AFS2019,
  author = {Am{\'e}ndola, Carlos and Friz, Peter and Sturmfels, Bernd},
 title = {Varieties of signature tensors},
 fjournal = {Forum of Mathematics, Sigma},
 journal = {Forum Math. Sigma},
 issn = {2050-5094},
 volume = {7},
 pages = {54},
 note = {Id/No e10},
 year = {2019},
 doi = {10.1017/fms.2019.3},
 zbMATH = {7049009},
 Zbl = {1411.14068},
 eprint = {1804.08325}
}

@phdthesis{Racinetthesis,
  title={Séries génératrices non-commutatives de polyzêtas et associateurs de Drinfeld},
  author={Racinet, Georges},
  school={Laboratoire Amiénois de Mathématique Fondamentale et Appliquée, Université d’Amiens},
  year={2000},
  url = {https://theses.hal.science/tel-00110891/document}
}

@book{ZhaoMZVbook,
   author = {Zhao, Jianqiang},
 title = {Multiple zeta functions, multiple polylogarithms and their special values},
 fseries = {Series on Number Theory and its Applications},
 series = {Ser. Number Theory Appl.},
 issn = {1793-3161},
 volume = {12},
 isbn = {978-981-4689-39-7; 978-981-4689-41-0},
 year = {2016},
 publisher = {Hackensack, NJ: World Scientific},
 doi = {10.1142/9634},
 zbMATH = {6560232},
 Zbl = {1367.11002}
}

@article{FuMZVAss,
  author = {Furusho, Hidekazu},
 title = {Double shuffle relation for associators},
 fjournal = {Annals of Mathematics. Second Series},
 journal = {Ann. Math. (2)},
 issn = {0003-486X},
 volume = {174},
 number = {1},
 pages = {341--360},
 year = {2011},
 doi = {10.4007/annals.2011.174.1.9},
 zbMATH = {5960703},
 Zbl = {1321.11088},
 eprint = {0808.0319}
}

@article {Hoffman,
    author = {Hoffman, Michael E.},
 title = {Quasi-shuffle products},
 fjournal = {Journal of Algebraic Combinatorics},
 journal = {J. Algebr. Comb.},
 issn = {0925-9899},
 volume = {11},
 number = {1},
 pages = {49--68},
 year = {2000},
 doi = {10.1023/A:1008791603281},
 zbMATH = {1428365},
 Zbl = {0959.16021},
 eprint = {math/9907173},
}

@article{Ihara_Kaneko_Zagier_2006, 
author = {Ihara, Kentaro and Kaneko, Masanobu and Zagier, Don},
 title = {Derivation and double shuffle relations for multiple zeta values},
 fjournal = {Compositio Mathematica},
 journal = {Compos. Math.},
 issn = {0010-437X},
 volume = {142},
 number = {2},
 pages = {307--338},
 year = {2006},
 doi = {10.1112/S0010437X0500182X},
 zbMATH = {5033724},
 Zbl = {1186.11053}
}

@article{roughveronese,
author = {Galuppi, Francesco},
 title = {The rough {Veronese} variety},
 fjournal = {Linear Algebra and its Applications},
 journal = {Linear Algebra Appl.},
 issn = {0024-3795},
 volume = {583},
 pages = {282--299},
 year = {2019},
 doi = {10.1016/j.laa.2019.08.029},
 zbMATH = {7125955},
 Zbl = {1432.14041},
 eprint = {1809.02522}
}

@article{bellingeri2025discrete,
 author = {Bellingeri, Carlo and Penaguiao, Raul},
 title = {Discrete signature varieties},
 fjournal = {SIAM Journal on Applied Algebra and Geometry},
 journal = {SIAM J. Appl. Algebra Geom.},
 issn = {2470-6566},
 volume = {9},
 number = {3},
 pages = {640--665},
 year = {2025},
 doi = {10.1137/23M1606629},
 zbMATH = {8087507},
 eprint = {2303.13377}
}

@article{diehleftapia2020acta,
	 author = {Diehl, Joscha and Ebrahimi-Fard, Kurusch and Tapia, Nikolas},
 title = {Time-warping invariants of multidimensional time series},
 fjournal = {Acta Applicandae Mathematicae},
 journal = {Acta Appl. Math.},
 issn = {0167-8019},
 volume = {170},
 pages = {265--290},
 year = {2020},
 doi = {10.1007/s10440-020-00333-x},
 zbMATH = {7344875},
 Zbl = {1460.62145},
 eprint = {1906.05823}
}

@misc{chevyrevkormilitzin2016,
	author = "Chevyrev, Ilya and Kormilitzin, Andrey",
	title = "{A Primer on the Signature Method in Machine Learning}",
	year = 2016,
    eprint = {1603.03788v1}
}

@article{brown2012mixed,
author = {Brown, Francis},
 title = {Mixed {Tate} motives over {{\(\mathbb{Z}\)}}},
 fjournal = {Annals of Mathematics. Second Series},
 journal = {Ann. Math. (2)},
 issn = {0003-486X},
 volume = {175},
 number = {2},
 pages = {949--976},
 year = {2012},
 doi = {10.4007/annals.2012.175.2.10},
 zbMATH = {6025005},
 Zbl = {1278.19008},
 eprint = {1102.1312}
}

@incollection {brown2012decomposition,
    author = {Brown, Francis C. S.},
 title = {On the decomposition of motivic multiple zeta values},
 booktitle = {Galois-Teichm\"uller theory and arithmetic geometry. Selected papers based on the presentations at the workshop and conference, Kyoto, Japan, October 25--30, 2010},
 isbn = {978-4-86497-014-3},
 pages = {31--58},
 year = {2012},
 publisher = {Tokyo: Mathematical Society of Japan},
 zbMATH = {6490730},
 Zbl = {1321.11087},
 eprint = {1102.1310}
}

@article{colmenarejo2020toric,
  author = {Colmenarejo, Laura and Galuppi, Francesco and Micha{\l}ek, Mateusz},
 title = {Toric geometry of path signature varieties},
 fjournal = {Advances in Applied Mathematics},
 journal = {Adv. Appl. Math.},
 issn = {0196-8858},
 volume = {121},
 pages = {35},
 note = {Id/No 102102},
 year = {2020},
 doi = {10.1016/j.aam.2020.102102},
 zbMATH = {7304610},
 Zbl = {1467.14149},
 eprint = {1903.03779}
}

@article{Blumlein:2009cf,
   author = {Bl{\"u}mlein, J. and Broadhurst, D. J. and Vermaseren, J. A. M.},
 title = {The multiple zeta value data mine},
 fjournal = {Computer Physics Communications},
 journal = {Comput. Phys. Commun.},
 issn = {0010-4655},
 volume = {181},
 number = {3},
 pages = {582--625},
 year = {2010},
 doi = {10.1016/j.cpc.2009.11.007},
 zbMATH = {5832387},
 Zbl = {1221.11183},
 eprint = {0907.2557}
}

@book{coxlittlescheck,
  author = {Cox, David A. and Little, John B. and Schenck, Henry K.},
 title = {Toric varieties},
 fseries = {Graduate Studies in Mathematics},
 series = {Grad. Stud. Math.},
 issn = {1065-7339},
 volume = {124},
 isbn = {978-0-8218-4819-7},
 year = {2011},
 publisher = {Providence, RI: American Mathematical Society (AMS)},
 zbMATH = {5934474},
 Zbl = {1223.14001}
}

@book{Reutenauer,
 author = {Reutenauer, Christophe},
 title = {Free {Lie} algebras},
 fseries = {London Mathematical Society Monographs. New Series},
 series = {Lond. Math. Soc. Monogr., New Ser.},
 volume = {7},
 isbn = {0-19-853679-8},
 year = {1993},
 publisher = {Oxford: Clarendon Press},
 zbMATH = {417855},
 Zbl = {0798.17001}
}

@misc{telen2022introductiontoricgeometry,
      title={Introduction to Toric Geometry}, 
      author={Simon Telen},
      year={2022},
      archivePrefix={arXiv},
      primaryClass={math.AG},
      eprint = {2203.01690}, 
}

@unpublished{zagierMZVtalk,
    author       = {Don Zagier},
    title        = {Relations among multiple zeta values {II}},
    note = {Grothendieck-Teichmüller groups, deformation and operads programme, Isaac Newton Institute, Cambridge.
    \url{https://www.newton.ac.uk/seminar/3542/}}
}

@article {lehmer,
    AUTHOR = {Lehmer, D. H.},
     TITLE = {On the maxima and minima of {B}ernoulli polynomials},
   JOURNAL = {Amer. Math. Monthly},
  FJOURNAL = {American Mathematical Monthly},
    VOLUME = {47},
      YEAR = {1940},
     PAGES = {533--538},
      ISSN = {0002-9890,1930-0972},
   MRCLASS = {33.0X},
  MRNUMBER = {2378},
MRREVIEWER = {W.\ E.\ Milne},
       DOI = {10.2307/2303833},
       URL = {https://doi.org/10.2307/2303833},
}

@book{friz2010multidimensional,
 author = {Friz, Peter K. and Victoir, Nicolas B.},
 title = {Multidimensional stochastic processes as rough paths. {Theory} and applications.},
 fseries = {Cambridge Studies in Advanced Mathematics},
 series = {Camb. Stud. Adv. Math.},
 volume = {120},
 isbn = {978-0-521-87607-0},
 year = {2010},
 publisher = {Cambridge: Cambridge University Press},
 doi = {10.1017/CBO9780511845079},
 zbMATH = {5677244},
 Zbl = {1193.60053}
}

@misc{mzvfromnumberstomotives,
title={Multiple zeta values:
from numbers to motives},
author={{Burgos Gil}, José Ignacio and Fresán, Javier},
year=2026,
url = {http://javier.fresan.perso.math.cnrs.fr/mzv.pdf}
}

@article{BoutetDeMonvel2004DivergentMultizeta,
  author    = {Boutet de Monvel, L.},
  title     = {Remark on Divergent Multizeta Series},
  fjournal   = {S\={u}ri Kaiseki Kenky\={u}sho K\={o}ky\={u}roku},
  journal = {RIMS K\={o}ky\={u}roku (Lecture Notes of the RIMS)},
  volume    = {1397},
  year      = {2004},
  pages     = {1--9},
  publisher = {Research Institute for Mathematical Sciences, Kyoto University},
  url       = {http://hdl.handle.net/2433/25985},
  note = {Workshop: Microlocal Analysis and Asymptotic Analysis}
  }

@article{lll,
  title = {Factoring Polynomials with Rational Coefficients},
  author = {Lenstra, A. K. and Lenstra, H. W. and Lovász, L.},
  date = {1982-12-01},
  journal = {Mathematische Annalen},
  shortjournal = {Math. Ann.},
  volume = {261},
  number = {4},
  year = {1982},
  pages = {515--534},
  issn = {1432-1807},
  doi = {10.1007/BF01457454},
}

@incollection{apery,
     author = {Ap\'ery, Roger},
     title = {Irrationalit\'e de $\zeta{2}$ et $\zeta{3}$},
     booktitle = {Journ\'ees Arithm\'etiques de Luminy},
     series = {Ast\'erisque},
     pages = {11--13},
     year = {1979},
     publisher = {Soci\'et\'e math\'ematique de France},
     number = {61},
     mrnumber = {3363457},
     zbl = {0401.10049},
     url = {https://www.numdam.org/item/AST_1979__61__11_0/}
}

@article{beukers1,
 author = {Beukers, F.},
 title = {A note on the irrationality of {{\(\zeta (2)\)}} and {{\(\zeta (3)\)}}},
 fjournal = {Bulletin of the London Mathematical Society},
 journal = {Bull. Lond. Math. Soc.},
 issn = {0024-6093},
 volume = {11},
 pages = {268--272},
 year = {1979},
 doi = {10.1112/blms/11.3.268},
 zbMATH = {3654230},
 Zbl = {0421.10023}
}

@article{poorten,
author = {Apostol, Tom M.},
 title = {A proof that {Euler} missed: evaluating {{\(\zeta(2)\)}} the easy way},
 fjournal = {The Mathematical Intelligencer},
 journal = {Math. Intell.},
 issn = {0343-6993},
 volume = {5},
 number = {3},
 pages = {59--60},
 year = {1983},
 doi = {10.1007/BF03026576},
 zbMATH = {3855198},
 Zbl = {0538.10001}
}

@incollection{Brown2014MotivicPeriods,
 author = {Brown, Francis},
 title = {Motivic periods and {{\(\mathbb{P}^1\setminus\{0,1,\infty\}\)}}},
 booktitle = {Proceedings of the international congress of mathematicians (ICM 2014), Seoul, Korea, August 13--21, 2014. Vol. II: Invited lectures},
 isbn = {978-89-6105-805-6; 978-89-6105-803-2},
 pages = {295--318},
 year = {2014},
 publisher = {Seoul: KM Kyung Moon Sa},
 zbMATH = {6797661},
 Zbl = {1373.11064},
 eprint = {1407.5165}
}

@misc{BCK2024FormalMZV,
  author        = {Burmester, Annika and Confurius, Niclas and K{\"u}hn, Ulf},
  title         = {{AGZT}-Lectures on Formal Multiple Zeta Values},
  year          = {2024},
  eprint        = {2406.13630},
  archivePrefix = {arXiv},
  primaryClass  = {math.NT}
}

@unpublished{WaldschmidtMZV,
  author       = {Waldschmidt, Michel},
  title        = {Lectures on Multiple Zeta Values},
  year         = {2011},
  note         = {Lecture notes, Special Year in Number Theory, The Institute of Mathematical Sciences (IMSc), Chennai, April 2011.},
  url          = {https://webusers.imj-prg.fr/~michel.waldschmidt/articles/pdf/MZV2011IMSc.pdf}
}

@article{Hof97,
 author = {Hoffman, Michael E.},
 title = {The algebra of multiple harmonic series},
 fjournal = {Journal of Algebra},
 journal = {J. Algebra},
 issn = {0021-8693},
 volume = {194},
 number = {2},
 pages = {477--495},
 year = {1997},
 doi = {10.1006/jabr.1997.7127},
 zbMATH = {1066135},
 Zbl = {0881.11067}
}

@article{BroadhurstKreimer,
   title={Association of multiple zeta values with positive knots via Feynman diagrams up to 9 loops},
   volume={393},
   ISSN={0370-2693},
   url={http://dx.doi.org/10.1016/S0370-2693(96)01623-1},
   DOI={10.1016/s0370-2693(96)01623-1},
   number={3-4},
   journal={Physics Letters B},
   publisher={Elsevier BV},
   author={Broadhurst, D.J. and Kreimer, D.},
   year={1997},
   month=Feb, pages={403–412},
   eprint = {hep-th/9609128}
   }

@article{terasomamzv,
   title={Mixed Tate motives and multiple zeta values},
   volume={149},
   ISSN={1432-1297},
   url={http://dx.doi.org/10.1007/s002220200218},
   DOI={10.1007/s002220200218},
   number={2},
   journal={Inventiones Mathematicae},
   publisher={Springer Science and Business Media LLC},
   author={Terasoma, Tomohide},
   year={2002},
   month=Aug, pages={339–369},
   eprint = {math/0104231},
   }

@misc{delignegoncharov,
      title={Groupes fondamentaux motiviques de Tate mixte}, 
      author={P. Deligne and A. B. Goncharov},
      year={2003},
      eprint={math/0302267},
      archivePrefix={arXiv},
      primaryClass={math.NT}
}

@incollection{gkz,
 author = {Gangl, Herbert and Kaneko, Masanobu and Zagier, Don},
 title = {Double zeta values and modular forms},
 booktitle = {Automorphic forms and zeta functions. In memory of Tsuneo Arakawa. Proceedings of the conference, Rikkyo University, Tokyo, Japan, September 4--7, 2004},
 isbn = {981-256-632-5},
 pages = {71--106},
 year = {2006},
 publisher = {Hackensack, NJ: World Scientific},
 zbMATH = {5064852},
 Zbl = {1122.11057},
 nourl = {https://people.mpim-bonn.mpg.de/zagier/files/mpim/05-96/fulltext.pdf}
}

@incollection{zagier1994,
 author = {Zagier, Don},
 title = {Values of zeta functions and their applications},
 booktitle = {First European congress of mathematics (ECM), Paris, France, July 6-10, 1992. Volume II: Invited lectures (Part 2)},
 isbn = {3-7643-2799-5},
 pages = {497--512},
 year = {1994},
 publisher = {Basel: Birkh{\"a}user},
 zbMATH = {727429},
 Zbl = {0822.11001},
 nourl = {https://people.mpim-bonn.mpg.de/zagier/files/scanned/ValuesofZetaFunsAndApplications/fulltext.pdf}
}

@misc{willwacher,
      title={The {D}eligne--{D}rinfeld conjecture after {OpenAI}}, 
      author={Thomas Willwacher},
      year={2026},
      eprint={2610.15592},
      archivePrefix={arXiv},
      primaryClass={math.QA},
}

@misc{hirose,
      title={Confluence relations, {H}offman basis over $\mathbb{Z}_{(2)}$, and the {D}eligne--{D}rinfeld conjecture}, 
      author={Minoru Hirose},
      year={2026},
      eprint={2610.23719},
      archivePrefix={arXiv},
      primaryClass={math.NT},
}

@incollection{youngintegration,
 author = {Young, L. C.},
 title = {An inequality of the {H{\"o}lder} type connected with {Stieltjes} integration.},
 booktitle = {Comptes rendus du congr\`es international des math\'ematiciens, Oslo, 1936. Tome 2: Conf\'erences de sections},
 year = {1937},
 publisher = {Oslo: A.W. Br{{\o}}ggers Boktrykkeri},
 language = {English},
 zbMATH = {2520391},
 JFM = {63.0191.03}
}

@article{Zudilin2001,
  author  = {Zudilin, Wadim},
  title   = {One of the numbers $\zeta(5)$, $\zeta(7)$, $\zeta(9)$, $\zeta(11)$ is irrational},
  journal = {Russian Mathematical Surveys},
  volume  = {56},
  number  = {4},
  pages   = {774--776},
  year    = {2001},
  doi     = {10.1070/RM2001v056n04ABEH000427},
  note    = {Translated from Uspekhi Mat. Nauk 56 (2001), no. 4(340), 149--150}
}

@article{ballrioval,
  title={Irrationalit{\'e} d’une infinit{\'e} de valeurs de la fonction z{\^e}ta aux entiers impairs},
  author={Ball, Keith and Rivoal, Tanguy},
  journal={Inventiones mathematicae},
  volume={146},
  number={1},
  pages={193--207},
  year={2001},
 doi = {10.1007/s002220100168},
  publisher={Springer Berlin Heidelberg, Berlin/Heidelberg}
}

@article{akhilesh,
author = {Akhilesh, P.},
 title = {Double tails of multiple zeta values},
 fjournal = {Journal of Number Theory},
 journal = {J. Number Theory},
 issn = {0022-314X},
 volume = {170},
 pages = {228--249},
 year = {2017},
 doi = {10.1016/j.jnt.2016.06.020},
 zbMATH = {6626762},
 Zbl = {1398.11111},
  publisher={Elsevier},
  eprint = {2105.12156}
}

@book{Sturmfels1996,
  author    = {Sturmfels, Bernd},
  title     = {{Gr{\"o}bner Bases and Convex Polytopes}},
  series    = {University Lecture Series},
  volume    = {8},
  publisher = {American Mathematical Society},
  address   = {Providence, RI},
  year      = {1996}
}

@misc{IbadanlecturesSotille,
    author = {Sotille, Frank},
    title = {Ibadan Lectures on Toric Varieties},
    year = {2024},
    eprint = {1708.01842v2},
    noteb = {\url{https://arxiv.org/abs/1708.01842v2}}
}
\bibliographystyle{habbrv2}

\appendix

\stoptoc

\section{An example on transforming to different bases}\label{sec:depthminimal}

{
\makeatletter
\origaddtocontents{toc}{\protect\contentsline {section}{Appendices}{\thepage }{\@currentHref }{\protected@file@percent }}
\makeatother
}

To transform from the Hoffman basis to a depth minimal basis, one can use the data mine \cite{Blumlein:2009cf}. The numerical and arithmetic advantage is that of such depth minimal bases is that the relations are often much simpler. We do not list this here but only provide an example of how the relation simplifies.  At weight $10$ and degree $2$, from the unique exact toric relation as in \eqref{eq:exmot10ideal}, we convert to the Hoffman basis to find the relation: 

 \begin{footnotesize}

\begin{alignat*}{2}
 \smallnumprint{428071864754152140800} \, \mzvb{2, 2, 2, 2, 2}^2 
&& {} -  \smallnumprint{108594930391617472000} \, \mzvb{2, 2, 2, 2, 2} \mzvb{3, 3, 2, 2} \\
 {} -  \smallnumprint{67275182473279359360} \, \mzvb{2, 2, 2, 2, 2} \mzvb{3, 2, 3, 2}
&& {} -  \smallnumprint{86109773838803347200} \, \mzvb{2, 2, 2, 2, 2} \mzvb{2, 3, 3, 2} \\
 {} -  \smallnumprint{101607970288175904000} \, \mzvb{2, 2, 2, 2, 2} \mzvb{3, 2, 2, 3}
&& {} -  \smallnumprint{53086451554806691840} \, \mzvb{2, 2, 2, 2, 2} \mzvb{2, 3, 2, 3} \\
 {} -  \smallnumprint{114012872158098144000} \, \mzvb{2, 2, 2, 2, 2} \mzvb{2, 2, 3, 3}
&& {} +  \smallnumprint{4113414477444870000} \, \mzvb{3, 3, 2, 2}^2 \\
 {} +  \smallnumprint{11168840049596326200} \, \mzvb{3, 3, 2, 2} \mzvb{3, 2, 3, 2}
&& {} +  \smallnumprint{9351168337459974000} \, \mzvb{3, 3, 2, 2} \mzvb{2, 3, 3, 2} \\
 {} +  \smallnumprint{7780916664213930000} \, \mzvb{3, 3, 2, 2} \mzvb{3, 2, 2, 3}
&& {} +  \smallnumprint{10361147299189752800} \, \mzvb{3, 3, 2, 2} \mzvb{2, 3, 2, 3} \\
 {} +  \smallnumprint{6658321356517230000} \, \mzvb{3, 3, 2, 2} \mzvb{2, 2, 3, 3}
&& {} +  \smallnumprint{7679117027914369203} \, \mzvb{3, 2, 3, 2}^2 \\
 {} +  \smallnumprint{12773030553349878120} \, \mzvb{3, 2, 3, 2} \mzvb{2, 3, 3, 2}
&& {} +  \smallnumprint{10515846056424462900} \, \mzvb{3, 2, 3, 2} \mzvb{3, 2, 2, 3} \\
 {} +  \smallnumprint{14171810536210292214} \, \mzvb{3, 2, 3, 2} \mzvb{2, 3, 2, 3}
&& {} +  \smallnumprint{8943099438566381400} \, \mzvb{3, 2, 3, 2} \mzvb{2, 2, 3, 3} \\
 {} +  \smallnumprint{5329635054385411200} \, \mzvb{2, 3, 3, 2}^2
&& {} +  \smallnumprint{8822017078080588000} \, \mzvb{2, 3, 3, 2} \mzvb{3, 2, 2, 3} \\
 {} +  \smallnumprint{11812627355170667280} \, \mzvb{2, 3, 3, 2} \mzvb{2, 3, 2, 3} 
&& {} +  \smallnumprint{7527773153838888000} \, \mzvb{2, 3, 3, 2} \mzvb{2, 2, 3, 3} \\
 {} +  \smallnumprint{3679026192079087500} \, \mzvb{3, 2, 2, 3}^2 
&& {} +  \smallnumprint{9748921500042269100} \, \mzvb{3, 2, 2, 3} \mzvb{2, 3, 2, 3} \\
 {} +  \smallnumprint{6312640813867440000} \, \mzvb{3, 2, 2, 3} \mzvb{2, 2, 3, 3} 
&& {} +  \smallnumprint{6528625962706643283} \, \mzvb{2, 3, 2, 3}^2 \\
 {} +  \smallnumprint{8317541161037889600} \, \mzvb{2, 3, 2, 3} \mzvb{2, 2, 3, 3} 
&& {} +  \smallnumprint{2715496736380560000} \, \mzvb{2, 2, 3, 3}^2 \quad = \quad 0
\end{alignat*}
\end{footnotesize}

If we choose a depth minimal basis with the following choice of depth 2 and depth 3 basis elements
\(    \{ \mzvb{10},\mzvb{1,9},\mzvb{3,7},\mzvb{5,5},\mzvb{1,1,8},\mzvb{1,2,7},\mzvb{1,4,5}\}
\)
in weight 10, we find that the relation simplifies to the relation in Example \ref{ex:mot-weight-10-expansion}.

\section{Generating weight \texorpdfstring{$12$}{12} relations in a depth-minimal basis}
\label{sec:wt12deg2ex}
We express four generating quadric relations for weight $12$ motivic MZVs (see Table \ref{tab:mot} and Example \ref{ex:mot-w12}) in the linear basis 
\[
\begin{gathered}
\{\ \ \zeta^{\mathfrak m}(12),\quad \zeta^{\mathfrak m}(1,11),\quad \zeta^{\mathfrak m}(2,10),\quad \zeta^{\mathfrak m}(3,9),\quad
\zeta^{\mathfrak m}(1,1,10),\quad \zeta^{\mathfrak m}(1,2,9), \\[2pt] \zeta^{\mathfrak m}(2,1,9),\quad \zeta^{\mathfrak m}(1,3,8),\quad
\zeta^{\mathfrak m}(2,2,8),\quad \zeta^{\mathfrak m}(3,1,8),\quad \zeta^{\mathfrak m}(1,1,1,9),\quad \zeta^{\mathfrak m}(1,1,2,8)\ \ \}
\end{gathered}
\]
of $\motZ$. This basis was chosen such that the number of depth $\leq r$ MZVs appearing is maximal for all $r$. We used the \texttt{LLL} algorithm to obtain a linear basis of the integer relations arising from the toric formalism that has `small' coefficients. Explicit generating relations for other weights are provided in the accompanying github repository.

\subsection*{Relation 1}
\let\smallnumprint\numprint

\begin{scriptsize}
\begin{equation*}
\begin{alignedat}{3}
\smallnumprint{34129881289} \, \zeta(12)^2
&& {}  - \smallnumprint{31297424801520} \, \zeta(12)\zeta(1,11)
&& {}  + \smallnumprint{373276310529600} \, \zeta(1,11)^2 
\\ {}  - \smallnumprint{16759810980744} \, \zeta(12)\zeta(2,10)
&& {}  + \smallnumprint{129124766646720} \, \zeta(1,11)\zeta(2,10)
&& {}  + \smallnumprint{11166798488976} \, \zeta(2,10)^2
\\ {}  - \smallnumprint{1300542288672} \, \zeta(12)\zeta(3,9)
&& {}  + \smallnumprint{14097989410560} \, \zeta(1,11)\zeta(3,9)
&& {}  + \smallnumprint{2438407610496} \, \zeta(2,10)\zeta(3,9)
\\ {}  + \smallnumprint{133114063104} \, \zeta(3,9)^2
&& {}  - \smallnumprint{1164661682688} \, \zeta(12)\zeta(1,1,10)
&& {}  - \smallnumprint{430629494722560} \, \zeta(1,11)\zeta(1,1,10)
\\ {}  - \smallnumprint{74482268829696} \, \zeta(2,10)\zeta(1,1,10)
&& {}  - \smallnumprint{8132059127808} \, \zeta(3,9)\zeta(1,1,10)
&& {}  + \smallnumprint{124198721224704} \, \zeta(1,1,10)^2
\\ {}  + \smallnumprint{10605975802080} \, \zeta(12)\zeta(1,2,9)
&& {}  - \smallnumprint{108939009081600} \, \zeta(1,11)\zeta(1,2,9)
&& {}  - \smallnumprint{18842240626560} \, \zeta(2,10)\zeta(1,2,9)
\\ {}  - \smallnumprint{2057217338880} \, \zeta(3,9)\zeta(1,2,9)
&& {}  + \smallnumprint{62838638714880} \, \zeta(1,1,10)\zeta(1,2,9)
&& {}  + \smallnumprint{7948339718400} \, \zeta(1,2,9)^2
\\ {}  + \smallnumprint{10213509803520} \, \zeta(12)\zeta(2,1,9)
&& {}  - \smallnumprint{102530832076800} \, \zeta(1,11)\zeta(2,1,9)
&& {}  - \smallnumprint{17733873530880} \, \zeta(2,10)\zeta(2,1,9)
\\ {}  - \smallnumprint{1936204554240} \, \zeta(3,9)\zeta(2,1,9)
&& {}  + \smallnumprint{59142248202240} \, \zeta(1,1,10)\zeta(2,1,9)
&& {}  + \smallnumprint{14961580646400} \, \zeta(1,2,9)\zeta(2,1,9)
\\ {}  + \smallnumprint{7040743833600} \, \zeta(2,1,9)^2
&& {}  + \smallnumprint{1098904795968} \, \zeta(12)\zeta(1,3,8)
&& {}  - \smallnumprint{17942895613440} \, \zeta(1,11)\zeta(1,3,8)
\\ {}  - \smallnumprint{3103427867904} \, \zeta(2,10)\zeta(1,3,8)
&& {}  - \smallnumprint{338835796992} \, \zeta(3,9)\zeta(1,3,8)
&& {}  + \smallnumprint{10349893435392} \, \zeta(1,1,10)\zeta(1,3,8)
\\ {}  + \smallnumprint{2618276613120} \, \zeta(1,2,9)\zeta(1,3,8)
&& {}  + \smallnumprint{2464260341760} \, \zeta(2,1,9)\zeta(1,3,8)
&& {}  + \smallnumprint{215622779904} \, \zeta(1,3,8)^2
\\ {}  + \smallnumprint{784931997120} \, \zeta(12)\zeta(2,2,8)
&& {}  - \smallnumprint{12816354009600} \, \zeta(1,11)\zeta(2,2,8)
&& {}  - \smallnumprint{2216734191360} \, \zeta(2,10)\zeta(2,2,8)
\\ {}  - \smallnumprint{242025569280} \, \zeta(3,9)\zeta(2,2,8)
&& {}  + \smallnumprint{7392781025280} \, \zeta(1,1,10)\zeta(2,2,8)
&& {}  + \smallnumprint{1870197580800} \, \zeta(1,2,9)\zeta(2,2,8)
\\ {}  + \smallnumprint{1760185958400} \, \zeta(2,1,9)\zeta(2,2,8)
&& {}  + \smallnumprint{308032542720} \, \zeta(1,3,8)\zeta(2,2,8)
&& {}  + \smallnumprint{110011622400} \, \zeta(2,2,8)^2
\\ {}  + \smallnumprint{941918396544} \, \zeta(12)\zeta(3,1,8)
&& {}  - \smallnumprint{15379624811520} \, \zeta(1,11)\zeta(3,1,8)
&& {}  - \smallnumprint{2660081029632} \, \zeta(2,10)\zeta(3,1,8)
\\ {}  - \smallnumprint{290430683136} \, \zeta(3,9)\zeta(3,1,8)
&& {}  + \smallnumprint{8871337230336} \, \zeta(1,1,10)\zeta(3,1,8)
&& {}  + \smallnumprint{2244237096960} \, \zeta(1,2,9)\zeta(3,1,8)
\\ {}  + \smallnumprint{2112223150080} \, \zeta(2,1,9)\zeta(3,1,8)
&& {}  + \smallnumprint{369639051264} \, \zeta(1,3,8)\zeta(3,1,8)
&& {}  + \smallnumprint{264027893760} \, \zeta(2,2,8)\zeta(3,1,8)
\\ {}  + \smallnumprint{158416736256} \, \zeta(3,1,8)^2
&& {}  + \smallnumprint{15736215306240} \, \zeta(12)\zeta(1,1,1,9)
&&  \quad = \quad 0 .
\end{alignedat}
\end{equation*}
\end{scriptsize}

\subsection*{Relation 2}

\begin{scriptsize}
\begin{equation*}
\begin{alignedat}{3}
\smallnumprint{7354665389989} \, \zeta(12)^2
&& {}  - \smallnumprint{9643451714788656} \, \zeta(12)\zeta(1,11)
&& {}  -\smallnumprint{4813139118801600} \, \zeta(1,11)^2
\\ {}  - \smallnumprint{2614130369617128} \, \zeta(12)\zeta(2,10)
&& {}  -\smallnumprint{15586214812137792} \, \zeta(1,11)\zeta(2,10)
&& {}  - \smallnumprint{2668659392160432} \, \zeta(2,10)^2
\\ {}  -\smallnumprint{66012670707360} \, \zeta(12)\zeta(3,9)
&& {}  - \smallnumprint{5278741583314176} \, \zeta(1,11)\zeta(3,9)
&& {}  -\smallnumprint{1229174708644224} \, \zeta(2,10)\zeta(3,9)
\\ {}  - \smallnumprint{104039091419904} \, \zeta(3,9)^2
&& {}  + \smallnumprint{4466744176066560} \, \zeta(12)\zeta(1,1,10)
&& {}  - \smallnumprint{7001597297885184} \, \zeta(1,11)\zeta(1,1,10)
\\ {}  + \smallnumprint{6969913950633984} \, \zeta(2,10)\zeta(1,1,10)
&& {}  +  \smallnumprint{2842207470747648} \, \zeta(3,9)\zeta(1,1,10)
&& {}  + \smallnumprint{5591477122891776} \, \zeta(1,1,10)^2
\\ {}  +  \smallnumprint{2482739522208096} \, \zeta(12)\zeta(1,2,9)
&& {}  -\smallnumprint{99500011879680} \, \zeta(1,11)\zeta(1,2,9)
&& {}  +  \smallnumprint{2150628757517952} \, \zeta(2,10)\zeta(1,2,9)
\\ {}  + \smallnumprint{772979062934016} \, \zeta(3,9)\zeta(1,2,9)
&& {}  +  \smallnumprint{1801303902388224} \, \zeta(1,1,10)\zeta(1,2,9)
&& {}  + \smallnumprint{77189654856960} \, \zeta(1,2,9)^2
\\ {}  +  \smallnumprint{2454831863855616} \, \zeta(12)\zeta(2,1,9)
&& {}  - \smallnumprint{1751649056501760} \, \zeta(1,11)\zeta(2,1,9)
&& {}  +  \smallnumprint{1732823867676672} \, \zeta(2,10)\zeta(2,1,9)
\\ {}  + \smallnumprint{695167842410496} \, \zeta(3,9)\zeta(2,1,9)
&& {}  +  \smallnumprint{2654642055020544} \, \zeta(1,1,10)\zeta(2,1,9)
&& {}  + \smallnumprint{389881189785600} \, \zeta(1,2,9)\zeta(2,1,9)
\\ {}  +  \smallnumprint{298527538544640} \, \zeta(2,1,9)^2
&& {}  + \smallnumprint{315504879316032} \, \zeta(12)\zeta(1,3,8)
&& {}  +  \smallnumprint{750351872438784} \, \zeta(1,11)\zeta(1,3,8)
\\ {}  + \smallnumprint{499513272088320} \, \zeta(2,10)\zeta(1,3,8)
&& {}  +  \smallnumprint{144682885315584} \, \zeta(3,9)\zeta(1,3,8)
&& {}  - \smallnumprint{153769845325824} \, \zeta(1,1,10)\zeta(1,3,8)
\\ {}  - \smallnumprint{93857515766784} \, \zeta(1,2,9)\zeta(1,3,8)
&& {}  - \smallnumprint{47806650630144} \, \zeta(2,1,9)\zeta(1,3,8)
&& {}  - \smallnumprint{17896690732032} \, \zeta(1,3,8)^2
\\ {}  + \smallnumprint{293178752634048} \, \zeta(12)\zeta(2,2,8)
&& {}  - \smallnumprint{571367363258880} \, \zeta(1,11)\zeta(2,2,8)
&& {}  + \smallnumprint{165269360215296} \, \zeta(2,10)\zeta(2,2,8)
\\ {}  +  \smallnumprint{82433908896768} \, \zeta(3,9)\zeta(2,2,8)
&& {}  + \smallnumprint{528900676780032} \, \zeta(1,1,10)\zeta(2,2,8)
&& {}  +  \smallnumprint{94543988290560} \, \zeta(1,2,9)\zeta(2,2,8)
\\ {}  + \smallnumprint{117932459212800} \, \zeta(2,1,9)\zeta(2,2,8)
&& {}  +  \smallnumprint{1047310645248} \, \zeta(1,3,8)\zeta(2,2,8)
&& {}  + \smallnumprint{9879043691520} \, \zeta(2,2,8)^2
\\ {}  +  \smallnumprint{304341815975040} \, \zeta(12)\zeta(3,1,8)
&& {}  + \smallnumprint{89492254589952} \, \zeta(1,11)\zeta(3,1,8)
&& {}  +  \smallnumprint{332391316151808} \, \zeta(2,10)\zeta(3,1,8)
\\ {}  + \smallnumprint{113558397106176} \, \zeta(3,9)\zeta(3,1,8)
&& {}  +  \smallnumprint{187565415727104} \, \zeta(1,1,10)\zeta(3,1,8)
&& {}  + \smallnumprint{343236261888} \, \zeta(1,2,9)\zeta(3,1,8)
\\ {}  +  \smallnumprint{35062904291328} \, \zeta(2,1,9)\zeta(3,1,8)
&& {}  - \smallnumprint{17373035409408} \, \zeta(1,3,8)\zeta(3,1,8)
&& {}  + \smallnumprint{10402699014144} \, \zeta(2,2,8)\zeta(3,1,8)
\\ {}  - \smallnumprint{1742584098816} \, \zeta(3,1,8)^2
&& {}  + \smallnumprint{1793928544911360} \, \zeta(12)\zeta(1,1,1,9)
&&  \quad = \quad 0 .
\end{alignedat}
\end{equation*}
\end{scriptsize}

\subsection*{Relation 3}

\begin{scriptsize}
\begin{overlong}
\begin{alignedat}{3}
\smallnumprint{4968133332796} \, \zeta(12)^2
&& {}  + \smallnumprint{6813472901257056} \, \zeta(12)\zeta(1,11)
&& {}  + \smallnumprint{35290958384563200} \, \zeta(1,11)^2
\\ {}  - \smallnumprint{6194594944090248} \, \zeta(12)\zeta(2,10)
&& {}  + \smallnumprint{6903017258197824} \, \zeta(1,11)\zeta(2,10)
&& {}  + \smallnumprint{7434234917260128} \, \zeta(2,10)^2
\\ {}  - \smallnumprint{999219911737120} \, \zeta(12)\zeta(3,9)
&& {}  - \smallnumprint{2946402045260544} \, \zeta(1,11)\zeta(3,9)
&& {}  + \smallnumprint{5731192250045184} \, \zeta(2,10)\zeta(3,9)
\\ {}  + \smallnumprint{1323461086385664} \, \zeta(3,9)^2
&& {}  - \smallnumprint{18438625367455872} \, \zeta(12)\zeta(1,1,10)
&& {}  + \smallnumprint{17494589451230208} \, \zeta(1,11)\zeta(1,1,10)
\\ {}  + \smallnumprint{33274469548528128} \, \zeta(2,10)\zeta(1,1,10)
&& {}  + \smallnumprint{14432741576128512} \, \zeta(3,9)\zeta(1,1,10)
&& {}  - \smallnumprint{555937133101056} \, \zeta(1,1,10)^2
\\ {}  + \smallnumprint{706246022073216} \, \zeta(12)\zeta(1,2,9)
&& {}  + \smallnumprint{16122724717810176} \, \zeta(1,11)\zeta(1,2,9)
&& {}  + \smallnumprint{10512653980558464} \, \zeta(2,10)\zeta(1,2,9)
\\ {}  + \smallnumprint{3854117109318144} \, \zeta(3,9)\zeta(1,2,9)
&& {}  - \smallnumprint{8449411855177728} \, \zeta(1,1,10)\zeta(1,2,9)
&& {}  - \smallnumprint{2278653132911616} \, \zeta(1,2,9)^2
\\ {}  + \smallnumprint{449855295142656} \, \zeta(12)\zeta(2,1,9)
&& {}  + \smallnumprint{10135881225639936} \, \zeta(1,11)\zeta(2,1,9)
&& {}  + \smallnumprint{11271742425990144} \, \zeta(2,10)\zeta(2,1,9)
\\ {}  + \smallnumprint{4541796093886464} \, \zeta(3,9)\zeta(2,1,9)
&& {}  - \smallnumprint{809544725987328} \, \zeta(1,1,10)\zeta(2,1,9)
&& {}  - \smallnumprint{2485259360243712} \, \zeta(1,2,9)\zeta(2,1,9)
\\ {}  - \smallnumprint{132788428701696} \, \zeta(2,1,9)^2
&& {}  - \smallnumprint{155993001040992} \, \zeta(12)\zeta(1,3,8)
&& {}  + \smallnumprint{5824465237391616} \, \zeta(1,11)\zeta(1,3,8)
\\ {}  + \smallnumprint{1438948170585216} \, \zeta(2,10)\zeta(1,3,8)
&& {}  + \smallnumprint{310416861282816} \, \zeta(3,9)\zeta(1,3,8)
&& {}  - \smallnumprint{4988736645156864} \, \zeta(1,1,10)\zeta(1,3,8)
\\ {}  - \smallnumprint{1705930426464768} \, \zeta(1,2,9)\zeta(1,3,8)
&& {}  - \smallnumprint{1452188619399168} \, \zeta(2,1,9)\zeta(1,3,8)
&& {}  - \smallnumprint{194865786989568} \, \zeta(1,3,8)^2
\\ {}  - \smallnumprint{361105582585440} \, \zeta(12)\zeta(2,2,8)
&& {}  + \smallnumprint{1034990443655424} \, \zeta(1,11)\zeta(2,2,8)
&& {}  + \smallnumprint{2046218926930560} \, \zeta(2,10)\zeta(2,2,8)
\\ {}  + \smallnumprint{860560048937472} \, \zeta(3,9)\zeta(2,2,8)
&& {}  + \smallnumprint{1123157058195456} \, \zeta(1,1,10)\zeta(2,2,8)
&& {}  - \smallnumprint{48292902001152} \, \zeta(1,2,9)\zeta(2,2,8)
\\ {}  + \smallnumprint{323557382873088} \, \zeta(2,1,9)\zeta(2,2,8)
&& {}  - \smallnumprint{186738128326656} \, \zeta(1,3,8)\zeta(2,2,8)
&& {}  + \smallnumprint{55371049786368} \, \zeta(2,2,8)^2
\\ {}  - \smallnumprint{258549291813216} \, \zeta(12)\zeta(3,1,8)
&& {}  + \smallnumprint{3429727840523520} \, \zeta(1,11)\zeta(3,1,8)
&& {}  + \smallnumprint{1742583548757888} \, \zeta(2,10)\zeta(3,1,8)
\\ {}  + \smallnumprint{585488455110144} \, \zeta(3,9)\zeta(3,1,8)
&& {}  - \smallnumprint{1932789793480704} \, \zeta(1,1,10)\zeta(3,1,8)
&& {}  - \smallnumprint{877111664232960} \, \zeta(1,2,9)\zeta(3,1,8)
\\ {}  - \smallnumprint{564315618263040} \, \zeta(2,1,9)\zeta(3,1,8)
&& {}  - \smallnumprint{288234851152896} \, \zeta(1,3,8)\zeta(3,1,8)
&& {}  - \smallnumprint{37998014376960} \, \zeta(2,2,8)\zeta(3,1,8)
\\ {}  - \smallnumprint{81558216382464} \, \zeta(3,1,8)^2
&& {}  + \smallnumprint{8355930327613440} \, \zeta(12)\zeta(1,1,1,9)
&&  \quad = \quad 0 .
\end{alignedat}
\end{overlong}
\end{scriptsize}

\subsection*{Relation 4}

\begin{scriptsize}
\begin{overlong}
\begin{alignedat}{3}
\smallnumprint{28397571319364} \, \zeta(12)^2
&& {}  - \smallnumprint{56570522483150496} \, \zeta(12)\zeta(1,11)
&& {}  + \smallnumprint{19581243700032000} \, \zeta(1,11)^2
\\ {}  - \smallnumprint{4561542679353384} \, \zeta(12)\zeta(2,10)
&& {}  + \smallnumprint{37303810511390784} \, \zeta(1,11)\zeta(2,10)
&& {}  + \smallnumprint{18556262526326112} \, \zeta(2,10)^2
\\ {}  + \smallnumprint{1405933281392288} \, \zeta(12)\zeta(3,9)
&& {}  + \smallnumprint{17421525965624064} \, \zeta(1,11)\zeta(3,9)
&& {}  + \smallnumprint{13475578121329152} \, \zeta(2,10)\zeta(3,9)
\\ {}  + \smallnumprint{2241464070664704} \, \zeta(3,9)^2
&& {}  + \smallnumprint{47957893340009856} \, \zeta(12)\zeta(1,1,10)
&& {}  + \smallnumprint{25491730325326848} \, \zeta(1,11)\zeta(1,1,10)
\\ {}  + \smallnumprint{19500107928279552} \, \zeta(2,10)\zeta(1,1,10)
&& {}  + \smallnumprint{3048809297614848} \, \zeta(3,9)\zeta(1,1,10)
&& {}  - \smallnumprint{5746936746737664} \, \zeta(1,1,10)^2
\\ {}  + \smallnumprint{13470356212968960} \, \zeta(12)\zeta(1,2,9)
&& {}  - \smallnumprint{28710833213952} \, \zeta(1,11)\zeta(1,2,9)
&& {}  + \smallnumprint{1103364317151360} \, \zeta(2,10)\zeta(1,2,9)
\\ {}  - \smallnumprint{459487010099712} \, \zeta(3,9)\zeta(1,2,9)
&& {}  - \smallnumprint{4043419975268352} \, \zeta(1,1,10)\zeta(1,2,9)
&& {}  - \smallnumprint{616104689624064} \, \zeta(1,2,9)^2
\\ {}  + \smallnumprint{13231446140463360} \, \zeta(12)\zeta(2,1,9)
&& {}  + \smallnumprint{7791014297438208} \, \zeta(1,11)\zeta(2,1,9)
&& {}  + \smallnumprint{5528194037222400} \, \zeta(2,10)\zeta(2,1,9)
\\ {}  + \smallnumprint{948306052374528} \, \zeta(3,9)\zeta(2,1,9)
&& {}  - \smallnumprint{4410462752243712} \, \zeta(1,1,10)\zeta(2,1,9)
&& {}  - \smallnumprint{1655402088296448} \, \zeta(1,2,9)\zeta(2,1,9)
\\ {}  - \smallnumprint{809263096233984} \, \zeta(2,1,9)^2
&& {}  + \smallnumprint{1953129725130528} \, \zeta(12)\zeta(1,3,8)
&& {}  - \smallnumprint{4080616004918016} \, \zeta(1,11)\zeta(1,3,8)
\\ {}  - \smallnumprint{2039524719707520} \, \zeta(2,10)\zeta(1,3,8)
&& {}  - \smallnumprint{749778345180672} \, \zeta(3,9)\zeta(1,3,8)
&& {}  - \smallnumprint{324736707465216} \, \zeta(1,1,10)\zeta(1,3,8)
\\ {}  + \smallnumprint{26928644931072} \, \zeta(1,2,9)\zeta(1,3,8)
&& {}  - \smallnumprint{249981609811968} \, \zeta(2,1,9)\zeta(1,3,8)
&& {}  + \smallnumprint{66222596219904} \, \zeta(1,3,8)^2
\\ {}  + \smallnumprint{1762001667126048} \, \zeta(12)\zeta(2,2,8)
&& {}  + \smallnumprint{2175164099603712} \, \zeta(1,11)\zeta(2,2,8)
&& {}  + \smallnumprint{1500339056349312} \, \zeta(2,10)\zeta(2,2,8)
\\ {}  + \smallnumprint{376456104798720} \, \zeta(3,9)\zeta(2,2,8)
&& {}  - \smallnumprint{618370929045504} \, \zeta(1,1,10)\zeta(2,2,8)
&& {}  - \smallnumprint{311625522307584} \, \zeta(1,2,9)\zeta(2,2,8)
\\ {}  - \smallnumprint{220480893149184} \, \zeta(2,1,9)\zeta(2,2,8)
&& {}  - \smallnumprint{89083011354624} \, \zeta(1,3,8)\zeta(2,2,8)
&& {}  - \smallnumprint{8083654013952} \, \zeta(2,2,8)^2
\\ {}  + \smallnumprint{1857565696128288} \, \zeta(12)\zeta(3,1,8)
&& {}  - \smallnumprint{952725952657152} \, \zeta(1,11)\zeta(3,1,8)
&& {}  - \smallnumprint{269592831679104} \, \zeta(2,10)\zeta(3,1,8)
\\ {}  - \smallnumprint{186661120190976} \, \zeta(3,9)\zeta(3,1,8)
&& {}  - \smallnumprint{471553818255360} \, \zeta(1,1,10)\zeta(3,1,8)
&& {}  - \smallnumprint{142348438688256} \, \zeta(1,2,9)\zeta(3,1,8)
\\ {}  - \smallnumprint{235231251480576} \, \zeta(2,1,9)\zeta(3,1,8)
&& {}  + \smallnumprint{21681090542592} \, \zeta(1,3,8)\zeta(3,1,8)
&& {}  - \smallnumprint{52625159691264} \, \zeta(2,2,8)\zeta(3,1,8)
\\ {}  - \smallnumprint{7736017287168} \, \zeta(3,1,8)^2
&& {}  + \smallnumprint{1369050731642880} \, \zeta(12)\zeta(1,1,1,9)
&&  \quad = \quad 0 .
\end{alignedat}
\end{overlong}
\end{scriptsize}

\section{Relation algebras of the stuffle algebra}
\label{sec:unregstufthm}

Similar to the proof of Theorem \ref{thm:unregstuffalg}, we have the following analogous statement for the whole stuffle algebra.

\begin{theorem}\begin{enumerate}
        \item For $j\geq 2$, $M_j$ generates the sublattice 
        \[ L_j=\Big\{\sum_i a_i w_i\,|\,a_i\in\mathbb Z,\,w_i\in\lyndon (\mathbb N)_{\leq j},\,\sum_i a_i\weight(w_i)\in j\mathbb Z\Big\}\] of index $j$ in $\mathbb Z^{\mathcal L(\mathbb N)_{\leq j}}$.
        \item For $j\geq 4$, we have \begin{equation*}
            \relvarr_j(\StC)=\relvar_j(\StC_{\langle \leq j-2\rangle})
        \end{equation*}
        \item For $j\geq 4$, we have
    \[ \dim \relvar_j(\StC) = \sum_{i=1}^j \lambda_{\mathbb N, i} - 1 \qquad \text{and} \qquad \dim \relvarr_j(\StC) = \sum_{i=1}^{j-2} \lambda_{\mathbb N, i} - 1.\]
        \item $\relvarr_j(\StC)$ is not normal for $j=16, 18$ and $j \geq 20$.
    \item For every $n$ there is a $j$ such that $I_{\StC,j}$ is not generated in degree $\leq n$.
    \end{enumerate}
    
\end{theorem}

\begin{proof}
    \begin{enumerate}
        \item For $j=2$, the statement is obvious.
        For $j\geq 3$, let $\alpha=j\cdot\letter1$ and $\beta=(j-2)\cdot \letter 1+\letter 2$, both in $M_j$. 
        Then $\beta-\alpha =-2\cdot\letter 1+\letter 2$, thus $\alpha,\beta$ generate $L_j\cap \mathbb Z^{\{\letter 1,\letter 2\}}$, as whenever for $n,c_1,c_2,d_1,d_2\in\mathbb Z$ we have
        \begin{equation*}
            c_1+2c_2=nj=d_1+2d_2
        \end{equation*}
        then $(d_1-c_1,d_2-c_2)$ is a $\mathbb Z$ multiple of $(2,-1)$.
        
          Let now $w\in\mathcal{L}(\mathbb N)$ with $3\leq\weight(w)\leq j$ and assume $M_j\cap \mathbb Z^{\{v\in\lyndon(\mathbb N)|v< w\}}$ generates $L_j\cap \mathbb Z^{\{v\in\lyndon(\mathbb N)|v< w\}}$.
    Let furthermore $y\in L_j\cap \mathbb Z^{\{v\in\lyndon(\mathbb N)|v\leq w\}}$ arbitrary and let $n\in\mathbb Z$ be the coefficient of $w$ in $y$.
    The integer vector $x=(j-\weight(w))\cdot \letter 1+w$ has weight $j$ and thus $x\in M_j$.
    Therefore
    \begin{equation*}
        y=nx+(y-nx)
    \end{equation*}
    is in the lattice generated by $M_j\cap \mathbb Z^{\{v\in\lyndon(\mathbb N)|v\leq w\}}$ as $y-nx\in L_j\cap \mathbb Z^{\{v\in\lyndon(\mathbb N)|v<w\}}$.

     $L_j$ has index $j$ as by definition, $L_j$ is the kernel of the morphism
        \[\mathbb Z^{\lyndon(\mathbb N)_{\leq j}} \to \mathbb Z / j\mathbb Z, \ w_i \mapsto \mathrm{wt}(w_i)\!\!\!\!\mod j\]
    which is surjective for $j\geq 1$.

    \item This can be argued as in Proposition \ref{prop:rel0-is-rel}, using that now there are $2$ distinct monomials of weight $2$.
    \item The first assertion follows from 1.\ and Theorem~\ref{thm:toric-recall}.
    Note that the induction argument used to prove point 1.\ shows also that $M_j\cap \mathbb Z^{\lyndon(\mathbb N)_{\leq j-2}}$ generates $L_j\cap \mathbb Z^{\lyndon(\mathbb N)_{\leq j-2}}$,
    which then implies the second assertion.
    \item Analogous to the proof of Theorem \ref{thm:mot-non-normal}. We just replace $\lyndon(\letter 2,\letter 3)_{\leq j}$ by $\lyndon(\mathbb N)_{\leq j}$, $\monlyndon(\letter 2,\letter 3)_{\leq j}$ by $\monlyndon(\mathbb N)_{\leq j}$, $A_{\ZmC,j}$ by $A_{\StC,j}$ and $\letter{23}$ by $\letter 5$.
    \item Completely analogous to the proof of Proposition \ref{prop:not_gen_in_fixed_deg}. However, since $\lambda_{\mathbb N,j}\geq 1$ for all $j\geq 1$,  the entries of the rigid matrices may now be arbitrary positive integers.
    \end{enumerate}
\end{proof}

\begin{figure}[h]
    \centering
    \includegraphics[width=0.4\linewidth]{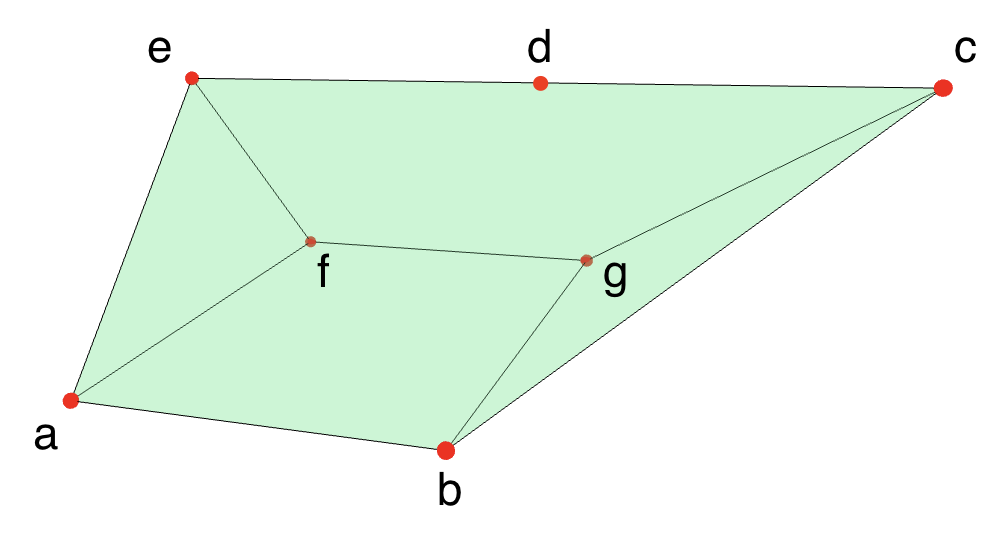}
    \caption{The three-dimensional polytope associated to the variety $\mathcal R_{5}^0(\StC)$.}\label{fig:poly-st-d3}
\end{figure}

\begin{table}
    \centering
    \renewcommand{\arraystretch}{1.2}
    \begin{tabular}{c|ccccccccc}
        Property $\backslash$ $j$ & 4 & 5 & 6 & 7 & 8 & 9 & 10 & 11 & 12\\\hline\vspace{-1em} \\ 
        $\mathrm{dim} \ \relvar_j(\StC)$ & 6 & 12 & 21 & 39 & 69 & 125 & 224 & 410 & 745\\
        amb.\ $\dim \ \relvar_j(\StC)$ & 7 & 15 & 31 & 63 & 127 & 255 & 511 & 1023 & 2047 \!\! \\
        $\mathrm{dim} \ \relvarr_j(\StC)$ & 1 & 3 & 6 & 12 & 21 & 39 & 69 & 125 & 224\\
        amb.\ $\dim \relvarr_j(\StC)$ & 2 & 6 & 16 & 36 & 79 & 169 & 356 & 738 & 1526\\
        $\mathrm{deg} \ \relvar_j(\StC)$ & 2 & 4 & 24 & 160 & 4456 & 314590 & \!\!\! 226690112 \!\!\! \\
        $\!\!\!\!\relvar_j(\StC)$ is proj.\ normal & \true & \true & \true & \true & \true\\
        $\relvar_j(\StC)$ is normal & \true & \true & \true & \true & \true \\
        $\!\!\!\!\relvarr_j(\StC)$ is proj.\ normal & \true & \true & \true & \true & \true & \true & \true  & \true & \true\\
        $\relvarr_j(\StC)$ is normal & \true & \true & \true & \true & \true & \true & \true  & \true & \true\\
        $\relvarr_j(\StC)$ is smooth & \true & \true & \false & \false & \false & \false & \false \\
        generators of $I_{\StC, j}$ &
        $1_2$ & $6_2$ & $42_2$ & $211_2$ & $1060_2$ & $4946_2$ & $22342_2$ &
    \end{tabular}

    \caption{Properties of the projective varieties $\relvar_j(\StC)$ and $\relvarr_j(\StC)$.
    Note that $\relvarr_j(\StC)$ is also projectively normal (and thus normal) for $j=13,14,15,17$. The empty cells in the table represent calculations which did not terminate in reasonable time due to computational limits.}
    \label{tab:stuffle}.
\end{table}

\section{Proof of Lemma \ref{lem:rhoisom}}
\label{app:proofrho}

We show that the inverse of
\begin{align*}
\rho: S^{\stuffle} \to S^{\shuffle}, \quad 
w \mapsto \sum_{w =v \letter 1^k} \lambda_k \shuffle \gamma(v)
\end{align*}
is given by 
\begin{align*}
\rho^{-1}: S^{\shuffle} \to S^{\stuffle},\quad
w\mapsto \sum_{w=v \letter 1^k} \gamma^{-1}(\mu_k\shuffle v),
\end{align*}
where the $\lambda_k$ are defined in \eqref{eq:def_lambda} and the $\mu_k$ in Lemma \ref{lem:rhoisom}

\begin{proof} 

By construction, we have
$\lambda_0=\mu_0=1$ and
\begin{align} \label{eq:lambda_mu}
\sum_{k=0}^n \lambda_k\shuffle \mu_{n-k}=0,\qquad n\geq1.
\end{align}
Let $w=v \letter 1 ^n$ where $v$ is a word in $\letter 0, \letter 1$ which does not end in $\letter 1$, then we have
\begin{align*}
\rho^{-1}(\rho(w))&=\sum_{k=0}^n \rho^{-1}(\lambda_k \shuffle \gamma(v)\letter 1^{n-k}) \\
&=\rho^{-1}(\gamma(v)\letter 1^n)+\sum_{k=2}^n \Big(\sum_{l=0}^{n-k-1}\rho^{-1}\big(\lambda_k'\shuffle \gamma(v) \letter 1^{n-k-l})\letter 0\letter 1^l + (\lambda_k\shuffle \gamma(v))\letter 1^{n-k}\big)\Big),
\end{align*}
as $\lambda_0=1$, $\lambda_1=0$, and $\lambda_k=\lambda_k'\letter 0$ for $k\geq 2$. Thus, we get
\begin{align*}
&\rho^{-1}(\rho(w))=\gamma^{-1}(\gamma(v)\letter 1^n)+\sum_{m=2}^n \gamma^{-1}(\mu_m\shuffle \gamma(v)\letter 1^{n-m}) \\
&+\sum_{k=2}^n\Big(\sum_{l=0}^{n-k-1}\Big(\gamma^{-1}((\lambda_k'\shuffle \gamma(v)\letter 1^{n-k-l})\letter 0 \letter 1^l)+\sum_{m=2}^l \gamma^{-1}(\mu_m\shuffle (\lambda_k' \shuffle \gamma(v)\letter 1^{n-k-l})\letter 0 \letter 1^{l-m} ))\Big) \\
&+\gamma^{-1}((\lambda_k\shuffle\gamma(v))\letter 1^{n-k})+\sum_{m=2}^{n-k}\gamma^{-1}(\mu_m\shuffle ((\lambda_k\shuffle \gamma(v))\letter 1^{n-k-m}))\Big),
\end{align*}
as again $\mu_0=1$, $\mu_1=0$, and $\mu_m$ ends in $\letter 0$ for $m\geq2$. In particular, also every word in $\lambda_k\shuffle \mu_m$ ends in $\letter 0$ for $k,m\geq2$. Hence, we deduce
\begin{align*}
\rho^{-1}(\rho(w))&= w+\sum_{m=2}^n \gamma^{-1}(\mu_m\shuffle \gamma(v)\letter 1^{n-m})+\sum_{k=2}^n \gamma^{-1}(\lambda_k\shuffle \gamma(v)\letter 1^{n-k}) \\
&\hspace{0.4cm}+\sum_{k=2}^n\sum_{m=2}^{n-k} \gamma^{-1}((\mu_m\shuffle\lambda_k)\shuffle \gamma(v)\letter 1^{n-k-m})\\
&=\sum_{k=0}^n\sum_{m=0}^{n-k} \gamma^{-1}((\mu_m\shuffle\lambda_k)\shuffle \gamma(v)\letter 1^{n-k-m}) \\
&=\sum_{l=0}^n \sum_{m+k=l} \gamma^{-1}((\mu_m\shuffle \lambda_k)\shuffle \gamma(v)\letter 1^{n-l}).
\end{align*}
By \eqref{eq:lambda_mu} the terms for $l>0$ vanish and hence
\begin{align*}
\rho^{-1}(\rho(w))= \gamma^{-1}((\mu_0\shuffle \lambda_0)\shuffle \gamma(v) \letter 1^n)=\gamma^{-1}(\gamma(v)\letter 1^n)=v\letter 1^n=w.
\end{align*}
A similar argument shows that $\rho(\rho^{-1}(w))=w$.
\end{proof}

\section{Numerical evaluation of MZVs in \macaulay}\label{app:macaulayzeta}

We indicate a straight-forward approach to numerical evaluation of MZVs to high precision, and give an implementation in \macaulay.  The method is due to D. Zagier, as explained in his talk \cite[11m50--22m]{zagierMZVtalk}.  Define the $M$-th order tail of an MZV by
\[
    \zeta_{\gg M}(s_1,\ldots,s_r) \coloneqq \sum_{M < k_1 < \cdots < k_r} \frac{1}{k_1^{s_1} \cdots k_r^{s_r}} \,, \quad \zeta_{\gg M}(\emptyset) \coloneqq 1 \,.
\]
The goal is to accurately compute \( \zeta(s_1,\ldots,s_r)  = \zeta_{\gg 0}(s_1,\ldots,s_r) \).  This is done in two steps: fix \( M > 0 \), then i) use Euler-Maclaurin summation to estimate the asymptotic tail \( \zeta_{\gg M}(s_1,\ldots,s_r) \), ii) explicitly sum the remaining \( M \) initial terms recursively. \medskip

\paragraph*{Summation of initial terms:} In \( \zeta_{\gg i-1}(s_1,\ldots,s_r) \), the condition \( k_1 > i-1 \) implies either \( k_1 > i \) or \( k_1 = i \), so following recursion holds:
\begin{equation}\label{eqn:mzvtail}
    \zeta_{\gg i-1}(s_1,\ldots,s_r) = \zeta_{\gg i}(s_1,\ldots,s_r) + \frac{1}{i^{s_1}} \zeta_{\gg i}(s_2,\ldots,s_{r}) \tag{$\ast_i$} \,.
\end{equation}
If we have (good) estimates for \( \zeta_{\gg M}(s_1,\ldots,s_j) \), with \( 1 \leq j \leq r \), then applying \eqref{eqn:mzvtail}, for \( i=M, M-1, \ldots, 1 \) and \( 1 \leq j \leq r \) calculates a (good) estimate \( \zeta_{\gg 0}(s_1,\ldots,s_r) = \zeta(s_1,\ldots,s_r) \).

To do this, first store the initial estimates 
\[
 [\zeta_{\gg M}(s_1,\ldots,s_r), \zeta_{\gg M}(s_2,\ldots,s_{r}), \ldots, \zeta_{\gg M}(s_r) , {\zeta_{\gg M}(\emptyset)} \coloneqq 1] \,.
\]
Then for \( 1 \leq j \leq r \), compute \(  \zeta_{\gg i-1}(s_j,\ldots,s_r) \) from \(  \zeta_{\gg i}(s_j, \ldots, s_r) \) and \(\zeta_{\gg i}(s_{j+1}, \ldots, s_{r}) \) using \eqref{eqn:mzvtail}, and store the  new estimates.  (Since \( \zeta_{\gg i-1}(s_j,\ldots,s_r) \) only requires  \( \zeta_{\gg i}(s_k, \ldots, s_r) \), \( k \geq j \), the already stored values can be updated in place.)  Repeat for \( M \geq i \geq 1 \), until \( \zeta_{\gg 0} \) is obtained.  This computes \( \zeta_{\gg 0} \) in time \( O(Mr) \).

\paragraph*{Estimate of asymptotic tail:} To estimate \( \zeta_{\gg M}(s_1,\ldots,s_r) \) we use Euler--Maclaurin summation recursively to find an asymptotic expansion for \( \zeta_{\gg M} \) as a power series in \( \frac{1}{M} \).

Euler-Maclaurin summation gives an expression for the difference between the integral and the sum of a given (\( 2N \)-times continuously differentiable) function \( f \):
\begin{align*}
    \sum_{n={a+1}}^b f(n) = {} & \int_a^b f(x) \mathrm{d}x + \frac{f(b) - f(a)}{2} + \sum_{k=1}^{N} \frac{B_{2k}}{(2k)!} \big( f^{(2k-1)}(b) - f^{(2k-1)}(a) \big) \\
    & {} - \frac{1}{(2N)!} \int_a^b B_{2N}(x - \lfloor x \rfloor) f^{(2N)}(x) \mathrm{d}x \,.
\end{align*}
Here \( B_{2N} \) is the $2N$-th Bernoulli polynomial, and the last term is an error term.  For  \( f(x) = \frac{1}{x^s} \), with \( s \geq 2 \), one has \( f^{(k)}(x) = (-1)^k k! \binom{s+k-1}{k} \frac{1}{x^{s+k}} \).  For \( a = M \) and \( b \to \infty \), we obtain
\begin{equation}\label{eqn:eulermac}
\begin{aligned}
    \sum_{n=M+1}^\infty \frac{1}{n^s} = {} & \frac{1}{(s-1)M^{s-1}} - \frac{1}{2M^s} + \sum_{k=1}^{N} \frac{B_{2k}}{2k} \binom{s+2k-2}{2k-1} \frac{1}{M^{s+2k-1}} \\
    & {} - \binom{s+2N-1}{2N} \int_{M}^\infty B_{2N}(x - \lfloor x \rfloor) x^{-s-2N} \mathrm{d}x
\end{aligned}
\end{equation}
Note \( B_{2N}(x - \lfloor x \rfloor) \) is bounded (lying in the image of \( [0,1] \) under the polynomial \( B_{2N}(x) \)).  Precise maxima and minima are known \cite{Lehmer}, but  \( -\lvert B_{2N} \rvert \leq B_{2N}(x - \lfloor x \rfloor) \leq \lvert B_{2N} \rvert \) is sufficient.  The error term above then satisfies
\begin{align*}
    \Bigg\lvert \binom{s+2N-1}{2N} \int_{M}^\infty B_{2N}(x - \lfloor x \rfloor) x^{-s-2N} \mathrm{d}x \Bigg\rvert & 
    \leq \binom{s+2N-1}{2N} \, \lvert B_{2N} \rvert \int_{M}^\infty x^{-s-2N} \mathrm{d}x \\
    & = \binom{s+2N-1}{2N} \frac{\lvert B_{2N} \rvert}{(s + 2N -1) M^{s + 2N - 1}} \,,
\end{align*}
which rapidly tends to 0 as \( M \to \infty \) ($N$ fixed). 

To evaluate \( \zeta(s) \) from this setup, fix \( N > 0 \), and choose \( M \) such that
\[
    \binom{s+2N-1}{2N} \frac{\lvert B_{2N} \rvert}{(s + 2N -1) M^{s + 2N - 1}} < \varepsilon \,.
\]  Then the main terms of \eqref{eqn:eulermac} gives an estimate for 
\[
    \zeta_{\gg M}(s) = \sum_{n=M+1}^\infty \frac{1}{n^s} \approx \frac{1}{(s-1)M^{s-1}} - \frac{1}{2M^s} + \sum_{k=1}^{N} \frac{B_{2k}}{2k} \binom{s+2k-2}{2k-1} \frac{1}{M^{s+2k-1}}
\]
to an error of at most \( \varepsilon \).  Directly summing the initial \( M \) remaining terms gives \( \zeta(s) \)  to an error of at most \( \varepsilon \).  More precisely: if \( M = 2N \), then the error is decreasing in \( s \) (the ratio of consecutive terms is \( \frac{2N + s-1}{Ms} < 1 \) ).  Setting \( s = 2 \), \( M = 2N \), and using the  evaluation \( \lvert B_{2n} \rvert = \frac{2(2N)!}{(2\pi)^{2N}} \zeta(2N) < \frac{4(2N)^{2N}}{(2\pi)^{2N}}   \), shows the error is \( < \frac{4}{N \cdot (2\pi)^{2N}} \) in this case, say. \medskip

Recursively applying Euler-Maclaurin summation to the polynomial part of \eqref{eqn:eulermac}, gives similar rapidly convergent series for the tails of MZVs.  Suppose it is known that
\[
    \sum_{M < n_1 < \cdots < n_d} \frac{1}{n_1^{s_1} \cdots n_d^{s_d}} = A_{s_1,\ldots,s_d}(M^{-1}) + O(M^{-\ell}) \,.
\]
for some polynomial \( A_{s_1,\ldots,s_d}(x) = \sum_{i=0}^{\ell-1} c_i x^i \).  Then with \( s_0 \geq 1 \),
\[
    \sum_{n_0 = M+1}^{\infty} \frac{1}{n_0^{s_0}} \sum_{n_0 < n_1 < \cdots < n_d}  \frac{1}{n_1^{s_1} \cdots n_d^{s_d}} = \sum_{n_0 = M+1}^{\infty} \frac{1}{n_0^{s_0}} A(n_0^{-1}) + \sum_{n_0 = a+1}^\infty \frac{1}{n_0^{s_0}} O(n_0^{-\ell}) \,.
\]
Application of Euler-Maclaurin summation to each term \( \frac{c_i}{n_0^{s_0+i}} \) produces another polynomial in \( M^{-1} \) and a error of \( O(M^{-j}) \) for some \( j \geq \ell \).  Finally, note
\[
    \sum_{n_0 = M+1}^\infty \frac{1}{n_0^{s_0}} O(n_0^{-\ell}) \leq \int_M^\infty O(n_0^{-\ell-s_0}) = O(M^{-\ell-s_0+1}) = O(M^{-\ell}) \,.
\]
Overall, we have
\[
    \sum_{M < n_0 < n_1 < \cdots < n_d} \frac{1}{n_0^{s_0} n_1^{s_1} \cdots n_d^{s_2}} = A_{s_0,s_1,\ldots,s_d}(M^{-1}) + O(M^{-\ell}) \,,
\]
for some new polynomial \( A_{s_0,s_1,\ldots,s_d}(x) \) of degree \( \ell-1 \).

To compute the polynomial effectively, one can use the recursion \eqref{eqn:mzvtail} to derive the following functional relation for \( A_{s_1,\ldots,s_d}(x) \),
\[
    A_{s_1,\ldots,s_d}\Big( \frac{x}{1-x} \Big) = A_{s_1,\ldots,s_d}(x) + x^{s_1} A_{s_2,\ldots,s_d}(x) \,.
\]
If we known \( A_{s_2,\ldots,s_d}(x) \) to order \( x^{\ell-1} \), then this sets up a system of equations for \( A_{s_1,\ldots,s_d} \).  Write \( F(x) = \sum_{i=0}^{\ell-1} a_i x^i \), and \( G(x) = \sum_{i=0}^{\ell-1} b_i x^i \).  Then comparing coefficients in
\[
    F\Big( \frac{x}{1-x} \Big) = F(x) + x^\alpha G(x) 
\]
gives 
\[
    \sum_{k=0}^{i-1} a_k \binom{i-1}{i-k} = b_{i-\alpha} 
    \implies a_i = \frac{1}{i} \Bigg\{ b_{i+1-\alpha} - \sum_{k=0}^{i-1} a_k \binom{i}{i+1-k} \Bigg\} \,.
\]
Since trivially \( \zeta_{\gg M}(\emptyset) = 1 \), we start with \( A_{\emptyset}(x) \equiv 1 \); the recursion will automatically compute the Bernoulli number coefficients appearing in \( A_{s_1}(x) \) as the first step.

\medskip
The \macaulay{} code available on the \github{} repository implements this numerical evaluation.  A reasonable choice for at least $\texttt{prec} = 100$ decimal digits of accuracy is typically $\texttt{ord} = 50$ for the asymptotic tail, and $\texttt{M} = 1100$, for the initial terms.

\bigskip
\hrule
\bigskip

 {\par\nobreak\vfil\penalty0\vfilneg
   \vtop\bgroup

\vspace{1cm}
\begin{small}
\noindent $\bullet$ Annika Burmester (\href{mailto:aburmester@math.uni-bielefeld.de}{aburmester@math.uni-bielefeld.de})\\
\phantom{lo}
\textit{Universit\"at Bielefeld, Bielefeld, Germany	} 
\\[0.2cm]
\noindent $\bullet$ Steven Charlton (\href{mailto:mail@stevencharlton.net}{mail@stevencharlton.net})
\\[0.2cm]
\noindent $\bullet$ Abhiram Kidambi (\href{mailto:mkabhiram@duck.com}{mkabhiram@duck.com})\\
\phantom{ }
\textit{ 
Albert Einstein Institute, Potsdam, Germany \& \\[-0.1cm] \phantom{lo}
Max Planck Institute for Mathematics in the Sciences, Leipzig, Germany} 
\\[0.2cm]
\noindent $\bullet$ Felix Lotter (\href{mailto:felix.lotter@mis.mpg.de}{felix.lotter@mis.mpg.de})\\\phantom{lo}
\textit{Max Planck Institute for Mathematics in the Sciences, Leipzig, Germany} 
\\[0.2cm]
\noindent $\bullet$ Rosa Prei\ss{} (\href{mailto:preiss@math.tu-berlin.de}{preiss@math.tu-berlin.de})\\\phantom{lo}
\textit{Technische Universit\"at Berlin, Berlin, Germany	} 
\end{small}
  \par\xdef\tpd{\the\prevdepth}\egroup
   \prevdepth=\tpd}
\end{document}